\documentclass[superscriptaddress,
 reprint,
 amssymb, amsmath,
 aip, cha 
]{revtex4-1}
\usepackage[utf8]{inputenc}
\usepackage[T1]{fontenc}
\usepackage{amsmath}
\usepackage{amsfonts}
\usepackage{amssymb}
\usepackage{graphicx}
\usepackage{bm}
\usepackage{grffile} %Allows figure file names with periods
\usepackage{dcolumn}
\usepackage{lipsum}
\usepackage{multirow}
\usepackage{tabularx}
\usepackage{booktabs}
\usepackage{tikz}
 
\usepackage{pifont}% http://ctan.org/pkg/pifont
\newcommand{\cmark}{\ding{51}}%
\newcommand{\xmark}{\ding{55}}%
\newcommand{\dmark}{\ding{79}}%
\usepackage{url}
\usepackage{listings}
\newcommand{\qj}{q^{(j)}}
\newcommand{\aj}{a^{(j)}}
\newcommand{\notilde}{}

\renewcommand*{\Xi}{{\boldsymbol{\xi}}}

\newcommand*{\Z}{{\mathbb{Z}}}

\newcommand*{\R}{{\mathbb{R}}}

\newcommand*{\Econd}{\mathbb{E}[\Phi_n|d_k;\,  k<n]}

\DeclareMathOperator{\sign}{sign}

\newcommand{\eps}{\varepsilon}
\renewcommand{\j}{^{(j)}}
\renewcommand{\O}{\mathcal{O}}
\newcommand{\E}{\mathbb{E}}
\DeclareMathOperator{\cov}{cov}

\graphicspath{{figures/}}

\newcommand{\Rone}[1]{#1}
\begin{document}

%\maketitle

\title{Linear response for macroscopic observables in high-dimensional systems}%arge-scale dynamics and linear response of weakly-coupled chaotic ensembles}

\author{Caroline L. Wormell}
 \email{ca.wormell@gmail.com}
 \affiliation{\mbox{School of Mathematics and Statistics, The University of Sydney, Sydney, NSW 2006, Australia}}
\author{Georg~A. Gottwald}
 \email{georg.gottwald@sydney.edu.au}
 \affiliation{\mbox{School of Mathematics and Statistics, The University of Sydney, Sydney, NSW 2006, Australia}}

\date{\today}

%%%%%%%%%%%%%%%%%%%%%%%%%%%%%%%%%%%%%%%%%%%%%%%%

\begin{abstract}
The long-term average response of observables of chaotic systems to dynamical perturbations can often be predicted using linear response theory, but not all chaotic systems possess a linear response. Macroscopic observables of complex dissipative chaotic systems, however, are widely assumed to have a linear response even if the microscopic variables do not, but the mechanism for this is not well-understood.

We present a comprehensive picture for the linear response of macroscopic observables in high-dimensional coupled deterministic dynamical systems, where the coupling is via a mean field and the microscopic subsystems may or may not obey linear response theory. We derive stochastic reductions of the dynamics of these observables from statistics of the microscopic system, and provide conditions for linear response theory to hold in finite dimensional systems and in the thermodynamic limit. In particular, we show that for large systems of finite size, linear response is induced via self-generated noise.

We present examples in the thermodynamic limit where the macroscopic observable satisfies LRT, although the microscopic subsystems individually violate LRT, as well a converse example where the macroscopic observable does not satisfy LRT despite all microscopic subsystems satisfying LRT when uncoupled. This latter, maybe surprising, example is associated with emergent non-trivial dynamics of the macroscopic observable. We provide numerical evidence for our results on linear response as well as some analytical intuition.
\end{abstract}

\maketitle

%%%%%%%%%%%%%%%%%%%%%%%%%%%%%%%%%%%%%%%%%%%%%%%%

\begin{quotation}
Since its introduction mid-last century, linear response theory (LRT) has been a cornerstone of statistical mechanics. If a system has a linear response, one can estimate the change of expectation values caused by a perturbation of a parameter using only information of the unperturbed system. LRT has been successfully applied in numerous areas, ranging from neurophysiology to climate science. It is widely believed that high-dimensional complex dynamical systems satisfy LRT, largely based on successes in applications. Separate efforts by mathematicians to understand the dynamical underpinnings of linear response theory, however, have found that many low-dimensional systems such as the logistic map do not obey LRT, but instead exhibit a rough dependency of their statistical properties with respect to perturbations. In this work we investigate the conditions for and mechanisms by which a system comprised of low-dimensional subunits---which individually may or may not obey LRT---may at a macroscopic scale have a linear response.
\end{quotation}

%%%%%%%%%%%%%%%%%%%%%%%%%%%%%%%%%%%%%%%%%%%%%%%%

\section{Introduction}

Since its introduction in the 1960s, linear response theory (LRT) has been widely used across numerous disciplines to quantify the change of the mean behaviour of observables in a perturbed environment. LRT is valid, in essence, provided the invariant measure varies differentiably with respect to the perturbation; consequently  LRT allows for a Taylor expansion of the perturbed invariant measure around the unperturbed invariant measure. Hence, when valid, LRT provides an expression of the average of some observable when subjected to small perturbations from an unperturbed state -- the system's so called {\em{response}} -- entirely in terms of statistical information from the unperturbed system. \\

Climate scientists in particular have successfully applied LRT to eke out valuable information about the change of certain atmospheric and oceanic observables under changed climatic conditions. Applications include atmospheric toy models \cite{MajdaEtAl10,LucariniSarno11,AbramovMajda07,AbramovMajda08,CooperHaynes11,CooperEtAl13}, barotropic models \cite{Bell80,GritsunDymnikov99,AbramovMajda09}, quasi-geostrophic models \cite{DymnikovGritsun01}, atmospheric models \cite{NorthEtAl93,CionniEtAl04,GritsunEtAl02,GritsunBranstator07,GritsunEtAl08,RingPlumb08,Gritsun10} and coupled climate models \cite{LangenAlexeev05,KirkDavidoff09,FuchsEtAl14,RagoneEtAl15}. The seminal work by Ruelle \cite{Ruelle97,Ruelle98,Ruelle09a,Ruelle09b} rigorously established that LRT is valid in uniformly hyperbolic Axiom A systems. Success in reliably estimating the response of a physical system, as exemplified by the above applications in the climate sciences, prompted scientists to believe that general  chaotic dynamical systems obeyed LRT. This belief was proven wrong by Baladi and co-workers \cite{BaladiSmania08,BaladiSmania10,Baladi14,BaladiEtAl15,DeLimaSmania18} who showed that simple dynamical systems such as the logistic map violate LRT and support an invariant measure that changes non-smoothly with respect to the perturbation. This raises the question of how a high-dimensional dynamical system, despite its constituent subsystems typically individually violating LRT, may exhibit linear response.\\

The majority of the scientific community, including the authors, believe that the interaction between the microscopic constituents in typical high-dimensional systems leads to an emergence of LRT at the macroscopic level. How exactly this is achieved and what the conditions are for the dynamical systems for which LRT is guaranteed, however, remains an open question. In the literature the validity of LRT in high-dimensional deterministic systems is often justified by appealing to the {\em{chaotic hypothesis}} of Gallavotti and Cohen \cite{GallavottiCohen95a,GallavottiCohen95b, Gallavotti19} according to which the attracting dynamics of high-dimensional system behaves for all practical purposes as an Anosov system. However, even under this hypothesis one cannot relate the equivalent Anosov systems for different perturbations, which is the focus of LRT. In particular, for dissipative systems the response of the attracting dynamics to perturbations depends on the properties of the flow outside the attractor as well as on it: off the attractor the flow may be non-hyperbolic, and hence leading to a breakdown of linear response. In stochastic systems, however, it is well established that LRT can be justified \cite{Haenggi78,HairerMajda10}.

In previous work, we argued that a combination of statistical limits of the high-dimensional system and a sufficient degree of heterogeneity in the system causes a high-dimensional system to obey LRT, even when the individual microscopic subsystems do not obey LRT   \cite{WormellGottwald18}. We considered a single resolved degree of freedom weakly coupled to $M$ unresolved uncoupled degrees of freedom, the so called heat bath, which evolve according to their own randomly drawn parameters. Both the distinguished degree of freedom as well as the heat bath were described by logistic maps, which individually violate LRT. In the thermodynamic limit we derived a stochastic limit system for the distinguished degree of freedom; it was shown, however, that the mere presence of stochasticity is not sufficient to guarantee LRT, but the microscopic subsystems need to be appropriately heterogeneous, with the parameters of the logistic map drawn from a sufficiently smooth distribution. The perturbations considered were homogeneous perturbations of the randomly-drawn logistic parameters in the microscopic heat bath system as well as general smooth perturbations in the evolution of the distinguished macroscopic variable.\\

%In \cite{GottwaldEtAl16} it was demonstrated that, given a time series of a system which provably violates LRT, linear response might be consistent with the underlying data and the breakdown of LRT may only be detectable with sufficient statistical significance for extremely long time series. Hence, it can not be excluded that the apparent observed validity of LRT in the climate sciences might in fact be a finite size effect. This lead to the construction of weakly coupled deterministic high-dimensional systems which obeyed LRT in the thermodynamic limit \cite{WormellGottwald18}. We considered a single resolved degree of freedom weakly coupled to $M$ unresolved uncoupled degrees of freedom, the so called heat bath. Both, the distinguished degree of freedom as well as the heat bath were described by logistic maps, violating individually LRT. The perturbations we considered were homogeneous perturbations of the logistic map parameters in the microscopic heat bath system as well as a perturbation in the logistic map parameter of the distinguished macroscopic variable. In our recent work we showed that LRT can be assured in such high-dimensional systems of coupling type, when the coupling is such that the macroscopic resolved variables exhibit effective stochastic dynamics and provided that the heat bath-like system is appropriately heterogeneous. \\

We continue this line of research and consider here macroscopic observables of high-dimensional dynamical systems whose microscopic constituents may violate LRT, rather than observables only of individual distinguished degrees of freedom. We extend our previous work to consider more general perturbations than homogeneously perturbing the parameters, and include the more realistic case when the microscopic dynamics are globally coupled via a mean field. This latter case has been well-studied, particularly in the case where the coupling is strictly attractive, and complex emergent dynamics at the level of the mean-field have been observed\cite{Kaneko90,Shibata99,Pikovsky94,Ershov95,Ershov97,Selley16}. We shall provide a systematic macroscopic reduction for the mean-field coupled dynamics, which we use to study a range of interesting dynamical scenarios in the context of linear response. We provide a comprehensive picture of the linear response behaviour of macroscopic observables, for uncoupled and for mean field coupled systems, and we find that the existence of LRT depends in an intricate way on the combination of effective stochastic behaviour of the macroscopic observable, the macroscopic dynamics of the thermodynamic limit, and on the smoothing property of heterogeneously distributed dynamical parameters of the microscopic subsystems. Indeed, we will present a case when all individual microscopic subsystems obey LRT when uncoupled, but the collective macroscopic dynamics violates LRT, and cases where this holds {\it vice versa}. To corroborate our findings, we will use a recently developed statistical test which allows to probe for the validity of LRT in a given time series \cite{GottwaldEtAl16}.\\

The paper is organized as follows. Section~\ref{s.LRT} briefly reviews LRT. We introduce the high-dimensional systems under consideration in Section~\ref{s.model} and summarise our results in Section \ref{s.results}. Sections~\ref{s.uncoupled} and~\ref{s.coupled} provide numerical evidence and an analytical treatment corroborating the results summarized in Table~\ref{t.result} for uncoupled and for mean field coupled systems. We conclude with a discussion and an outlook in Section~\ref{s.discussion}.

%%%%%%%%%%%%%%%%%%%%%%%%%%%%%%%%%%%%%%%%%%%%%%%%%%%%%%%%%%

\section{Linear response theory}
\label{s.LRT}
We briefly review some basic notation of linear response theory. Consider a family of dynamical systems $f_\varepsilon:D \to D$ on some space $D$ where the map $f_\varepsilon$ depends smoothly on the parameter $\varepsilon$ and where for each $\varepsilon$ the dynamical system admits a unique invariant physical measure $\mu_\varepsilon$. An ergodic measure is called physical if for a set of initial conditions of nonzero Lebesgue measure the temporal average of a typical observable converges to the spatial average over this measure. LRT is concerned with the change of the average of an observable $\phi:D\to\R$,  
\begin{align*}
\E^\varepsilon [\phi] = \int_D \phi\, d\mu_\varepsilon
\end{align*}
upon varying $\varepsilon$. A system exhibits {\em{linear response}} at $\eps=\eps_0$,  if the derivative
\begin{align*}
\E^{\varepsilon_0} [\phi]^\prime := \frac{\partial}{\partial\varepsilon} \E^\varepsilon [\phi]_{|_{\varepsilon_0}}
% = \lim_{\varepsilon \rightarrow \varepsilon_0} 
%\frac{\langle A\rangle_\varepsilon - \langle A\rangle_{\varepsilon_0}}{\varepsilon - \varepsilon_0}
\end{align*}
exists. A sufficient condition for this is that the invariant measure $\mu_\varepsilon$ is differentiable with respect to $\varepsilon$. This derivative can be expressed entirely in terms of the invariant measure $\mu_{\eps_0}$ of the unperturbed system using so-called linear response formulae \cite{Ruelle09a,Ruelle98,Baladi14}. The average of an observable of the perturbed state is then expressed to first order as
\begin{align*}
\E^\varepsilon [\phi]  \approx \E^{\varepsilon_0} [\phi] +(\eps-\eps_0)\, \E^{\varepsilon_0} [\phi]^\prime .
\end{align*}
If the derivative exists, then this expansion expresses the remarkable result that the average of the perturbed state is determined up to $o(\eps-\eps_0)$ by the properties of the unperturbed system.
%\begin{align*}
%\E^{\varepsilon_0} [\phi]^\prime  = - \sum_{n=0}^\infty \E^{\varepsilon_0} \left[ \divg_{\rho} \left. \frac{\partial f_\eps}{\partial \eps} \right|_{\eps = 0}\ \phi\circ f_\eps^n \right].
%%{\mathrm{d}x})}{\frac{\mathrm{d} \rho}
%\end{align*}
If however it does not exist, we say there is a breakdown of linear response, which manifests itself in a rough dependency of averages of the observable on the perturbation $\eps$ \cite{GottwaldEtAl16}. 
%We remark that the choice of observables is crucial for the applicability of LRT: if the observable does not capture sufficient dynamic information about the dynamical system, linear response formulae are not informative; for example, an odd observable on a system symmetric about $0$ is by construction identically zero regardless of whether the system exhibits linear response or not.

%%%%%%%%%%%%%%%%%%%%%%%%%%%%%%%%%%%%%%%%%%%%%%%%%%%%%%%%%

\section{Model}
\label{s.model}
%
% FIGURE
%
\begin{figure*}[htbp]
\centering
\includegraphics[width=0.9\textwidth]{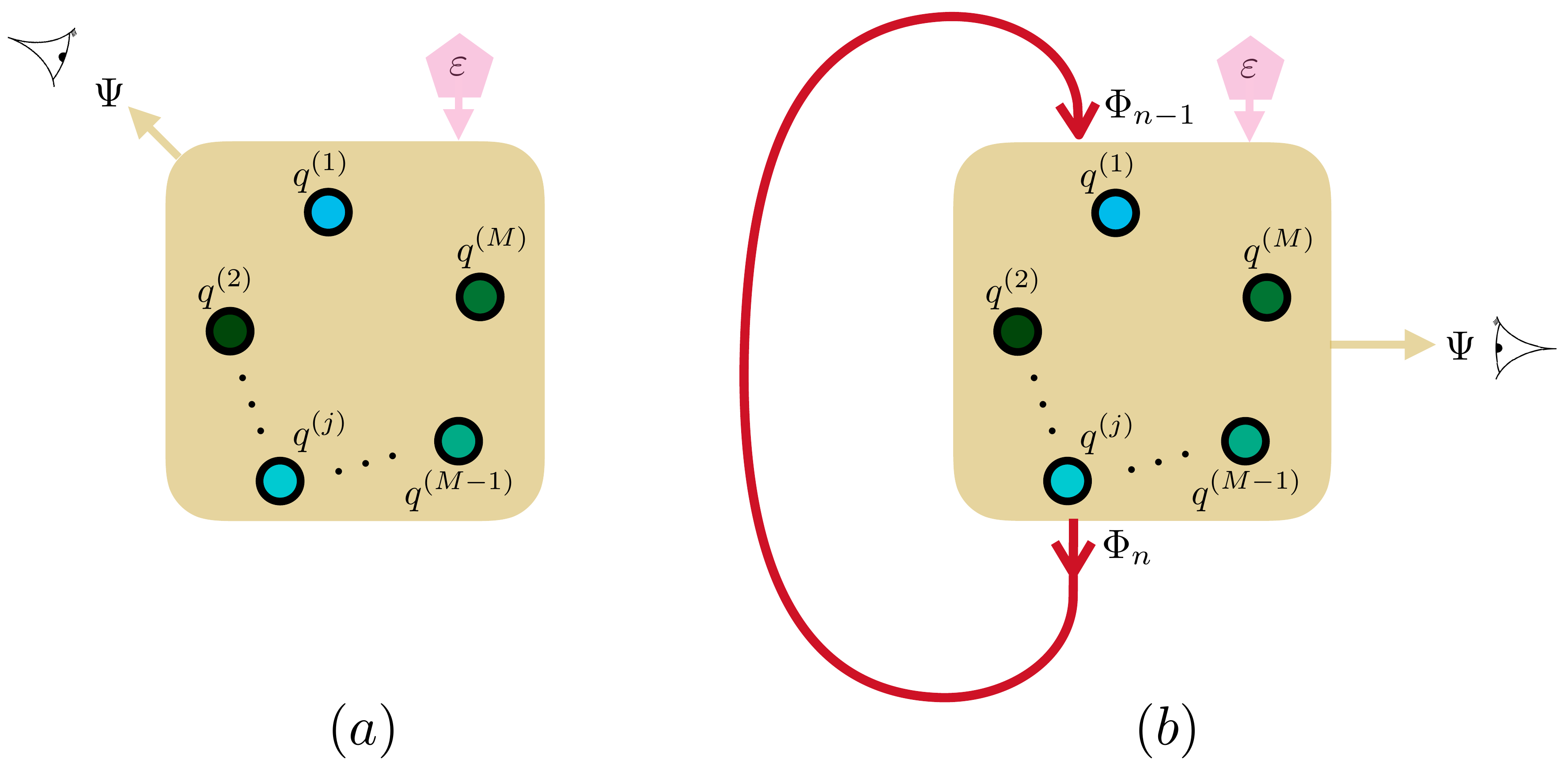}
%\centerline{(a) $\qquad\qquad\qquad\qquad\qquad\qquad\qquad\qquad$ (b)}
%(a) $\qquad\qquad\qquad\qquad\qquad\qquad\qquad\qquad$ (b)
\caption{General set-up. We consider the behaviour of macroscopic observables $\Psi$ which are constructed either  from (a): a large system of $M$ uncoupled microscopic units $q^{(j)}$, or (b) a large system in which the $M$ microscopic units are coupled via a mean field variable $\Phi$. Perturbations $\eps$ are applied globally to the dynamics of all heat bath variables $\qj$.}
\label{f.scenarios}
\end{figure*}

We consider high-dimensional systems composed of $M\gg 1$ chaotic microscopic degrees of freedom $q^{(j)}$, $j=1,\cdots,M$, which evolve in discrete time $n$ according to their individual parameters $a^{(j)}$. These degrees of freedom, evolving in isolation, may or may not obey LRT. We restrict our study of LRT to macroscopic observables
\begin{align}
\Psi_n = \Psi(q^{(1)}_n,q^{(2)}_n,\cdots,q^{(M)}_n).
\label{e.Psi0}
\end{align}
We consider in particular a mean field observable
\begin{align}
\Psi_n = \frac{1}{M} \sum_{j=1}^M \psi(q^{(j)}_n),
\label{e.Psi}
\end{align}
where $\psi$ is some observable of the microscopic variable. \Rone{We consider here smooth observables (at least H\"older continuous), which ensures that the statistical limit laws we later invoke exist.} We consider here two scenarios, illustrated in Figure~\ref{f.scenarios}, where the perturbations of size $\eps$ are globally applied to the dynamics: the case where the dynamics of the $\qj_n$ are uncoupled with
\begin{align}
\qj_{n+1}=f(\qj_n;\aj,\eps),
\label{e.scen1}
\end{align}
and the case where the dynamics of the $\qj_n$ are coupled by a mean field $\Phi$ with
\begin{align}
\qj_{n+1}=f(\qj_n,\Phi_n;\aj,\eps),
\label{e.scen2}
\end{align}
where the coupling mean field is given by
\begin{align}
\Phi_n = \frac{1}{M} \sum_{j=1}^M \phi(\qj_n)
\label{e.Phi}
\end{align}
for some function $\phi$. Note that we explicitly distinguish between the macroscopic observable $\Psi_n$ and the dynamic variable $\Phi_n$. In our numerical simulations we did not find any difference in their respective LRT properties; the distinction, however, is instructive for the theoretical considerations provided later.
\\

\begin{figure}[tbp]
	\centering
	\includegraphics[width=0.77\columnwidth]{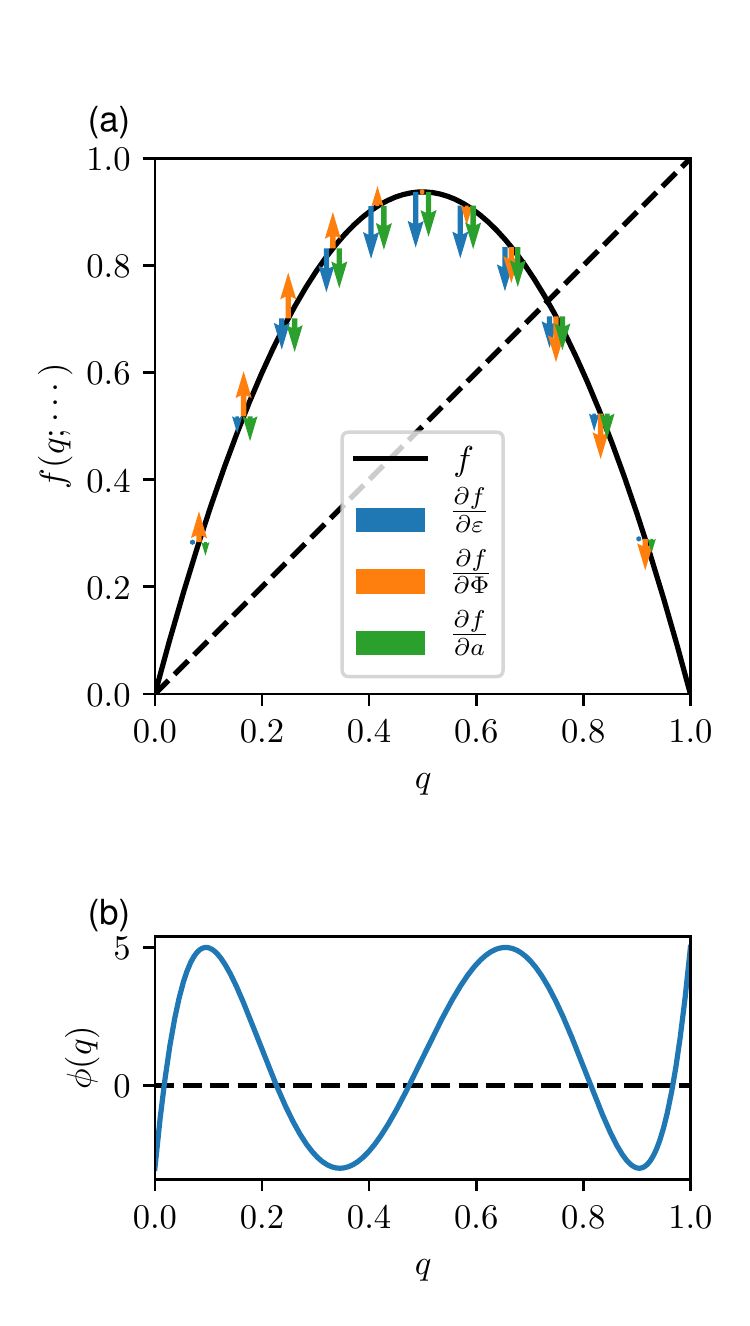}
	\caption{(a) Plot of $q\j_{n+1}$ under logistic dynamics  (\ref{e.logistic}) as a function of $q\j_n$ for $r\j_n > 1/2,\, a\j = 3.75,\, \Phi_n = 0$ and $\eps = 0$. Effects of perturbations in $\eps$, $\Phi$ and $a$ on the map are indicated by arrows. (b) Plot of coupling function $\phi(q)$.}
	\label{f.cestnepasunmaplogistique}
\end{figure}

For each of the two scenarios we study the linear response of $\Psi$ for three different types of the microscopic dynamics $f$. We consider the case of microscopic dynamics which when viewed in isolation obeys LRT, such as uniformly expanding maps (the specific maps we will study are described in Section \ref{s.coupled_TrivA}). We then consider the case when the microscopic dynamics when viewed in isolation does not obey LRT. The simplest such system is the logistic map as established by Baladi and co-workers \cite{BaladiSmania08,BaladiSmania10,Baladi14,BaladiEtAl15,DeLimaSmania18}. We shall distinguish two subcases here; one where the parameters of the logistic map are drawn from a smooth heterogeneous distribution and one where they are drawn from a non-smooth distribution.  For concreteness, we consider perturbations of the following modified logistic map,
\begin{widetext}
 \begin{align}
{\small{
\left( q_{n+1}^{(j)},r_{n+1}^{(j)}\right) 
 = 
\begin{cases}
\left( q_{n}^{(j)},2r_{n}^{(j)}\right) & r_n^{(j)}<\tfrac{1}{2}\\
\left(a^{(j)}\,q_n^{(j)}(1-q_n^{(j)}) + h(\qj_n,\Phi_n) + \eps g(\qj_n),2 r_n^{(j)}-1\right) &  r_n^{(j)}\ge\tfrac{1}{2}
\end{cases}
}},
\label{e.logistic}
\end{align}
\end{widetext}
where the logistic map parameters $\aj$ are sampled from a distribution $\nu(a)da$  \cite{WormellGottwald18}. The action of this map on $q$ is plotted in Figure~\ref{f.cestnepasunmaplogistique}(a). To investigate the linear response properties of this system we choose the perturbation function
\begin{align}
g(\qj_n) = 4(\qj_n(1-\qj_n))^2.
\label{e.g}
\end{align}

The function $h(\qj_n,\Phi_n)$ denotes the mean field coupling which is set to $h\equiv 0$ in the uncoupled scenario. In the coupled scenario we will consider the mean field coupling
\begin{align}
h(\qj_n,\Phi_n) = (1-2\qj_n)\qj_n(1-\qj_n) \tanh \Phi_n.
\label{e.h}
\end{align}
The effect of $g$ and $h$ are also plotted in Figure~\ref{f.cestnepasunmaplogistique}(a).

We remark that a naive choice of mean field coupling with $h=\Phi_n$ would just lead back to the standard logistic map for some $p^{(j)}_n=\alpha \qj_n +\beta$ with a modified logistic map parameter $\aj=\aj(\Phi_n)$. The mean field $\Phi_n$ is given by (\ref{e.Phi}) and is constructed using
\begin{align*}
\phi(q) = 4T_5(2q-1) + 1,
\end{align*}
where $T_5(x) = 16x^5 - 20x^3 + 5x$ is the 5th Chebyshev polynomial, which oscillates between $\pm 1$ in the domain (see Figure~\ref{f.cestnepasunmaplogistique}(b)). %For the external perturbation $\eps g$ we employ again (\ref{e.g}).

In our numerical simulations we use $\psi(q,r)=q$ for our mean-field observable $\Psi$.

The inclusion of the mixing doubling map dynamics $r_n$ ensures that the overall dynamics is mixing even when the logistic parameters $\aj$ correspond to regular dynamics. The inclusion of the cocycle $r_n$, however, does not alter the invariant measure of the logistic map for constant $\Phi_n$ and the marginal invariant measure of $\qj$  the invariant measure of a logistic map at parameter $\aj$. Hence, notwithstanding any dynamics of $\Phi_n$, the microscopic dynamics (\ref{e.logistic}) violates LRT while being mixing. 

\begin{figure}[tbp]
\centering
\includegraphics[width=0.6\columnwidth]{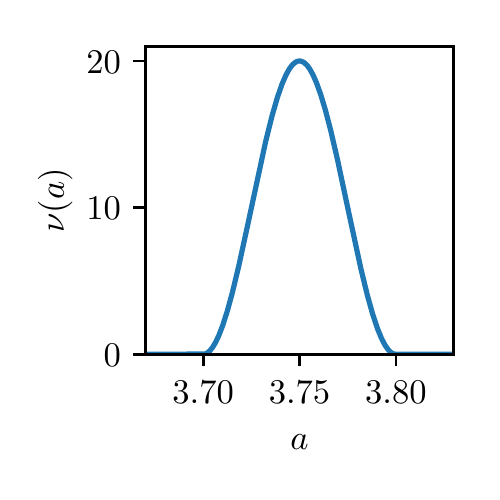}
\caption{Probability density function $\nu(a)$ of the raised cosine distribution (\ref{e.rosine}) with compact support on the interval $[3.7,3.8]$.}
\label{f.rosine}
\end{figure}

In \citet{WormellGottwald18} it was established that the heterogeneity of the parameters $\aj$, as exemplified by the regularity of $\nu(a)$, was crucial in establishing LRT (albeit in a different, less general setting). We therefore consider here two cases: the case when $\nu(a)$ is smooth, in particular at least once-differentiable with respect to $a$, and the case when $\nu(a)$ is non-smooth, for example when $\nu(a)$ is a linear combination of delta functions. Similar to \citet{WormellGottwald18} we choose as a smooth distribution the raised cosine distribution supported on the interval $[3.7,3.8]$, which is given by
\begin{align}
\nu(a) = \mathbf{1}_{[3.7,3.8]} \frac{1}{0.1} \left(1+\cos\left(\frac{a-3.75}{0.05}\pi\right)\right).
\label{e.rosine}
\end{align} 
We have chosen this distribution as it is both compactly supported and resembles a Gaussian distribution (see Figure~\ref{f.rosine}). For a non-smooth distribution we choose the discrete distribution
\begin{align} \nu = \frac{1}{3} (\delta_{3.72}+\delta_{3.75}+\delta_{3.78}) ,
\label{e.deltas}
\end{align}
which has a similar distribution of moments.
%
% FIGURE
%

%
%
% The smoothness of the distribution of the parameters $\nu(a)$ is crucial for the existence of linear response, as already established in \cite{WormellGottwald18}. 

\section{Summary of results}
\label{s.results}

Our main results for these different dynamical scenarios and cases are summarised in Table~\ref{t.result}. We differentiate between the thermodynamic limit $M = \infty$ and the case of a large, finite heat bath size $M$ (which may not necessarily approach a smooth limit as $M \to \infty$). We summarise the dynamical mechanisms leading to the comprehensive picture provided in Table~\ref{t.result}, which are to the best of our knowledge hitherto unknown. The following sections will establish these findings in detail.
\begin{itemize}
	\item Macroscopic mean field observables generated by an appropriate heterogeneous set of microscopic chaotic systems may exhibit linear response, even if the individual members of those systems may not individually have LRT (Section \ref{s.uncoupledB}).
	
%	``Certain'' forcing means either the perturbations are constant in time (Section \ref{s.uncoupledB}) or the perturbations are bounded and the system has been kicked at some (possibly long) time in the past (Section ??). The necessity for a kick in the latter case is due to the fact that noise-like perturbations will, often dramatically, collapse narrow periodic windows.
	
         \item In the thermodynamic limit, macroscopic observables obey a law of large numbers. If the microscopic dynamics is mixing, this leads in the case of no back-coupling to trivial macroscopic dynamics (Section~\ref{s.stochreduction}); if the microscopic dynamics is coupled via its mean field and provided the microscopic dynamics collectively obeys LRT, one can derive a smooth non-Markovian closure for macroscopic variables (Section~\ref{s.surrogate}). In the latter case, if the macroscopic dynamics converges to a fixed point or to a limit cycle, the macroscopic mean field observables satisfy LRT in the thermodynamic limit (Section~\ref{s.coupled_Triv}). However, the reduced macroscopic dynamics may also converge to a chaotic dynamical system which violates LRT. This is possible even if the individual microscopic dynamics is uniformly hyperbolic (Section~\ref{s.coupled_NonTriv}).  %In the latter case either converge to constant equilibria  and \ref{s.coupled_Triv}), or, when there is back-coupling, may have their own emergent nontrivial dynamics which may or may not obey LRT (Section~\ref{s.coupled_NonTriv}).

	\item In finite ensembles with $M<\infty$ the mean field involves an ${\mathcal{O}}(1/\sqrt{M})$ correction to the thermodynamic mean field dynamics, which may not obey LRT (Section \ref{s.stochreduction}). The possible violation of LRT of macroscopic observables, however, is not detectable for practical purposes, and the observed linear response is determined by the linear response property of the thermodynamic limit. We call this behaviour {\em{approximate LRT}}.
	
	\item In finite mean-field coupled systems such as (\ref{e.scen2}), macroscopic mean fields typically satisfy a central limit theorem. As a result, the back-coupling of the mean field introduces a small ``noise'' into the microscopic systems, which can induce linear response. The statistical properties of this dynamic self-generated noise and its linear response properties are determined by the linear response property of the thermodynamic limit. In the case of failure of linear response in the thermodynamic limit, the convergence to the thermodynamic limit is approached for finite large $M$ through the creation of saddle-node bifurcations (with associated multistability) which become increasingly dense in $\eps$ (Section \ref{s.coupled_TrivC}).
\end{itemize}

%
% TABLE
%
\begin{table*}
%\begin{tabular}{ |p{6.75cm}p{2.4cm}|p{3cm}|p{3cm} | }
\begin{tabular}{ |cc||c|c| }
%\begin{tabularx}\textwidth{ |X|X|X|X|  }\toprule
\hline
& &\multicolumn{2}{c|}{macroscopic observables} \\
\hline
microscopic subsystem& & uncoupled & coupled\\
\hline
\hline
\multirow{2}{0.45\textwidth}{$f$ satisfies LRT } & finite $M$\; &\cmark&\cmark\\
                                                 & $M\to\infty$ \;&\cmark&\dmark\\
\multirow{2}{0.45\textwidth}{$f$ violates LRT with smooth $\nu(a)da$} & finite $M$\; &(\cmark)&(\cmark)\\
                                                 & $M\to\infty$ \; &\cmark&\dmark\\
\multirow{2}{0.45\textwidth}{$f$ violates LRT with non-smooth $\nu(a)da$ } & finite $M$\;  &\xmark&(\cmark)\\
                                                 & $M\to\infty$ \;&\xmark&\xmark\\
\hline
\end{tabular}
\caption{Summary of our main result. The checkmarks \cmark\, denote cases when the macroscopic observable $\Psi$ enjoys LRT. The bracketed checkmarks (\cmark)\, denote cases of approximate LRT, when LRT is satisfied for practical purposes. The cross-marks \xmark\, denote cases when LRT is violated for the macroscopic observable. The star \dmark\, denotes cases when LRT may or may not be satisfied depending on the linear response of the limiting dynamics of the macroscopic observable (see Section~\ref{s.coupled}).}
\label{t.result}
\end{table*}

In the following we provide numerical evidence and theoretical arguments corroborating these results. We first consider the case of macroscopic observables of an uncoupled heat bath before considering the case of macroscopic observables of a mean field coupled heat bath.

%%%%%%%%%%%%%%%%%%%%%%%%%%%%%%%%%%%%%%%%%%%%%%%%%%%%%%%%%%

\section{Macroscopic observables of uncoupled microscopic subsystems}
\label{s.uncoupled}
We are concerned with the behaviour of averages of the macroscopic observable $\Psi$. We distinguish here two averages; the average with respect to initial conditions of $\qj$, %(which may depend on time when we address the scenario involving backcoupling of the mean field) 
which we denote by $\E$, and the average over the independently chosen logistic map parameters distributed according to $\nu(a)$ which we denote by angular brackets $\langle \cdot\rangle$. In real systems (for which the parameters $\aj$ are selected once only), the average relevant for linear response  is $\E$, the expectation with respect to initial conditions.

We describe a stochastic reduction of the mean field dynamics in Section \ref{s.stochreduction} and then in Section \ref{s.uncoupledLRT} discuss the linear response properties for each of the three kinds of microscopic subsystems that we outlined in Section \ref{s.model}: in Table~\ref{t.result} these are covered in the rows corresponding to the uncoupled macroscopic observables.

\subsection{Stochastic reduction of mean field dynamics}
\label{s.stochreduction}

The average with respect to initial conditions is written as
\begin{align*}
\E\psi(\qj) = \Rone{\int} \psi(q) d\mu^{a^{(j)}}(q), 
\end{align*}
where $\mu^{\aj}(\qj)$ is the invariant measure of $\qj$. The Law of Large Numbers then reads as
\begin{align} 
\langle \E \Psi \rangle = \iint \psi(q) d\mu^{a}(q) d\nu(a).
\label{eq:EEPhi} 
\end{align}
(In view of Section~\ref{s.coupled} where the mean field coupling is considered and the $q^{(j)}$ depend on a time-varying driver, we remark that in that case averages are computed with a time dependent measure $\mu^{\aj}_n(\qj)$.)

We first establish the case of LRT for a finite heat bath. For large but finite system size $M$, both averages are equipped with their own finite size correction, described by the central limit theorem. In equilibrium each ensemble member $q\j_n$, at a given time $n$, is an independent sample from the invariant measure $\mu^{a\j}$. Macroscopic observables $\Psi$, as defined in (\ref{e.Psi}), can be approximated using the central limit theorem and the independence of the $q\j$ by
\begin{align} 
\Psi_n = \E \Psi + \frac{1}{\sqrt{M}}\zeta_n +o(1/\sqrt{M}), 
\label{e.PhiLLN} 
\end{align}
where the expectation value
\begin{align*} 
\E \Psi = \frac{1}{M} \sum_{j=1}^M \int \psi(q) d\mu^{a\j}(q)  
\end{align*}
is over initial conditions $\qj$ at fixed $\aj$. The random mean-zero Gaussian process $\zeta_n$ has autocovariance function $C^\zeta(m)$ with
\begin{align}
C^\zeta(m)&=\cov(\zeta_n,\zeta_{n+m})
=\lim_{M\to\infty}\frac{1}{M}\sum_{j=1}^M \E[\psi^{(j)}_0\psi^{(j)}_m] \nonumber \\
&  = \langle \E[\psi_0\psi_m] \rangle .
\label{e.covzetan}
\end{align}
The existence of a central limit theorem is guaranteed for unimodal maps using results of Lyubich \cite{Lyubich02} who proved that almost every non-regular logistic parameter satisfies the so-called Collet-Eckmann condition \cite{colletEckmann83}, which then implies the existence of good statistical properties including the central limit theorem\cite{AlvesEtAl04,MelbourneNicol08}. We remark that the parameters determining the process $\zeta_n$ in the finite-$M$ case, such as the covariance (\ref{e.covzetan}), have the same LRT properties as the associated thermodynamic limit $\langle \E^\eps \Psi \rangle$.\\

The independent sampling of the $a^{(j)}$ allows for a further application of the central limit theorem, and we can write 
\begin{align}	
\E \Psi_n = \langle \E \Psi_n \rangle + \frac{1}{\sqrt{M}} \eta + o(1/\sqrt{M}),
\label{e.EPhiLLN} 
\end{align}
where the random variable $\eta$ is, for fixed $\eps$, a mean-zero Gaussian variable. As a function of $\eps$, $\eta$ is a Gaussian process with covariance 
\begin{align}
\langle \eta^\eps \eta^{\eps'} \rangle = \langle \E^\eps[\psi] \E^{\eps'}[\psi] \rangle - \langle \E^\eps[\psi] \rangle \langle \E^{\eps'}[\psi] \rangle,
\label{e.covetan}
\end{align}
and typically is no more differentiable with respect to $\eps$ than $\E^\eps[\psi]$, which implies that LRT is violated for finite $M$ if the microscopic subsystems do not individually satisfy LRT. However, for finite $M\gg 1$ the response term $\langle \E\Psi_n\rangle$ dominates over the contribution of $\eta$ and the violation of LRT can only be detected for vanishingly small values of $\eps$. We call this instance of LRT for all practical purposes {\em{approximate LRT}}.\\
%This addresses the bracketed nature of the (\cmark) checkmark in the first column in Table~\ref{t.result}.\\ 

We remark that, notwithstanding the rough parameter selection error discussed above, and recalling that the linear response of the process $\zeta_n$ is determined by the linear response property of the associated thermodynamic limit,  the overall linear response of $\Psi_n$ depends entirely on whether the thermodynamic limit $\langle \E\Psi_n\rangle$ satisfies LRT or not. We discuss this question in the next section.

%However, if the support of $\nu(a)da$ is sufficiently small so that the variation of $\E[\psi]$ over the parameters $\aj$ is small compared to the typical variance $R(0)=\E[(\psi-\E[\psi])^2]$ for these parameters, then the contribution of the smoothing Gaussian process $\zeta_n$ will be dominant, and the breakdown of LRT due to $\eta$ will not be detectable using statistical significance tests \cite{GottwaldEtAl16}. On the other hand, it is clear that if the support of $\nu(a)da$ is chosen too small, the lack of heterogeneity restricts the validity of LRT to a small range of perturbation sizes $\eps$ \cite{WormellGottwald18}. 
%where the angle brackets denote the expectation with respect to the logistic map parameters $a^{(j)}$, hence
%
%The finite size corrections of $\O\left(\frac{1}{\sqrt{M}}\right)$ will be discussed in detail in the following section. 

%%%%%%%%%%%%%%%%%%%%%%%%%%%%%%%%%%%%%%%%%%%%%%%%%%%%%%%%%%

\subsection{LRT of thermodynamic limit mean field observables of uncoupled microscopic subsystems}
\label{s.uncoupledLRT}

We now investigate the response of $\langle \E\Psi_n\rangle$, i.e. the thermodynamic limit. We distinguish between three cases: when the microscopic dynamics satisfies LRT, and when the the microscopic dynamics does not satisfy LRT and has a distribution $\nu(a)$ of the parameters which is either smooth or non-smooth.

%We remark that, the parameters of the process $\zeta_n$ in the finite-$M$ cases have the same LRT properties as the thermodynamic limit cases discussed below, and thus notwithstanding the rough parameter selection error discussed in the previous section, the finite $M$ case for any system has the same LRT properties as in the thermodynamic limit.

\subsubsection{The microscopic subsystems satisfy LRT}
\label{s.uncoupledA}
If the microscopic dynamics obeys LRT, as is the case for uniformly expanding maps such as (\ref{e.Lanford}), which will be considered in Section~\ref{s.coupled_TrivA}, it is clear that LRT also holds for macroscopic observables defined in (\ref{e.Psi}). For finite heat bath sizes $M$, we have
\begin{align*}
\frac{d}{d\eps}\E^\eps \Psi_n = \frac{1}{M}\sum_{j=1}^M \frac{d}{d\eps}\E \psi(\qj_n)
\end{align*} 
and the macroscopic observable $\Psi$ has LRT since the $M$ subsystems individually satisfy LRT with uniformly bounded  $\frac{d}{d\eps} d\mu^{(a_j,\eps)}$. The validity of LRT carries over to the thermodynamic limit with
\begin{align*}
\frac{d}{d\eps}\langle \E^\eps \Psi_n \rangle &= \iint \psi^\eps(q) \frac{d}{d\eps} d\mu^{(a,\eps)}(q) \nu(a)da.
\end{align*} 
Note that we may allow for a $\nu$-measure zero subset of subsystems at any given $\eps$ to individually violate LRT, and still obtain LRT for the macroscopic observable $\Psi$ in the $M\to\infty$ limit. In this case, however, the $\eta$ correction may not be differentiable, and we observe approximate LRT.

%%%%%%%%%%%%%%%%%%%%%%%%%%%%%%%%%%%%%%%%%%%%%%%%%%%%%%%%%%

\subsubsection{The microscopic subsystems do not satisfy LRT but are appropriately heterogeneous}
\label{s.uncoupledB}
As typical microscopic dynamics leading to violation of LRT we consider the modified logistic map (\ref{e.logistic}), which is perturbed in $\eps$ by the function (\ref{e.g}), and whose parameters are drawn from the smooth raised cosine distribution (\ref{e.rosine}) which is three-times continuously differentiable (i.e. $\nu(a)$ is $C^3$). 

\begin{figure}[htbp]
	\centering
	\includegraphics[width=\columnwidth]{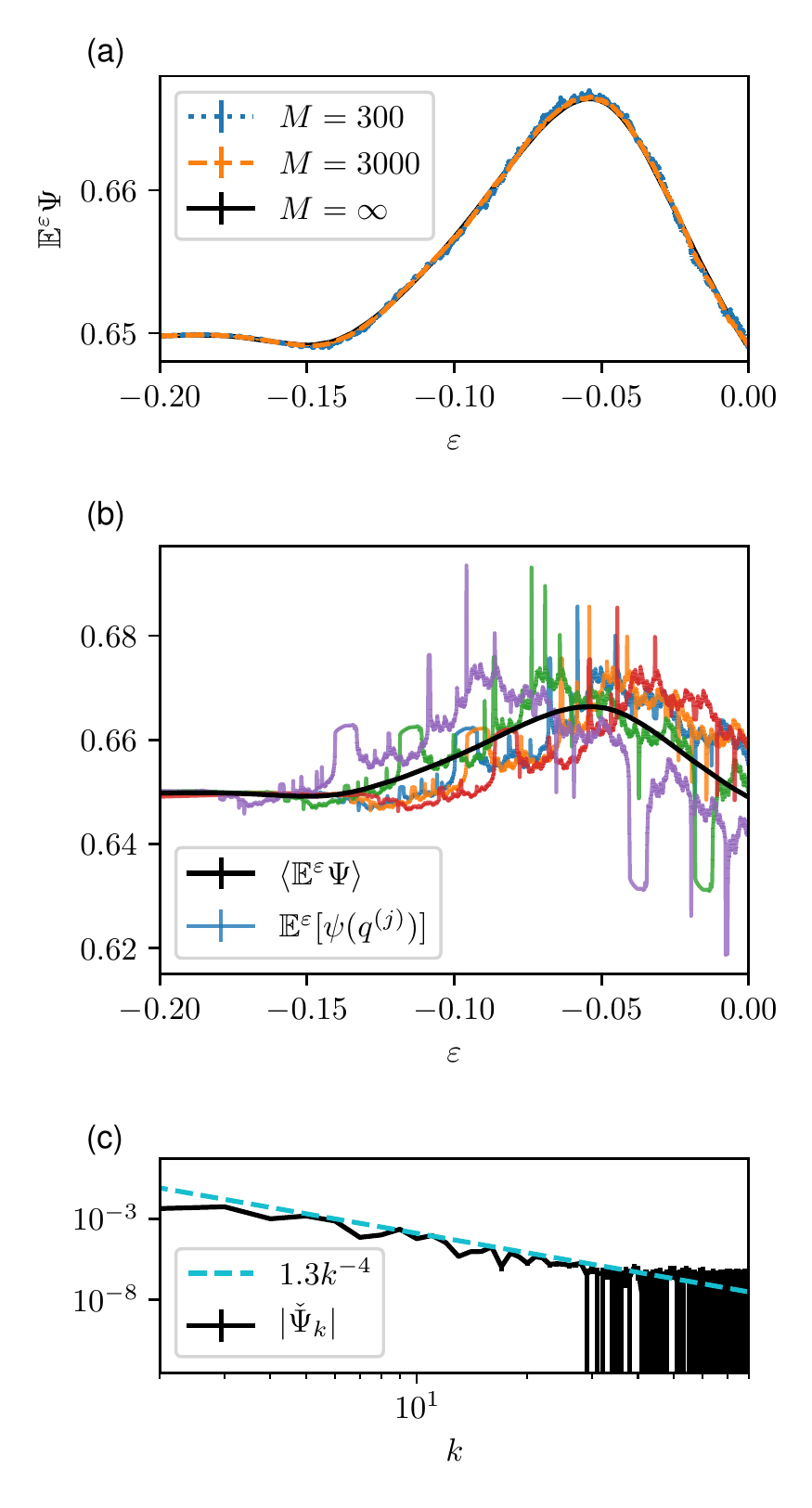}
	\caption{Response term $\E^\eps \Psi$ for a perturbation of the form (\ref{e.g}) for an uncoupled heat bath scenario in the case when the microscopic dynamics is given by the logistic map (\ref{e.logistic}), which does not satisfy LRT, and the logistic map parameters are sampled from the raised cosine distribution (\ref{e.rosine}). For different values of $\eps$ we employ a total of $10^5$ iterates to estimate $\E^\eps \Psi$ as a temporal average. (a) Plots for finite $M$: 95\% confidence intervals were estimated from $10$ realisations differing in the initial conditions of the heat bath, and are not visible. Thermodynamic limit curve (black), confidence intervals also not visible, was estimated from $50$ realisations of $3000$ iterates for $10^6$ parameters $a\j$ independently selected for each $\eps$. (b) Thermodynamic limit of $\E^\eps\Psi$ (black), with LRT-violating response of microscopic variables $\E^\eps\psi(q\j)$ (coloured lines), estimated from $10$ realisations of $10^6$ iterates each. (c) Estimate of Chebyshev coefficients $\sum_{k=0}^\infty \check \Psi_k T_k(0.1^{-1}(x+0.1)) := \langle \E^\eps \Psi \rangle$.
	}
	\label{f.Scenario1_Case2}
\end{figure}

Figure~\ref{f.Scenario1_Case2} provides numerical evidence that, for these maps, the macroscopic observable $\Psi$ with $\psi(x)=x$ has linear response for a wide range in $\eps$. To determine the smoothness of $\langle \E^\eps \Psi \rangle$, we determine its Chebyshev coefficients on a Chebyshev roots grid of $1000$ points. \Rone{It is well known that any smooth function $f(x)$ can be expressed as an infinite series of Chebyshev polynomials $T_k(x)$ as $f(x)=\sum_{k=0}^\infty f_k T_k(x)$ and the degree of differentiability of the function is given by the decay of its Chebyshev coefficients $f_k$ \cite{Trefethen13}.} We find that the Chebyshev coefficients of $\langle \E^\eps \Psi\rangle$ decay as $\O(k^{-4})$, which is indicative of $\langle \E^\eps \Psi\rangle$ being between $C^{3^-}$ and $C^4$ differentiable over a large interval: this level of differentiability, as we will see below, is connected to the smoothness of the raised-cosine distribution (\ref{e.rosine}), which is $C^3$\,\cite{WormellGottwald18}. 
%{\bf{Are we sure that in our situation of smearing there is a similar argument??? If not maybe we should write: It is curious that the degree of differentiability equals the smoothness of the raised-cosine distribution (\ref{e.rosine}), which is also $C^3$; in the case of LRT for a distinguished degree of freedom this was made rigorous\cite{WormellGottwald18} }}. 
We have also employed the test statistics for higher-order linear response developed in \citet{GottwaldEtAl16} to test the null-hypothesis of $\langle \E^\eps\Psi \rangle$ being well-approximated by a linear combination of $T_k(0.1^{-1}(\eps+0.1)), k = 0,\ldots,60$ for $\eps \in [-0.2,0]$, i.e. that the response is in fact a smooth function. This test is summarised in Appendix \ref{a.lrttest}. We used the aforementioned Chebyshev grid simulating $1,000,000$ different randomly selected parameters with $50$ runs of $3000$ timesteps each, and quantified Birkhoff variance both within parameters and between parameters (i.e. arising from the random parameter selection), using standard ANOVA methods \cite{Rice}. We obtain a $p$-value of $0.26$, consistent with the null-hypothesis of a smooth response. 

%
% FIGURE
%
We note that in Figure~\ref{f.Scenario1_Case2} the response of $\E^\eps \Psi$ for systems with finite $M$ have (barely) noticeable rough deviations from the $M\to\infty$ limit: these non-smooth deviations arise from the finite sampling of parameters $\aj$ from $\nu(a)$, approximated by the random variable $\eta$ defined in (\ref{e.EPhiLLN}), as discussed above.\\

\begin{figure}[htbp]
	\centering
	\includegraphics[width=\columnwidth]{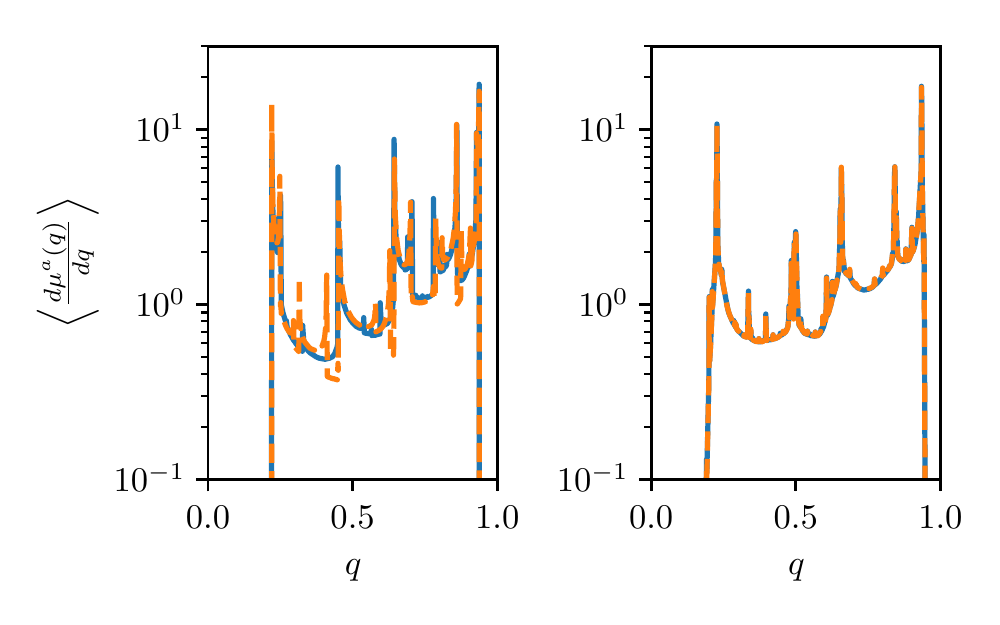}
	(a) $\qquad\qquad\qquad\qquad\qquad\qquad\qquad\qquad$ (b)
	\caption{Histogram of the averaged acim of the $q\j$ for the logistic map system (\ref{e.logistic}) with (a) $\nu = \delta_{3.75}$ and (b) $\nu$ the raised-cosine distribution (\ref{e.rosine}), for $\eps = -5 \times 10^{-4}$ (orange dashes) and $\eps = 0$ (blue line).}
	\label{f.AveragedAcim}
\end{figure}

We now provide a heuristic argument how averaging over a smooth distribution such as the raised cosine distribution (\ref{e.rosine}) can lead to LRT for the macroscopic observable $\Psi$. Let us first recall the dynamic reason of why LRT is violated in the logistic map. We follow here \citet{Ruelle09b} in our exposition. The critical point $q=c$ with  $f'(c) = 0$ leads to a non-uniform compression of the phase space around $q=c$: an initially smooth initial density which contains the critical point in its support is pushed forward under the dynamics to a non-smooth density with a spike with an inverse square-root singularity at $q=f(c)$. This compression is repeated to produce further inverse square-root singularities at locations $q_n=f^n(c)$ of amplitudes asymptotically proportional to $\alpha^{-n/2}$ (and thus contain a probability mass of order $\alpha^{n/2}$), where $1<\alpha$ denotes the Lyapunov multiplier of the logistic map. The result is that the invariant density contains an infinite number of spikes of decreasing amplitude. The effect of the perturbation, by modifying the forward orbit of the critical point $(f^n(c))_{n\in\mathbb{N}}$, is to displace these spikes. Because the map $f$ is chaotic and thus exponentially sensitive to perturbations, spikes move with an instantaneous speed of the order of $\alpha^n$ per unit change of the perturbation. This scenario is illustrated in Figure~\ref{f.AveragedAcim}(a) where we show the absolutely continuous invariant measure (acim), averaged over the heat bath, of an perturbed and of a slightly perturbed dynamics. The high speed of the small spikes (i.e. those with large $n$) in conjunction with their relatively large probability mass implies that the sum of their (distributional) derivatives diverges, leading to breakdown of linear response. The reader is referred to \citet{Ruelle09b} for details and to \citet{GottwaldEtAl16} for a numerical illustration. For comparison we also depict in Figure~\ref{f.AveragedAcim}(b) the averaged invariant measure for the case when the parameters are heterogeneously drawn from a raised-cosine distribution and for which we showed above that LRT is valid. The averaging over the heterogenously drawn parameters clearly implies a smoothed invariant measure of the logistic map. 

This may be seen as analogous to a recent heuristic argument for linear response in general non-hyperbolic systems\cite{Ruelle18}: rough contributions to the response caused by singularities in the physical measures that arise from stable manifold-unstable manifold tangencies average out if these singularities distribute themselves suitably evenly.

In the specific case where the microscopic dynamics evolves under unimodal maps such as those studied here, we can make a more concrete argument for the effect of the smearing out of the fast displacement of the small spikes upon perturbation. It is conjectured by Avila {\em{et al}} \cite{AvilaEtAl03}, there exists an $\eps$-dependent analytic function $\alpha(a,\eps)$ of the invariant measures, such that the map with parameters $(a,\eps)$ is topologically conjugate to the map with parameters $(\alpha(a,\eps),0)$. Unimodal maps, at least those of Benedicks-Carleson type, have linear response within topological conjugacy classes \cite{BaladiSmania12}, and as a result we can say
\begin{align}
\langle \E^\eps \Psi\rangle
&=\iint \psi(q)d\mu_n^{a,\eps}\nu(a)da \notag \\
&= \iint \psi(q)d\mu_n^{\alpha(a,\eps),0}\nu(a)da + \mathrm{h.o.t.} \notag \\
&= \iint \psi(q)d\mu_n^{\alpha,0}\nu(a(\alpha,\eps))\frac{da}{d\alpha}d\alpha  + \mathrm{h.o.t.}, \label{e.avila}
\end{align}
where the higher order terms capture the response of the $\mu_n$ under the topological conjugacy-preserving parameter change from $(\alpha(a,\eps),0)$ to $(a,\eps)$, and are thus similarly smooth. Hence the existence of linear response $d\langle \E^\eps \Psi\rangle/d\eps$ (and by a similar argument higher-order response) is guaranteed, provided that the distribution of the parameters of the logistic map $\nu(a)$ is at least once continuously differentiable (and provided $\alpha(a,\eps)$ is analytic in a sufficiently uniform way\cite{AvilaEtAl03}). The linear response of $\langle \E^\eps \Psi\rangle$ was numerically confirmed in Figure~\ref{f.Scenario1_Case2}. Since in at least one-dimensional systems topological conjugacy classes form manifolds of finite codimension, we believe this argument will generalise to more general maps, provided the space of parameters $a\j$ is sufficiently high-dimensional. 
\\

\begin{figure}[htbp]
	\centering
	\includegraphics[width=\columnwidth]{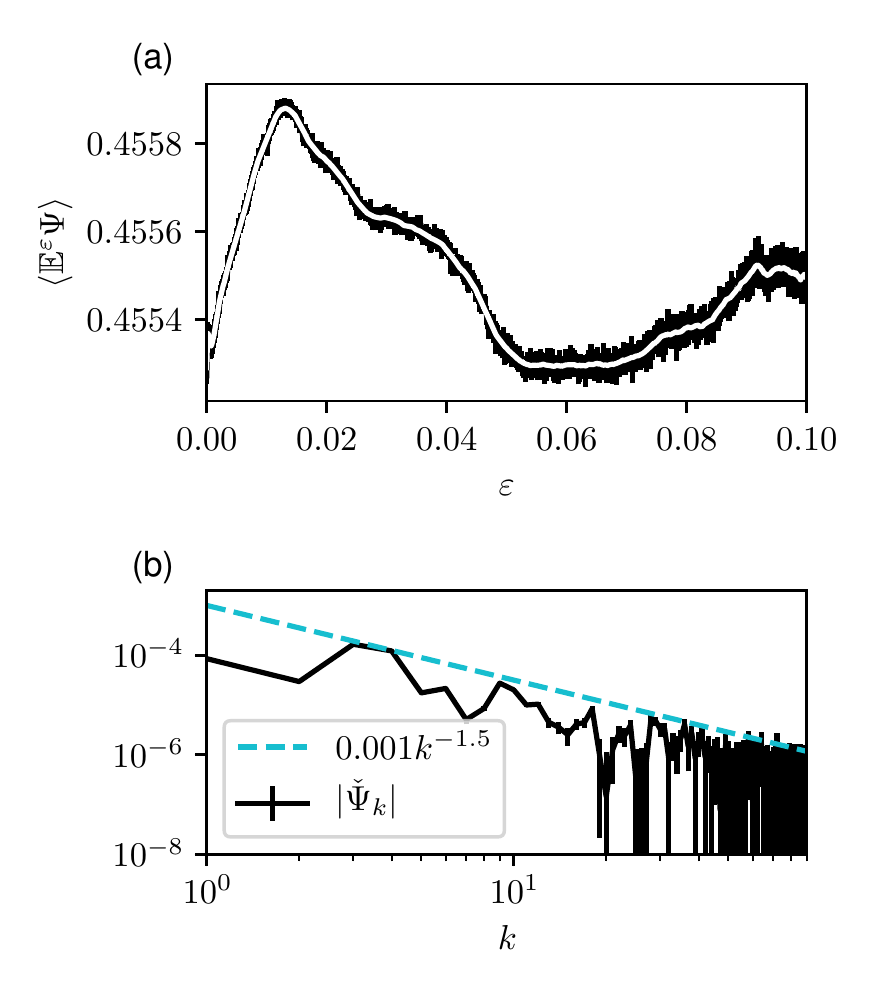}
%	(a) $\qquad\qquad\qquad\qquad\qquad\qquad\qquad\qquad$ (b)
\caption{Response term $\E^\eps \Psi$ for an uncoupled heat bath scenario for the map (\ref{e.torus}) where the parameters are sampled from a raised cosine distribution (\ref{e.rosine2}). (a) Infinite $M$ limit with confidence intervals (black) and $21$-point moving average with confidence intervals (white) from $15$ realisations of $6000$ iterates for $10^6$ parameters $a\j$ independently selected for each $\eps$. (b) Estimate of Chebyshev coefficients $\sum_{k=0}^\infty \check \Psi_k T_k(0.05^{-1}(x+-0.05)) := \langle \E^\eps \Psi \rangle$.}
	\label{f.Scenario1_qorus}
\end{figure}

\begin{figure}[htbp]
	\centering
	\includegraphics[width=\columnwidth]{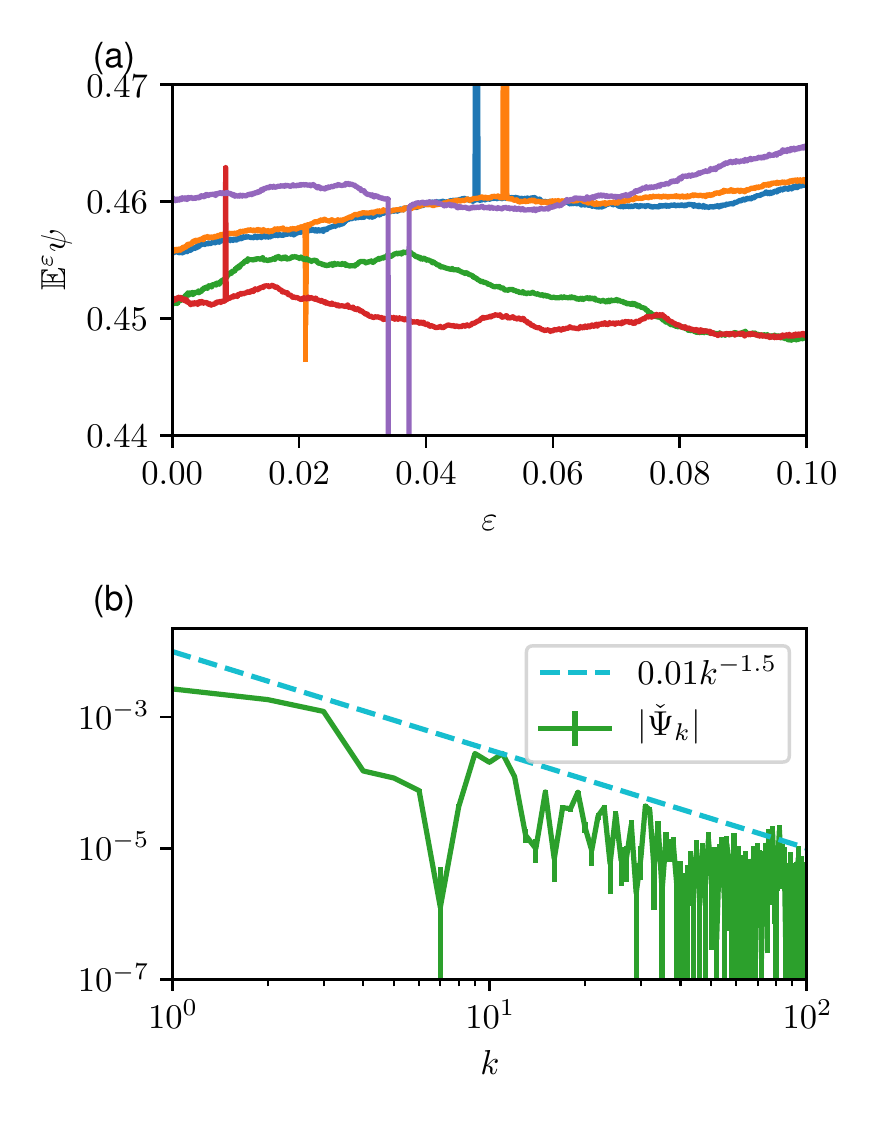}
	%	(a) $\qquad\qquad\qquad\qquad\qquad\qquad\qquad\qquad$ (b)
	\caption{Individual response terms $\E^\eps \psi(x\j,y\j)$ with confidence intervals for the map (\ref{e.torus}) where the parameters selected from the raised cosine distribution (\ref{e.rosine2}). (a) Response for five randomly selected microscopic variables. The large jumps of the response outside the figure correspond to regions of regular dynamics. (b) Estimate of Chebyshev coefficients $\sum_{k=0}^\infty \check \Psi_k T_k(0.05^{-1}(x+-0.05)) := \langle \E^\eps \Psi \rangle$ for one of the variables in (a).}
	\label{f.Scenario1_indiv}
\end{figure}

However, we caution that our smearing argument may not generalise to other systems, at least when the support of the parameters $a$ is one-dimensional. In Figure~\ref{f.Scenario1_qorus} we present numerical evidence demonstrating that mean-field averaging fails to improve the linear response of a unimodal map of the torus $(\R/\Z)^2$ for heterogeneously distributed parameters, given by
\begin{align} 
 %(x\j_{n+1}, y\j_{n+1}) = (x\j_n + a\j y \sin \pi x\j_n, y\j_n + a\j \sin \pi (x\j_n + y\j_n) + \eps) \mod 1. 
 x\j_{n+1} & =x\j_n + a\j y \sin \pi x\j_n \mod 1\\
 y\j_{n+1} &= y\j_n + a\j \sin \pi (x\j_n + y\j_n) + \eps \mod 1,
\label{e.torus} 
\end{align}
for $j = 1,\ldots, M$. The parameters $a\j$ are again distributed according to a raised cosine distribution with support on $[3.7,4.3]$\
\begin{equation}
\nu(a) = \mathbf{1}_{[3.7,4.3]} \frac{1}{0.1} \left(1+\cos\left(\frac{a-3}{0.3}\pi\right)\right),
\label{e.rosine2}
\end{equation} 
and the mean-field observable is given by $\Psi = \frac{1}{M} \sum_{j=1}^M \psi(x\j,y\j)$ where $\psi(x,y) = x$, as before. We tested the null-hypothesis of $\langle \E^\eps\Psi \rangle$ being well-approximated by a linear combination of $T_k(0.05^{-1}(\eps-0.05)), k = 0,\ldots,100$ for $\eps \in [0,0.1]$, and obtained a p-value $p = 0.49$, consistent with the null hypothesis. However in Figure~\ref{f.Scenario1_qorus} we see that the estimated Chebyshev coefficients decay approximately as $\O(k^{-1.5})$ which is slower than $\O(k^{-4})$ seen in the one-dimensional unimodal example, indicating a rather low-order differentiability. This level corresponds quite closely to that obtained for the expectation value $\E^\eps \psi(x\j,y\j)$ of a single microscopic systems, as illustrated in Figure~\ref{f.Scenario1_indiv}; hence the averaging over parameters appears only to be smoothing out the large jumps arising from periodic windows but does not improve the degree of differentiability.
 %Note that in many systems with sufficiently large stable dimension (for which a sufficiently large number of unstable Lyapunov exponents are necessary), some kind of low-order smooth response is generically expected \cite{Ruelle18}.

%
%\begin{figure}[htbp]
%\centering
%\includegraphics[width=0.7\linewidth]{figures/lcpet.eps}
%\includegraphics[width=0.7\linewidth]{figure04.eps}
%\caption{Overlaying two a.c.i.m.s of the logistic map with $a=3.8$ (maroon) and $a=3.8+6.05\times 10^{-4}$ (light blue). The displacement of a smaller spike with $n=17$ is illustrated by the two ovals, and that of a large spike with $n=3$ with a dashed box.}
%\label{fig.lcpet}
%\end{figure}
%

%

%%%%%%%%%%%%%%%%%%%%%%%%%%%%%%%%%%%%%%%%%%%%%%%%%%%%%%%%%%

\subsubsection{The microscopic subsystems do not satisfy LRT and are not appropriately heterogeneous}
\label{s.uncoupledC}
\begin{figure}[htbp]
\centering
\includegraphics[width=\columnwidth]{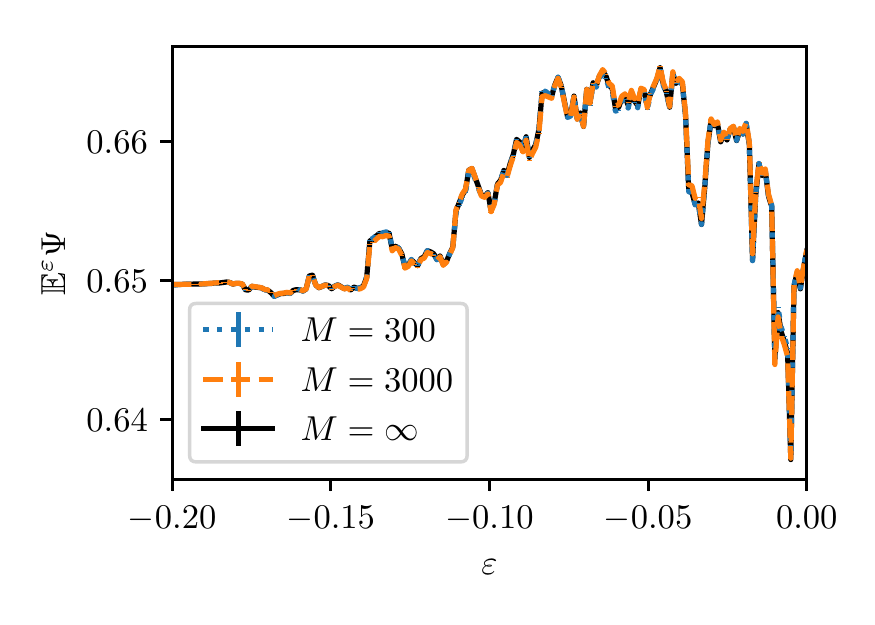}
\caption{ Response term $ \E^\eps \Psi$ for a perturbation of the form (\ref{e.g}) for an uncoupled heat bath scenario in the case when the microscopic dynamics is given by the logistic map (\ref{e.logistic}), which does not satisfy LRT, and the logistic map parameters are distributed as in (\ref{e.nurough}) with $\nu = \tfrac{1}{3}(\delta_{3.72} + \delta_{3.75} + \delta_{3.78})$. Error bars were estimated from $10$ realisations differing in the initial conditions of the heat bath, and are not visible.
}
\label{f.Scenario1_Case3}
\end{figure}

If the microscopic dynamics does not obey LRT and the logistic map parameters are non-smoothly distributed, then LRT fails for macroscopic observables (\ref{e.Psi}), independent of whether the heat bath is finite or infinite. In this case the averaging over the heat bath variables does not provide the necessary smearing of the non-smoothness of the perturbed invariant measures $\mu^{a,\varepsilon}$. To illustrate this, consider  the following non-smooth parameter distribution
\begin{align}
\nu(a) = \sum_{k=1}^p w_k \delta(a-a_k),
\label{e.nurough}
\end{align}
where at least one of the logistic map parameters $a_k$ corresponds to chaotic dynamics. The invariant measures $\mu^{a_k,\eps}$ for fixed parameters $a_k$ are not differentiable with respect to the perturbation size $\eps$ per assumption. The averaging over the independent heat bath variables only involves finitely many logistic parameter values, and hence in this situation averaging is not able to smear the effect of the non-smoothness of the finite number of associated invariant measures $\mu^{a_j}$. In Figure~\ref{f.Scenario1_Case3} we show the response $\E^\eps \Psi$ which as expected exhibits non-smooth behaviour upon varying the strength of the perturbation $\eps$. The response term $\E^\eps\Psi$ quickly converges and for $M=300$ is almost indistinguishable by eye from the response $\langle \E^\eps \Psi \rangle$ in the thermodynamic limit.

%
% FIGURE
%
%
%

%%%%%%%%%%%%%%%%%%%%%%%%%%%%%%%%%%%%%%%%%%%%%%%%%%%%%%%%%%

\section{Linear response of macroscopic observables of microscopic subsystems with mean field coupling}
\label{s.coupled}
We now consider the case when the heat bath variables $\qj$ couple via the mean field $\Phi_n$. In Section~\ref{s.surrogate} we derive a non-Markovian closure of the mean-field dynamics, along the lines of the reduction derived in Section~\ref{s.uncoupled}, that is deterministic in the thermodynamic limit and stochastic for finite $M$; in Sections \ref{s.coupled_Triv} and \ref{s.coupled_NonTriv} we study the mean-field dynamics and its linear response using this macroscopic closure.

\subsection{Surrogate approximation of the mean field dynamics}
\label{s.surrogate}

\begin{figure}[tbp]
\centering
\includegraphics[width=0.35\columnwidth]{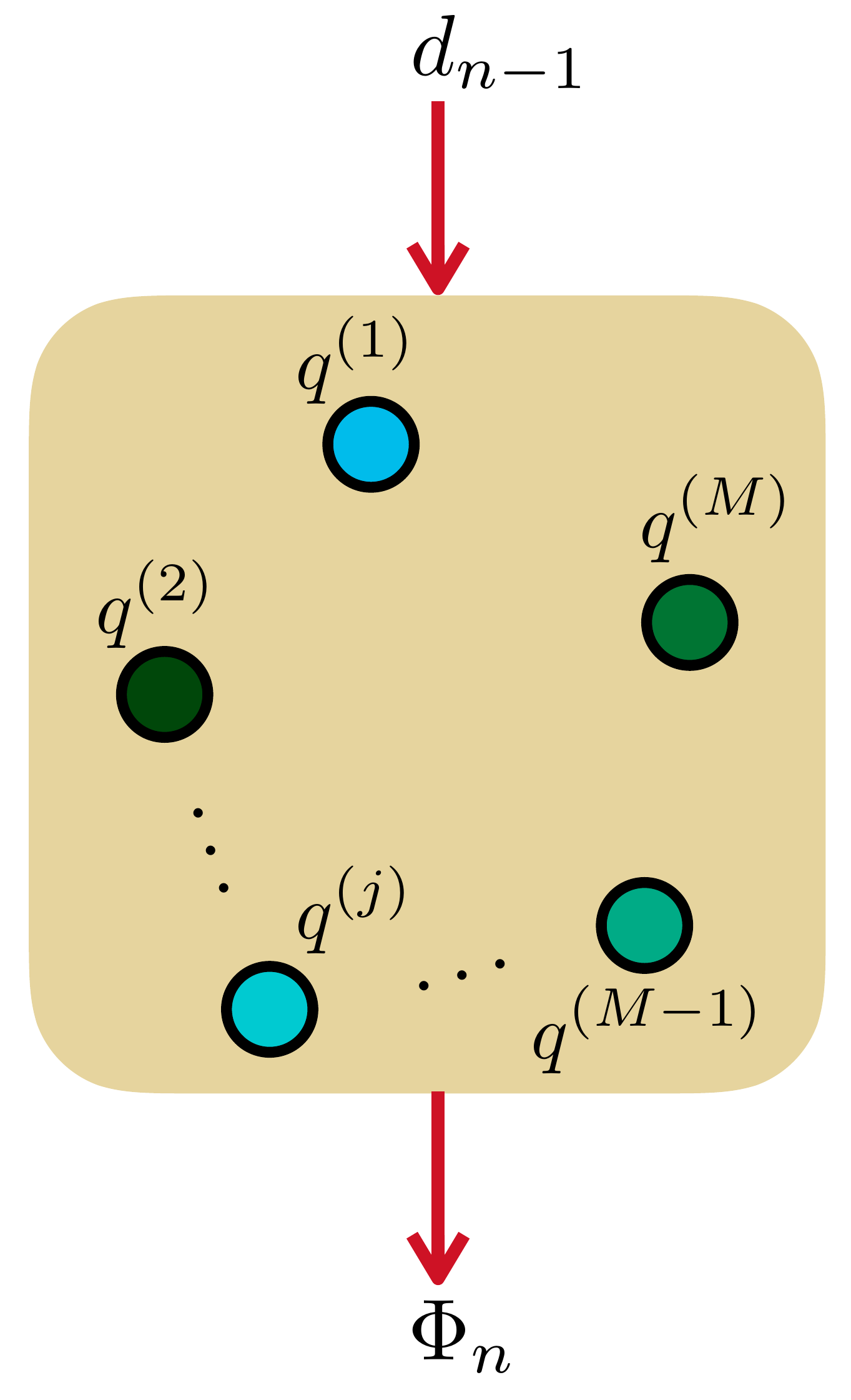}\\
\caption{Sketch of the macroscopic dynamics $\Phi_n$ mediated by the microscopic reservoir.}
\label{f.dynamicalfeedback}
\end{figure}

Before we can study the response of the mean field coupled system to external perturbations $\eps g$, we need to understand the implied macroscopic dynamics $\Phi_n$ generated by the system for the externally unperturbed system with $\eps=0$. To do so we view the system as driven by a prescribed time-dependent external driver $d_n$ rather than the mean field $\Phi_n$, as illustrated in Figure~\ref{f.dynamicalfeedback} (which should be compared with Figure~\ref{f.scenarios}(b)). Hence we replace the mean field coupled dynamics (\ref{e.scen2}) by
\begin{align}
\qj_{n+1}=f(\qj_n,d_n;\aj, \eps)
\label{e.scen2_dn}
\end{align}
for a prescribed external driver $d_n$. In the thermodynamic limit of the mean field coupled system (\ref{e.scen2}) we will see that the macroscopic mean field dynamics is deterministic (see (\ref{e.PhiF}) further down), and the driver $d_n=\Phi_n$ is indeed prescribed by the initial conditions, which is simply the initial distribution of the $q\j$. For large but finite $M$, the $\qj_n$ conditioned on the history of the mean field $(\Phi_{n-1},\ldots)$ can be considered as independently distributed and one can again view the mean field $\Phi_n$ as an externally prescribed noisy driver $d_n$ with specified statistical properties. The surrogate system (\ref{e.scen2_dn}) with the external driver $d_n$ chosen as a random draw of the Gaussian process (defined below in (\ref{e.Phiback})) provides an accurate representation of the statistical behaviour of the original mean field coupled system (\ref{e.scen2}). We have checked that both yield the same linear response $\E^\varepsilon\Psi$, and now set out to study the linear response of the original mean field coupled system via the surrogate system. 

%Because of the large ensemble size, we can understand the $\qj_n$ conditioned on the history of the mean field $(\Phi_{n-1},\ldots)$ as independently distributed: thus, we can study the $\Phi$-dynamics of our self-coupled system using the previous discussion by setting $d_n = \Phi_n$ for all $n$, and thus the CLT approximation (\ref{e.Phiback})  for the dynamics of $\Phi_n$ (and the equivalent reduction for $\Psi_n$) should hold. 

Let us now determine the statistical properties of a macroscopic mean field observable for the driven surrogate system (\ref{e.scen2_dn}). The mean fields $\Phi_n$ and $\Psi_n$ are again Gaussian process with (now time-dependent) statistical properties given again by statistical limit laws, and we write in particular   
\begin{align}
\Phi_n = \langle \mathbb{E}^d[\Phi_n]\rangle +\frac{1}{\sqrt{M}}\notilde\zeta_n + \frac{1}{\sqrt{M}}\notilde\eta_n + o\left(\frac{1}{\sqrt{M}}\right),
\label{e.Phiback}
\end{align}
where
\begin{align}
\mathbb{E}^d[\Phi_n]=\Econd 
\label{e.Econd}
\end{align}
denotes the conditional expectation over the past history of the driver and averages now involve time-dependent measures $\mu_n^{\aj}$. The autocovariances of the mean-zero Gaussian process $\notilde\zeta_n$ are given by a central limit theorem approximation of $\Phi_n - \mathbb{E}^d[\Phi_n]$ with
\begin{align}
\cov[\notilde\zeta_n,\notilde\zeta_{n-k}] = \langle \cov[\phi(q\j_n),\phi(q\j_{n-k})] \rangle, 
\label{e.coupledzetaacv}
\end{align}
where the covariance is defined using the conditional average over the history of the driver (cf. (\ref{e.covzetan})). Note that the autocovariance is not a function of $n-m$ due to the non-Markovian nature of the dynamics. Similarly, a central limit theorem approximation of $\mathbb{E}^d[\Phi_n]- \langle \mathbb{E}^d[\Phi_n]\rangle$, defines the mean-zero Gaussian process $\eta_n$ with autocovariance
\begin{widetext}
\begin{align}
\langle \notilde\eta_n^\eps, \notilde\eta_m^{\eps'}\rangle= \langle \E^{d,\eps}[\phi(q\j_n)] \E^{d,\eps'}[\phi(q\j_m)] \rangle - \langle \E^{d,\eps}[\phi(q\j_n)] \rangle \langle \E^{d,\eps'}[\phi(q\j_m)] \rangle, 
\label{e.covzetan_dn}
\end{align}
\end{widetext}
where again the conditional expectation values $\mathbb{E}$ are used (cf. (\ref{e.covetan})). Note that the Gaussian processes $\notilde\zeta_n$ and $\notilde\eta_n$ are independent. \\

The impulsive response of $\Phi_n$ at a given time to a perturbation of the driving process $d_n \mapsto d_n + \theta_n$, where $|\theta_n|\ll 1$, can be, at least formally, captured by the susceptibility function 
\begin{align}
R_n(z) = \sum_{k=1}^\infty \chi_{n,k} \, z^k,
\label{e.Rz}
\end{align}
defined for complex $z$ with $|z|\leq 1$ \cite{Ruelle04}. The fluctuation coefficients $\chi_{n,k}$ describe the change of the mean field induced by the drivers $\theta_n$ as
\begin{align}
\langle \E^{d+\theta}[ \Phi_n] \rangle - \langle \E^{d} [\Phi_n] \rangle= \sum_{k=1}^\infty \chi_{n,k} \, \theta_{n-k}.
\label{e.dphichi}
\end{align}
The fluctuation coefficients $\chi_{n,k}$ of $\Phi_n$ are given as an average over the microscopic fluctuation coefficients $\chi_{n,k}^a$ of $\phi(q\j_n)$ as
\begin{align}
\chi_{n,k}=\int\chi_{n,k}^a\, \nu(a)da.
\end{align}
A necessary condition for LRT with respect to a bounded driver $\theta_n$ is the summability of the coefficients $\chi_k$. Once LRT with respect to the driver $d_n$ can be shown, we can proceed to study the linear response with respect to the external perturbation with $\eps\neq 0$ (recall that $d_n=\Phi_n$ (cf. Figure~\ref{f.scenarios})). Note that if $\langle \E^d[ \Phi_n]\rangle$ does not satisfy LRT with respect to a perturbation of the driver $d_n$, then it cannot be expected to satisfy LRT with respect to external perturbation.\\
 
%
% FIGURE
%

In the thermodynamic limit, provided the microscopic dynamics is mixing, we can use that the measures $\mu^a_n$ are the physical invariant measures generated by the cocycle $f(\cdot,\Phi_n;a, \eps)$ \cite{Froyland10, Buzzi99} to create a closure of the dynamics of $\Phi_n$ as a deterministic recurrence relation:
\begin{align} 
\Phi_{n} = \langle\E^{d=\Phi}[\Phi_n]\rangle =: F(\Phi_{n-1},\Phi_{n-2},\ldots; \eps). 
\label{e.PhiF} 
\end{align}
If the mixing times of the $\qj$ are much shorter than a delay $k_*$, then the effect of the driving $\Phi_{n-k}$ for $k > k_*$ is minimal and the mean field dynamics is effectively Markovian in a space of dimension $k_*$ or less. The linear response with respect to perturbations, $\E^\varepsilon\Phi_n$, is now determined by the properties of the deterministic macroscopic dynamics (\ref{e.PhiF}). In the following sections we shall consider cases and conditions on the deterministic macroscopic dynamics (\ref{e.PhiF}) when $\E^\varepsilon\Phi_n$ enjoys linear response and when it does not.

For finite size $M$ the response approximates that of the thermodynamic limit, as in the case of an uncoupled heat bath discussed in Section~\ref{s.uncoupled}. However, The CLT approximation (\ref{e.Phiback}) assures that the microscopic dynamics (\ref{e.scen2_dn}), which is driven by $d_n=\Phi_n$, is essentially  stochastic with a noise process $\zeta_n$ that has decay of temporal correlations (since the $\qj$ exhibit decay of correlations). This self-generated dynamic noise induces linear response for finite size mean field coupled heat baths, which can be seen as a consequence of results by Hairer and Majda \cite{HairerMajda10}.\\%, with the complication that the parameters determining the self-generated noise process are self-induced. \\

In Section \ref{s.coupled_Triv}-\ref{s.coupled_NonTriv} we will consider the dynamics of the system in the thermodynamic limit. The first case is when $\Phi_n=\langle\E\Phi_n\rangle$ approaches a fixed point $\bar \Phi$ for $M\to \infty$, the second case is when the mean field $\Phi_n$ itself exhibits nontrivial dynamics. Whereas in the first case the linear response of the macroscopic observable $\Psi$ is determined by the properties of the microscopic dynamics, in the latter case it is entirely determined by the response of the macroscopic dynamics.

\subsection{Trivial dynamics of the mean field observable}
\label{s.coupled_Triv}

Let us first look at the case of the mean field at a stable fixed point $\bar\Phi$, in the sense that the mean field remains bounded when perturbed from $\bar\Phi$ and, when the collective dynamics has LRT, $\bar{\Phi}$ is a stable fixed point of $\Phi_n$ in the thermodynamic limit.

To understand the stability, we can apply the external driving framework expounded in the previous section to the dynamics of our system about the equilibrium $d_n \equiv \bar \Phi$.  Stability is in fact assured provided that the complex susceptibility function $R(z)$ does not have any roots inside the unit disk. This follows by considering  $\theta_n= \langle \E\Phi_n \rangle - \bar\Phi = \sum_{k=0}^\infty \chi_{n,k}\theta_{n-k}$ with $\theta_n\sim \lambda^n$ in (\ref{e.dphichi}) which leads to 
\begin{align*}
R(\lambda^{-1})-1=0,
\end{align*}
and hence, for unstable $|\lambda|>1$, to the above condition for instability for the susceptibility function $R(z)$.\\ 

If stability is ensured, the linear response of the fixed point $\bar\Phi $ with respect to external perturbations $\eps g$ is established by the implicit function theorem from the deterministic macroscopic dynamics (\ref{e.PhiF}). Once an external perturbation $\eps g$ is applied, the fixed point depends on $\eps$, and we write
\begin{align}
0 = F(\bar\Phi^\eps,\bar\Phi^\eps,\bar\Phi^\eps,\cdots;\eps) - \bar\Phi^\eps. \label{e.couplingfixedpoint}
\end{align}
In the following numerical experiments the $M\to\infty$ limit is computed by estimating a solution to this algebraic equation.

Differentiation with respect to the external perturbation yields
\begin{align*}
0 &= \frac{d\bar\Phi^\eps}{d\eps} \left( \frac{\partial}{\partial \bar\Phi^\eps} F(\bar\Phi^\eps,\bar\Phi^\eps,\bar\Phi^\eps,\cdots;\eps)-1\right)\\
&\hphantom{=\;} + \frac{\partial }{\partial \eps}F(\bar\Phi^\eps,\bar\Phi^\eps,\bar\Phi^\eps,\cdots;\eps)\\
&=\frac{d\bar\Phi^\eps}{d\eps} \left( \sum_{k=1}^\infty \chi_k -1\right)
+ \frac{\partial }{\partial \eps}F(\bar\Phi^\eps,\bar\Phi^\eps,\bar\Phi^\eps,\cdots;\eps).
\end{align*}
This immediately yields that
\begin{align}
\frac{d\bar\Phi^\eps}{d\eps} = \frac{ \frac{\partial }{\partial \eps}F}{1-R(1)},
\label{e.LRTMinfty}
\end{align}
and hence the existence of linear response, provided $R(1)-1 \neq 0$. 

%We remark that in the large, finite $M$ case one can linearly expand the macroscopic dynamics (\ref{e.PhiF}) about $\Phi_n \equiv \bar{\Phi}$ and from (\ref{e.Phiback}) obtain a CLT approximation:
%\[ \Phi_{n} = \bar\Phi + \frac{1}{\sqrt{M}} \sum_{k=1}^\infty S_k \zeta_{n-k} + \frac{1}{\sqrt{M}}S(1) \eta + o(1/\sqrt{M}),\]
%where $S(z) := 1/(1-R(z))$ and $\zeta_n$ is now a stationary process with 
%\[\Cov[\zeta_n, \zeta_{n-k}] = \langle \Cov^{d \equiv \bar{\Phi}} [\phi(\qj_n),\phi(\qj_{n-k})] \rangle.\]

As for the uncoupled scenario, we shall now discuss the linear response behaviour for the three different cases of the microscopic dynamics, which are covered by the rows in Table~\ref{t.result} corresponding to the coupled macroscopic observables.

%%%%%%%%%%%%%%%%%%%%%%%%%%%%%%%%%%%%%%%%%%%%%%%%%%%%%%%%%%

\subsubsection{The microscopic subsystems satisfy LRT}
\label{s.coupled_TrivA}

\begin{figure}[htbp]
	\centering
	\includegraphics[width=0.9\columnwidth]{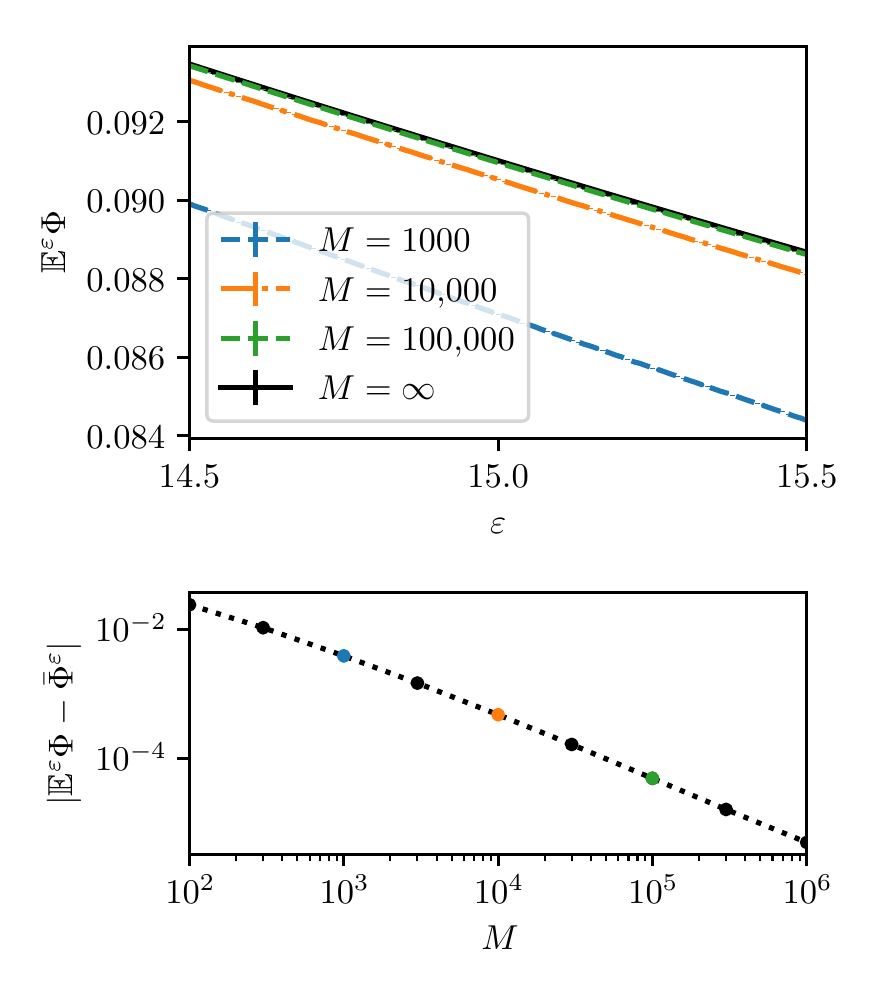}%\;\;
	\caption{(a) Response term $\E^\eps \Phi$ for the uniformly expanding map (\ref{e.Lanford}) with mean field coupling under trivial dynamics; (b) Difference $|\E^\eps \Phi - \bar\Phi^\eps|$ for $\eps = 15$, exhibiting $\O(1/M)$ convergence.
	% of $\E^\eps \Phi$ to thermodynamic limit for $\eps = 15$.% Right: normalised autocovariance function of $\Phi$ for finite $M$ response in trivial dynamics, including $M\to\infty$ limit calculated via (...) and spectral methods \cite{Poltergeist}.
	}
	\label{f.Scenario2_Case1_trivial}
\end{figure}

We consider here the case of uniformly expanding dynamics of the microscopic systems, such that each subsystem individually satisfies LRT. In particular, we choose the following uniformly expanding map %Lanford's map \cite{Lanford98}
\begin{widetext}
\begin{align}
q_{n+1}= \frac{T(q_n) + K_n \left(1 - \sqrt{0.03(1 - 0.97K_n^2) + 0.97(T(q_n) + K_n)^2}\right)}{1 - 0.97K_n^2},%2\qj_n + \qj_{n+1}= \frac{2q\j_n - \sign q\j_n + K_n \left(1 - \sqrt{0.03(1 - 0.97K_n^2) + 0.97(2q\j - \sign q\j_n + K_n)^2}\right)}{1 - 0.97K_n^2},%2\qj_n + \tfrac{1}{2}\qj_n(1-\qj_n) \qquad ({\rm{mod}}\; 1).
\label{e.Lanford}
\end{align}
\end{widetext}
where $K_n = \tanh(\eps \Phi_n - 2)$, $q_n\in [-1,1]$ and $T(q) = 2q - \sign q$ is the doubling map. All microscopic degrees of freedom $\qj$  evolve according to the same map but with randomly distributed initial conditions. (Note that having identical microscopic subsystems implies that $\langle \E^\eps\Phi\rangle=\E^\eps\Phi$.) This map is full-branch uniformly expanding for fixed $K_n$. It is carefully constructed to allow for nontrivial mean field dynamics for larger values of $\eps$ which will be discussed in Section~\ref{s.coupled_NonTriv}. We choose the %an even Lebesgue-measure zero 
coupling function $\phi(q) = -\tfrac{23}{30} + \tfrac{7}{2} q^2 - 2 q^4$ to generate the mean field $\Phi_n$. For simplicity, we choose the mean field observable $\Psi=\Phi$. The dynamics in the thermodynamic limit $M = \infty$ was computed using a spectral method \cite{Wormell19, Poltergeist} (see Appendix \ref{a.poltergeist} for more details). For small values of $\eps$ the mean field $\Phi_n$ converges to a stable fixed point, and the macroscopic observable $\Psi=\Phi$ satisfies LRT as shown in Figure~\ref{f.Scenario2_Case1_trivial}. The variation about $\Phi = \bar{\Phi}$ can be shown to converge to the limiting distribution of the mean zero stochastic process $\zeta_n$ with autocovariance (\ref{e.covzetan_dn}). We will see later in Section~\ref{s.coupled_NonTriv} that for larger values of $\eps$, the mean field exhibits nontrivial chaotic dynamics, violating LRT.

%
% FIGURE
%
%

%%%%%%%%%%%%%%%%%%%%%%%%%%%%%%%%%%%%%%%%%%%%%%%%%%%%%%%%%%

\subsubsection{The microscopic subsystems do not satisfy LRT but are appropriately heterogeneous}
\label{s.coupled_TrivB}

\begin{figure}[htbp]
 	\centering
 	\includegraphics[width=0.9\columnwidth]{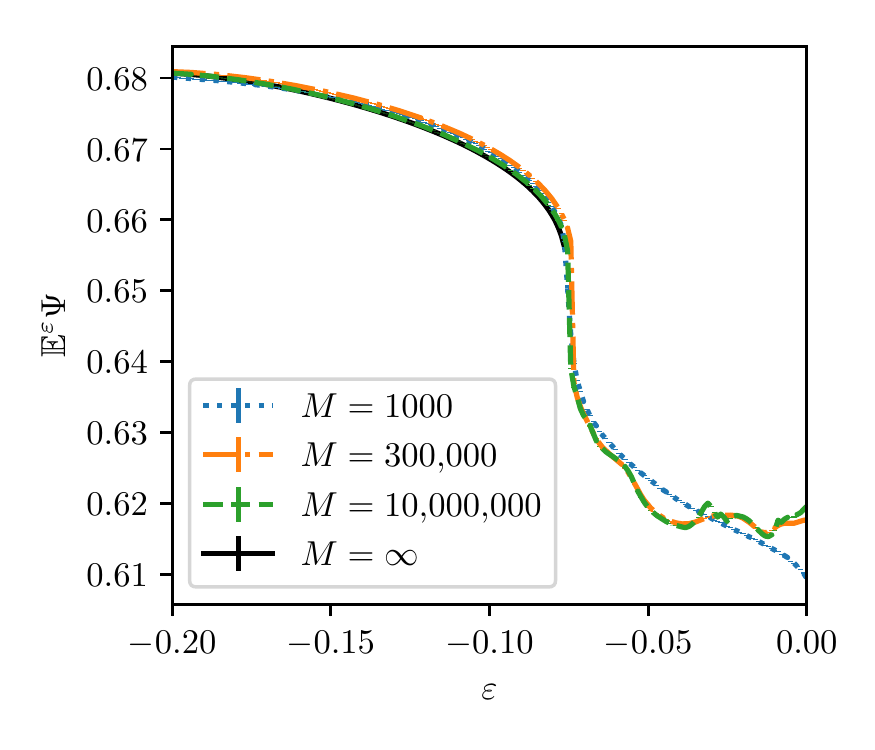}\;\;
 	\caption{Response term $\E^\eps \Psi$ for the modified logistic map (\ref{e.logistic}) with mean field coupling (\ref{e.h}). The parameters $\aj$ are drawn from the raised-cosine distribution (\ref{e.rosine}). Error bars were estimated from $10$ realisations of $10^5$ iterates, differing in the initial conditions of the heat bath, and are not visible.}
 	\label{f.Scenario2_Case2}
 \end{figure}

We consider a mean field coupled system of LRT-violating modified logistic maps (\ref{e.logistic}). We choose the mean field coupling (\ref{e.h}) and draw the parameters $a^{(j)}$ of the logistic map from the smooth raised-cosine distribution (\ref{e.rosine}). In this case, the macroscopic dynamics $\Phi_n$ converges to a stable fixed point $\Phi_n\to \bar\Phi$ for $\eps < -0.075$. The associated linear response is clearly visible in Figure~\ref{f.Scenario2_Case2}. In fact, as discussed in Section~\ref{s.uncoupledB}, nonlinear third order response holds for the three times continuously differentiable raised cosine distribution (\ref{e.rosine}). We remark that for $\eps> -0.075$ the mean field exhibits nontrivial dynamics in the thermodynamic limit, and we observe a breakdown of LRT to be discussed in Section~\ref{s.coupled_NonTriv}.
 
\begin{figure}[htbp]
	\centering
	\includegraphics[width=0.9\columnwidth]{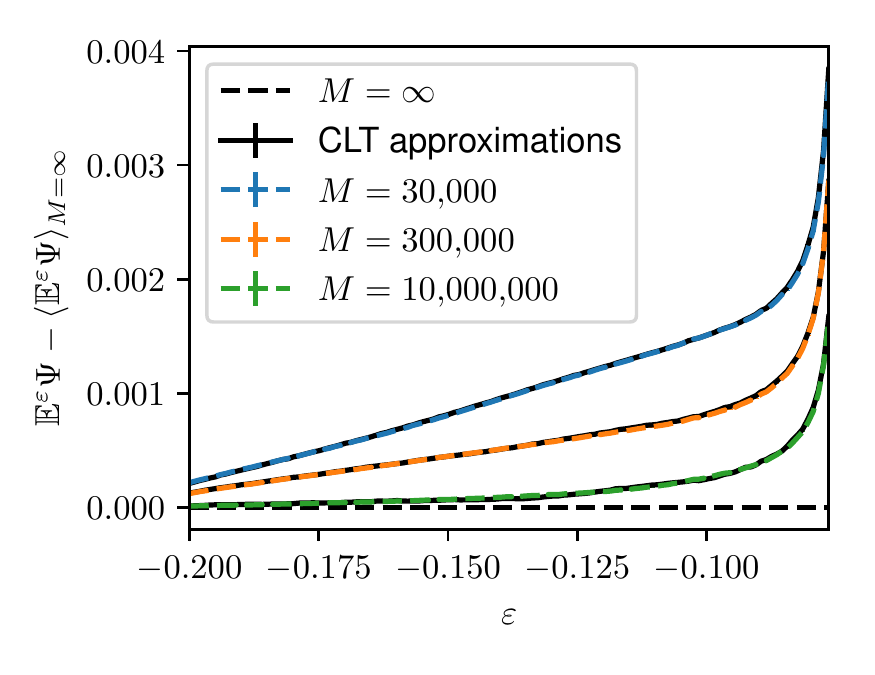}\;\;
	\caption{Difference between the response term $\E^\eps \Psi$ for finite $M$ and for the thermodynamic limit for the modified logistic map (\ref{e.logistic}) with mean field coupling (\ref{e.h}). The parameters $\aj$ of the logistic map are drawn from the raised-cosine distribution (\ref{e.rosine}). For each value of $M$ the response of the corresponding CLT approximation using noise estimated from $M = 10^6$ was used. Error bars were estimated from $10$ realisations of $10^5$ iterates, differing in the initial conditions of the heat bath, and are not visible.}
	\label{f.Scenario2_Case2_diff}
\end{figure}

\begin{figure}[htbp]
	\centering
	\includegraphics[width=0.9\columnwidth]{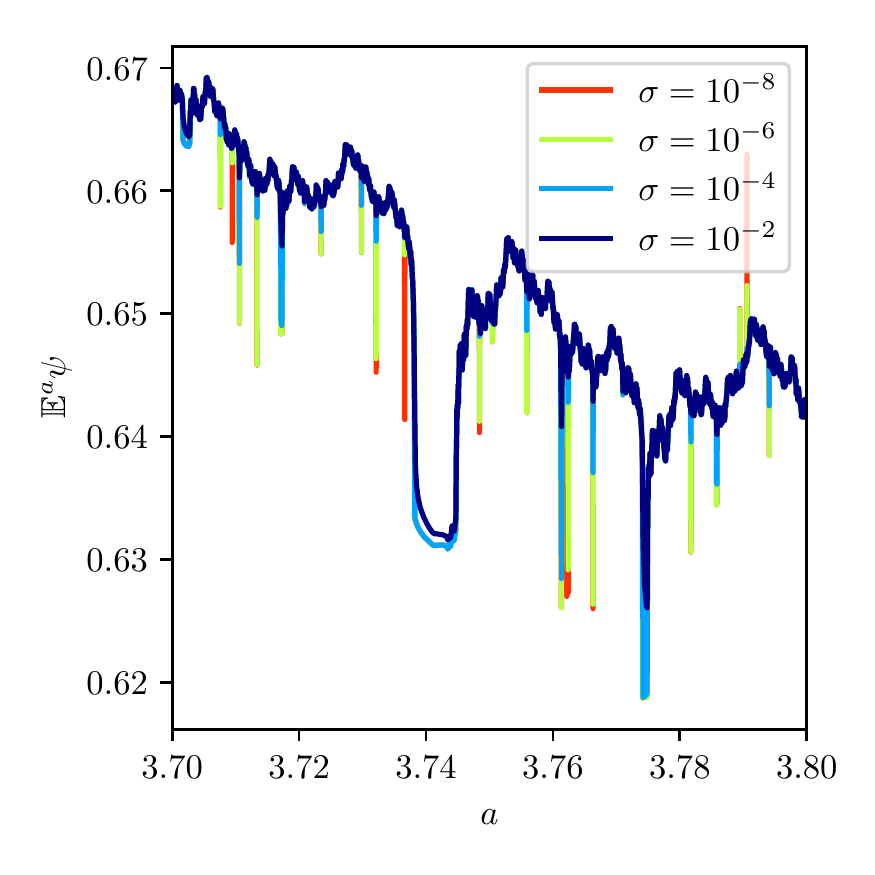}%\;\;
	\caption{Response term $\E^a \psi$ for a single stochastically driven logistic map $q_{n+1} = a q_n (1-q_n) + \sigma \xi_n$, where $\xi_n$ is {\it i.i.d.} Gaussian noise and the observable $\psi(q) = q$. The response is recorded at increments of $da = 10^{-8}$, thus for small $\sigma$ only a subset of narrow periodic windows are captured.}
	\label{f.noisedestruction}
\end{figure}

In Figure~\ref{f.Scenario2_Case2_diff} we see very slow convergence of the mean $\E^\eps \Psi$ to its limiting value: in particular, it is slower than the $\O(1/M)$ rate for uniformly expanding dynamics leading to trivial dynamics (cf Figure~\ref{f.Scenario2_Case1_trivial}), and seemingly slower still than the $\O(1/\sqrt{M})$ rate that we might expect from sampling errors of $\eta$. Although smooth families of microscopic logistic maps allow for linear response of macroscopic observables to {\it constant-in-time perturbations} as discussed in Section \ref{s.uncoupledB}, and in fact numerical experiments (not shown) suggest that the susceptibility function $\chi$ has summable decay, they do not appear to have linear response with respect to stochastic perturbations. We argue that this counter-intuitive lack of linear response with respect to the (self-generated) stochastic perturbations arises from the noise-induced destruction of narrow periodic windows that have ``extreme'' values of $\E^{a,\eps} \psi(q)$ compared with the neighbouring, more stochastically stable chaotic parameters. Thus at these periodic parameter values the macroscopic dynamics exhibits a disproportionately large response to the introduction of noise. We illustrate this in Figure~\ref{f.noisedestruction}, where we plot the response of a single logistic map with additive noise of variance $\sigma^2$. Here the noise models the finite size effects of the heat bath with $\sigma \sim 1/\sqrt{M}$. One sees clearly that periodic windows can be destroyed by very small amounts of noise ($\sigma = 10^{-6}$). One also sees that associated with the destruction of these periodic windows is a very large response in the average $\E^a\psi$.\\

% of ... width and deviation from the chaotic response... are 

%....{\bf{actually not clear what I see here....the red curves are not always more extreme.....no they are, they're just hidden}}.\\

The statistical properties of the macroscopic observable for dynamics of a finite size heat bath can be modelled again by a surrogate system. Writing
\[ 
\Phi_n = \bar \Phi^{\eps,M} + \frac{1}{\sqrt{M}} \zeta_n,
\]
where $\zeta_n$ is the Gaussian CLT correction term to $\Phi_n \equiv \bar \Phi^{\eps,\infty}$ with covariance given by (\ref{e.coupledzetaacv}). The macroscopic dynamics (\ref{e.couplingfixedpoint}) then becomes 
\[ 
0 = \E F(\bar \Phi^{\eps,M} + \frac{1}{\sqrt{M}} \zeta_{n-1}, \bar \Phi^{\eps,M} + \frac{1}{\sqrt{M}} \zeta_{n-2}, \ldots) 
- \bar \Phi^{\eps,M}.
\]
The response $\E^\eps \Psi$ for this surrogate macroscopic dynamics is shown in Figure~\ref{f.Scenario2_Case2_diff}, labelled CLT approximation, and is barely distinguishable from the response of the original macroscopic dynamics.

%{\bf{Something on $R(1)$? We have computed the susceptibility function $R(1)$ and found that $R(1)=???$ for $\eps \in[27.5,30]$. Put in prev section?}}

%
% FIGURE
%

%%%%%%%%%%%%%%%%%%%%%%%%%%%%%%%%%%%%%%%%%%%%%%%%%%%%%%%%%%

\subsubsection{The microscopic subsystems do not satisfy LRT and are not appropriately heterogeneous}
\label{s.coupled_TrivC}

\begin{figure}[htbp]
	\centering
	\includegraphics[width=0.9\columnwidth]{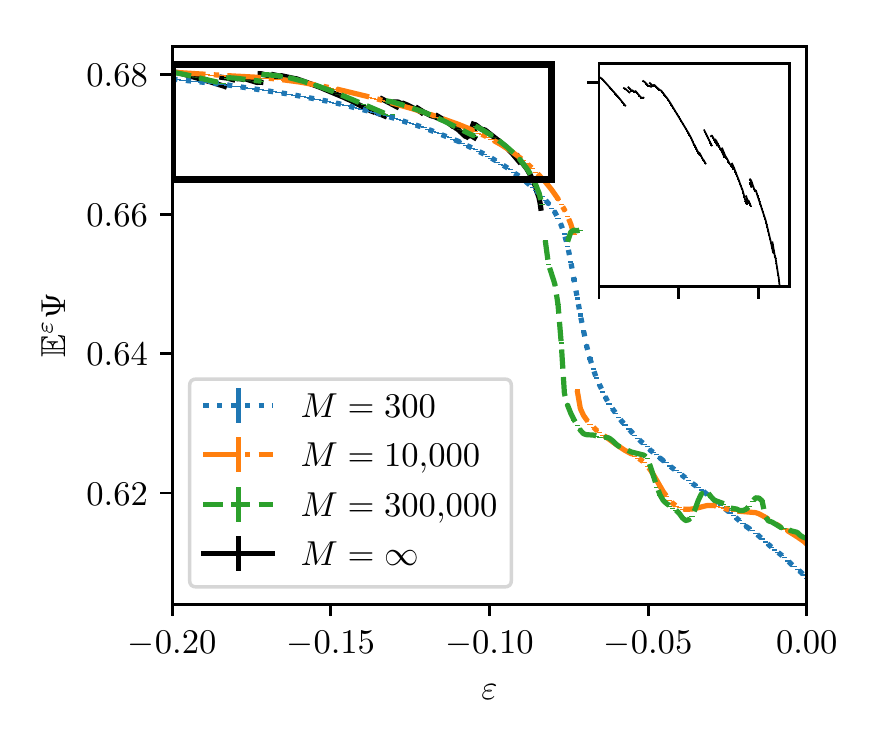}\;\;
	\caption{Response term $\E^\eps \Psi$, including multistability, for the modified logistic map (\ref{e.logistic}) with mean field coupling (\ref{e.h}). The parameters $\aj$ are drawn from the discrete distribution (\ref{e.deltas}). %$\nu = \frac{1}{3} (\delta_{3.72}+\delta_{3.75}+\delta_{3.78})$. 
	Error bars were estimated from $200$ realisations differing in the initial conditions (including initialising distributions) of the heat bath, and are not visible. The inset illustrates the occurrence of saddle-node bifurcations in the infinite-dimensional limit.}
\label{f.Scenario2_Case3}
\end{figure}

For non-smooth distributions of the logistic map parameters $\aj$ such as the discrete distribution (\ref{e.deltas}),  
%$\nu = \frac{1}{3} (\delta_{3.72}+\delta_{3.75}+\delta_{3.78})$ 
the dynamics also converges to a stable fixed point $\bar\Phi^{\eps,M}$ for $\eps<-0.075$. The mean field $\bar\Phi^{\eps,M}$ varies smoothly with respect to $\eps$ for almost all $\eps$, but $\E^\eps\Psi$ experiences saddle node bifurcations on increasingly dense sets as $M$ approaches the thermodynamic limit $M\to\infty$. This is illustrated in Figure~\ref{f.Scenario2_Case3}. Looking at (\ref{e.LRTMinfty}), we see that linear response is violated where the fixed point loses stability and $\partial F(\bar\Phi^\eps,\bar\Phi^\eps,\cdots)/\partial \bar\Phi^\eps=R(1)=1$. 
%The breakdown of LRT as $M\to\infty$ is linked to growth in $\partial F(\bar\Phi^\eps,\bar\Phi^\eps,\cdots)/\partial \bar\Phi^\eps$. 

For finite $M$, the macroscopic equation can be modelled as $\Phi_n^{\eps,M}=\bar\Phi^{\eps,M} + \zeta_n/\sqrt{M}$ where
\[  \bar{\Phi}^{\eps,M} = \E^\eps F(\Phi_n^{\eps,M},\Phi_n^{\eps,M},\cdots). \] 
The derivatives $\partial \E^\eps F(\Phi_n^\eps,\Phi_n^\eps,\cdots)/\partial \bar\Phi^\eps$ are well defined and thus for appropriate bath sizes $M$ we observe again approximate LRT for all practical purposes. However, as $M\to \infty$, saddle-node bifurcations become visible as a result of the diminishing effect of stability-providing noise, and in the thermodynamic limit we observe failure of LRT (see inset in Figure~\ref{f.Scenario2_Case3}). The failure of linear response through saddle-node bifurcations is accompanied by $\E^\eps\Psi$ experiencing multistability, with multiple very close stable equilibria, demonstrated in Figure~\ref{f.Scenario2_Case3}\footnote{When randomly searching for equilibria it is important to make sure that, as well as randomly initialising the $q\j_0$, the {\it distribution} from which the $q\j_0$ are sampled is also randomly initialised, as up to an error term of $\O(M^{-1/2})$ the macroscopic dynamics are deterministic functions of the initial measures of the microscopic variables $\mu^a_0$.}. This can be understood as coming from the fact that $\E^\eps F(\Phi_n^{\eps,M},\Phi_n^{\eps,M},\cdots)$ is essentially a smoothed out version of the rough logistic map response $\E^\eps F(\bar\Phi,\bar\Phi,\cdots)$: as $M \to \infty$ the smoothing decreases, leading to increasing numbers of roots of the equation.

 %For finite heat baths, the finite-size induced noise causes the dynamics to hop between the stable equilibria and, if sufficiently close, smears out the multistability, leading to approximate LRT. 

%The poor behaviour of the derivative corresponds to the breakdown of LRT in the limiting systems as the noise is decreased to zero.
 
%{\bf{more; relate to "proof" above. You wrote: 
%"$\partial F(\bar\Phi^\eps,\bar\Phi^\eps,?)/\partial \bar\Phi^\eps$ is not defined, but $\partial F(\bar\Phi^\eps + \zeta_n/\sqrt{M},\bar\Phi^\eps+ \zeta_{n-1}/\sqrt{M},?)/\partial \bar\Phi^\eps$ is well-defined, at least in some annealed sense. The problem is that as $M \to \infty$ this derivative tends to get bigger, which makes the fixed point more unstable, inducing more bifurcations." But here we are only considering the $M=\infty$ case; the finite $M$ case is dealt with above.}}

%
% FIGURE
%
%

%	* Actually, we could think about even a very simplified version of the tangent space dynamics:
%	\[\rho \mapsto \mathcal{L} \rho + \nu \int \phi \rho dx, \]
%	where $\mathcal{L}$ is a transfer operator, $\int \nu dx = 0$, and suppose the leading eigenvalue of this operator is of norm $>1$. You can separate out different eigenspaces using the left eigenfunctions of the operator, which except for the 1-eigenvalue (which we aren't considering coz tangent space) are all in funky hyperdistribution spaces. Can we then get a nice splitting so $|(a,b)| \leq C^{-1} \|a\| \|b\|$ for $C<1$ and $b$ in the kernel of the left eigenfunction and $a$ a right eigenfunction?

%%%%%%%%%%%%%%%%%%%%%%%%%%%%%%%%%%%%%%%%%%%%%%%%%%%%%%%%%%

\subsection{Nontrivial dynamics of the mean field observable}
\label{s.coupled_NonTriv}

The mean field $\Phi$ or any macroscopic observable $\Psi$ may itself exhibit non-trivial dynamics of varying complexity in the thermodynamic limit $M=\infty$. The overall response behaviour is then determined by the macroscopic dynamics rather than by the properties of the microscopic subsystems. We show the emergence of non-trivial chaotic macroscopic dynamics which violates LRT. The first one, surprisingly, involves a heat bath which evolves under uniformly expanding dynamics when uncoupled, and the second one involves microscopic dynamics that individually violate LRT.

%%%%%%%%%%%%%%%%%%%%%%%%%%%%%%%%%%%%%%%%%%%%%%%%%%%%%%%%%%

%\subsubsection{The microscopic subsystems satisfy LRT}
%\label{s.coupled_NonTrivA}

\begin{figure}[htbp]
\centering
\includegraphics[width=0.77\columnwidth]{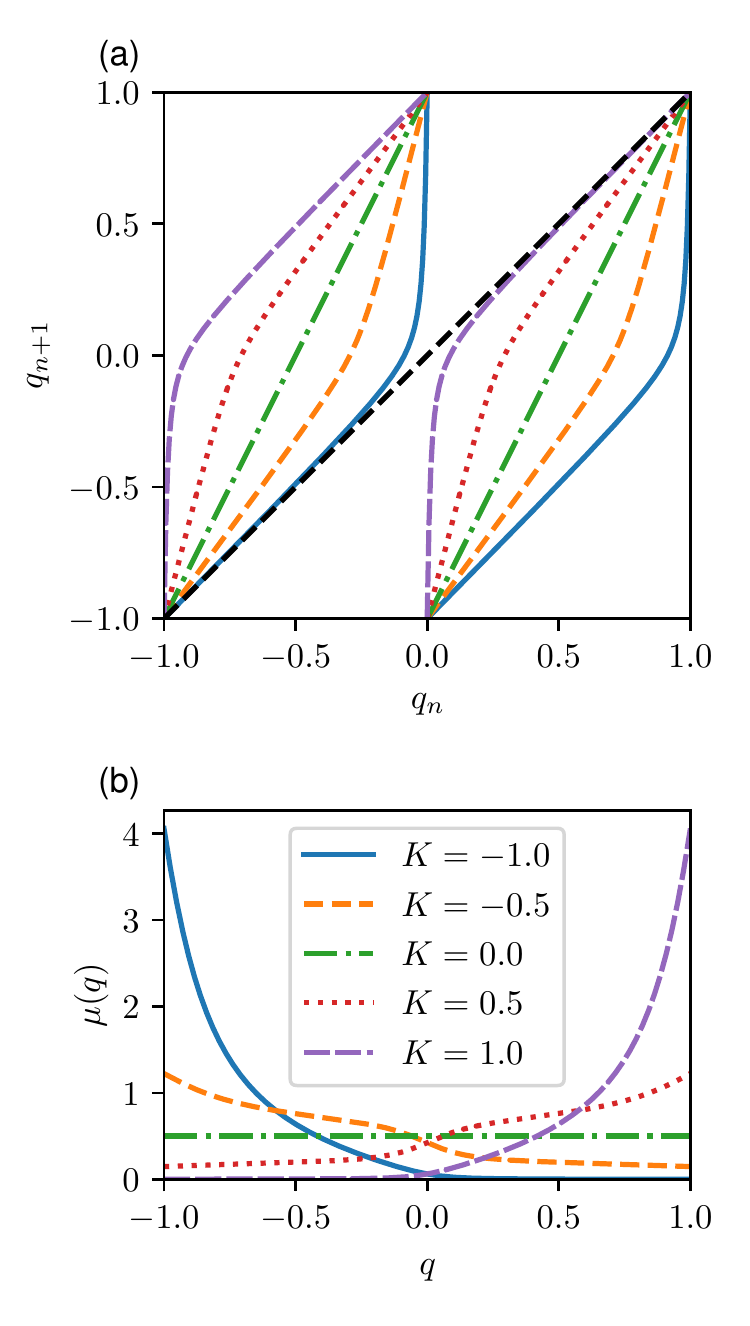}%\\
%\includegraphics[width=0.97\columnwidth]{figures/cestnepasunmapLanford2.pdf}
%(a) $\qquad\qquad\qquad\qquad\qquad\qquad\qquad\qquad$ (b)
\caption{(a) Plot of the uniformly expanding map (\ref{e.Lanford}) and (b) its invariant measure for various values of $K_n \equiv K$.}
\label{f.cestnepasunmapLanford}
\end{figure}

To generate emergent nontrivial macroscopic dynamics of the mean field, we again use the uniformly expanding map (\ref{e.Lanford}) with the even Lebesgue-measure zero coupling function $\phi(q) = -\tfrac{23}{30} + \tfrac{7}{2} q^2 - 2 q^4$ to generate the mean field $\Phi_n$. We show in Figure~\ref{f.cestnepasunmapLanford} the map and its invariant measure, where the dynamics in the thermodynamic limit $M = \infty$ was computed using a spectral method \cite{Wormell19, Poltergeist} (see Appendix \ref{a.poltergeist} for more details). The map and coupling function $K_n$ are judiciously chosen to yield nontrivial dynamics for the mean field $\Phi_n$, mediating dynamics akin to a unimodal map for $\Phi_n$. 

The map is constructed such that when the $\qj$ are approximately evenly distributed, $\Phi_n\approx 0$, causing an extreme value $K_{n} \approx \tanh(-2) \approx -0.96$: this pushes the $\qj$ strongly towards $q = -1$ which leads to a larger value $\Phi_{n+1}$, concentrating around $\Phi=2/\eps$. For these values of $\Phi_{n+1}$, $K_{n+1}\approx 0$, and thus in the next step the $\qj$ are spread more evenly over the interval $[-1,1]$, mapping $\Phi_{n+1}$ back around zero. The concentration in the first step provides the folding and the sensitivity of $K_n$ to small changes in $\Phi_n$ for large $\eps$ provides the stretching necessary for chaotic dynamics. 

\begin{figure}[htbp]
\centering
\includegraphics[width=0.97\columnwidth]{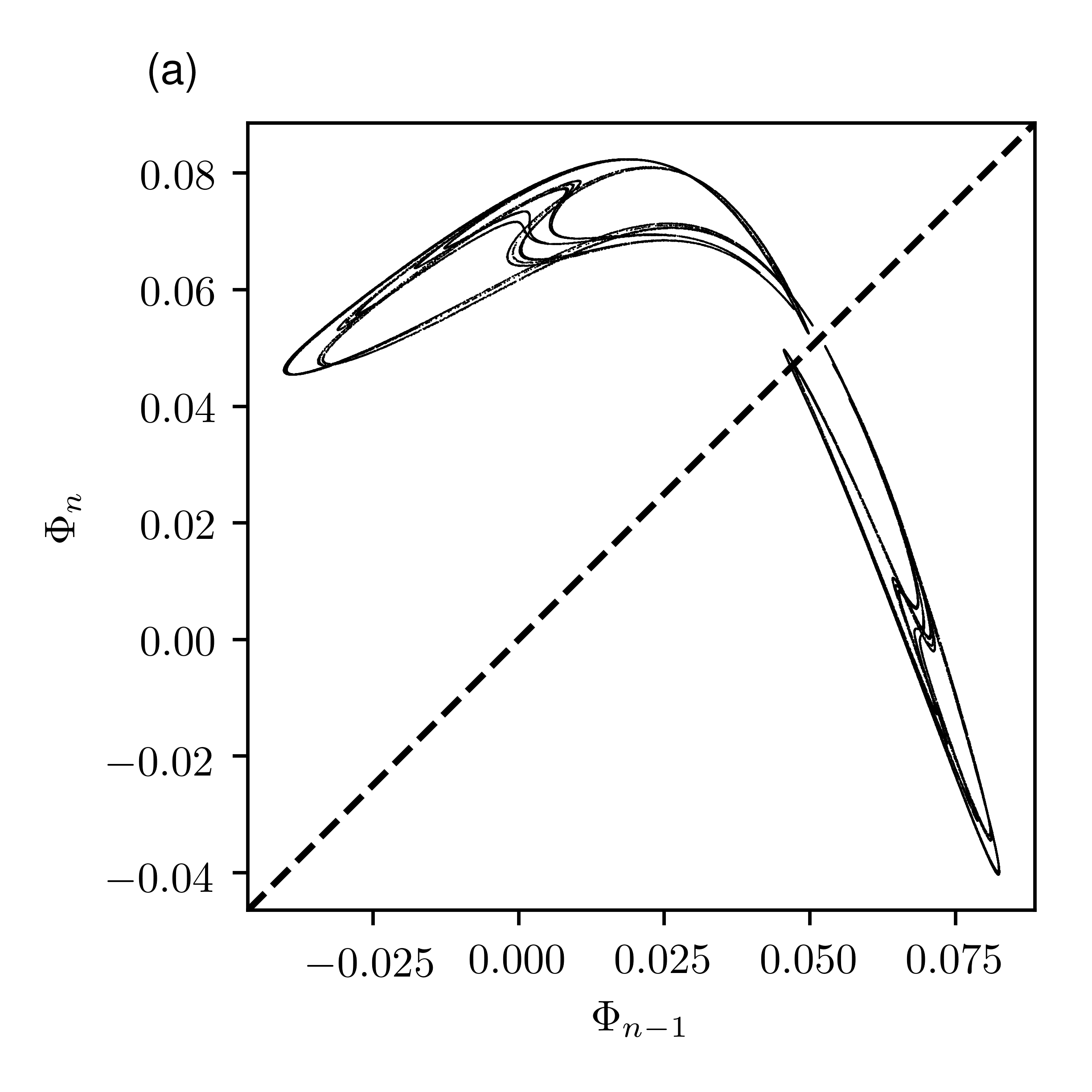}\\
\includegraphics[width=0.97\columnwidth]{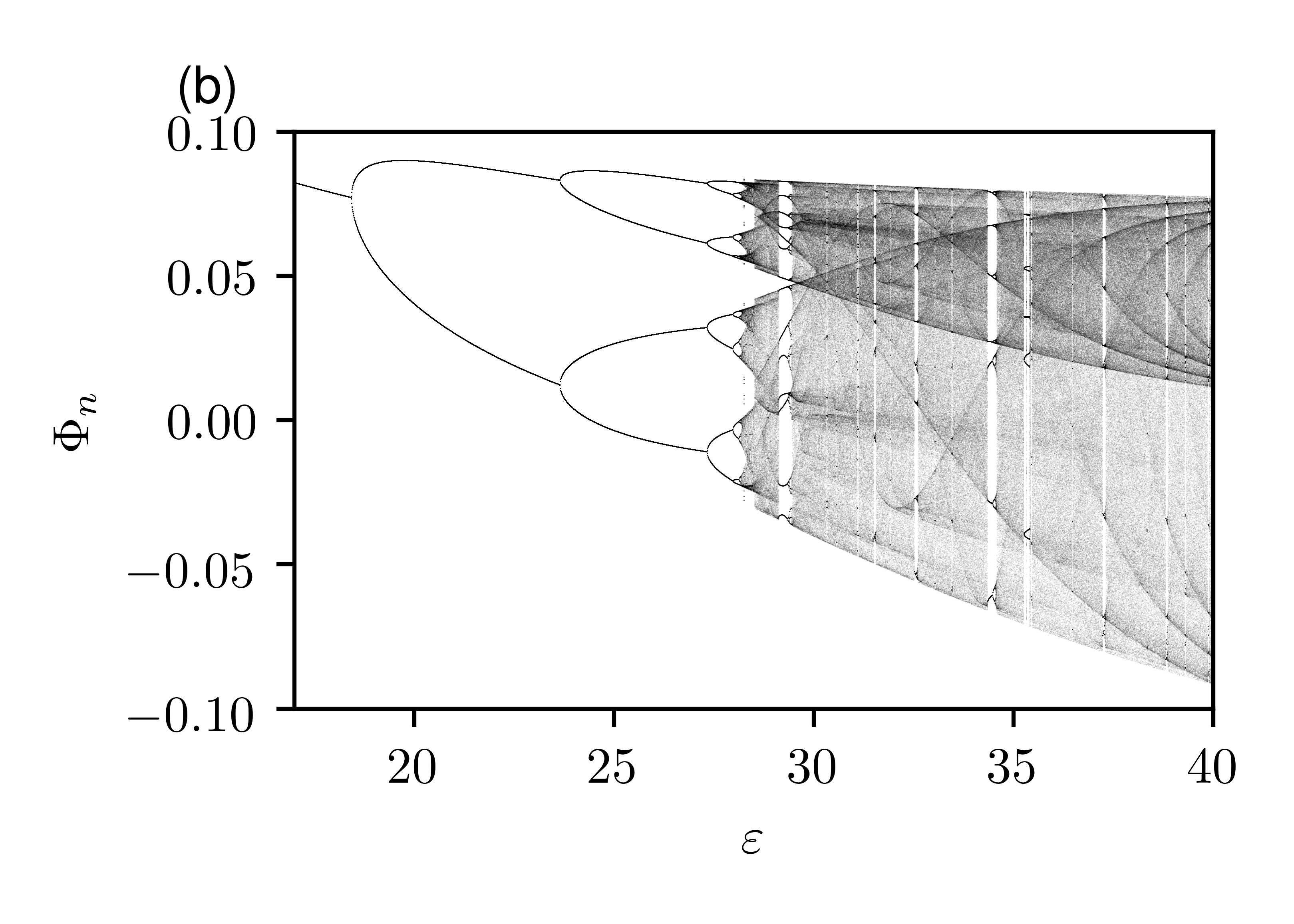}
\caption{(a): 2D projection of the attractor onto delay coordinates of the macroscopic map $\Phi_{n+1}=F(\Phi_n,\Phi_{n-1},\ldots)$ 
	generated by the uniformly expanding map (\ref{e.Lanford}) for $\eps = 30$. The system has two periodic components separated by a gap around the unstable fixed point $\Phi_{n-k} \equiv 0.51258$. (b): Bifurcation diagram of the map (\ref{e.Lanford}) showing period doubling bifurcations and chaotic dynamics.}
\label{f.PhinoLanford}
\end{figure}

In Figure~\ref{f.PhinoLanford} we show the map $\Phi_{n+1}=F(\Phi_n, \Phi_{n-1}, \ldots)$ generated by the dynamical system (\ref{e.Lanford}) in the thermodynamic limit $M = \infty$ for $\eps=30$. The dynamics is clearly chaotic with the leading Lyapunov exponent $\lambda_1=0.18>0$ ($\lambda_2=-0.43$ and $\lambda_3=-0.81$). The dynamics of the macroscopic observable $\Psi_n=\Phi_n$ exhibits a complex bifurcation cascade upon varying $\eps$, depicted in Figure~\ref{f.PhinoLanford}. For $\eps \leq 18.4159$, the macroscopic dynamics has a stable fixed point; upon increasing the perturbation $\eps$ a period-doubling cascade leads to chaotic, apparently unimodal-like, dynamics intermingled with periodic windows for values of $\eps> 26.1649$. One can clearly see dark scars in the bifurcation diagram in the chaotic region of $\eps>26.1649$. This is reminiscent of the logistic map \cite{ColletEckmann07} where the scars denote narrow intervals of $\langle\Phi_n \rangle$ with increased probability, corresponding to large spikes in the invariant measure, which (unlike small spikes) vary smoothly with respect to perturbations $\eps$. \\

\begin{figure}[htbp]
	\centering
	\includegraphics[width=0.97\columnwidth]{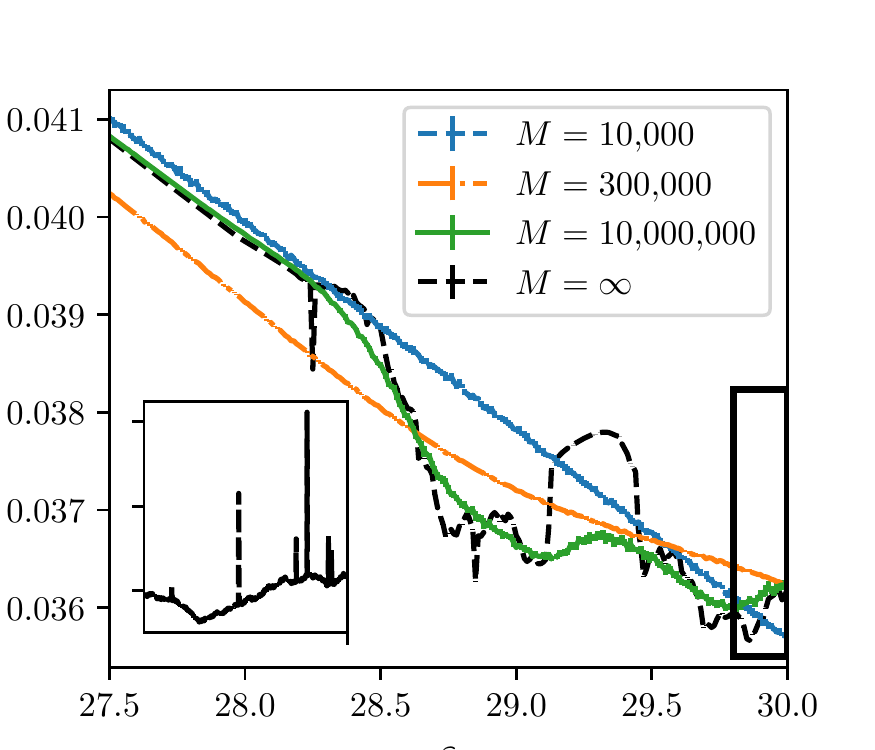}
	%	\includegraphics[width=0.45\columnwidth]{figures/cestnepasunmapLanford2.pdf}\\
	%	(a) $\qquad\qquad\qquad\qquad\qquad\qquad\qquad\qquad$ (b)
	\caption{Response term $\E^\eps \Psi$ of the uniformly expanding map (\ref{e.Lanford}) for finite $M$ response, showing convergence to the thermodynamic limit $M\to\infty$. The black box shows the region which is magnified on the right. Error bars were estimated from $10$ realisations differing in the initial conditions of the heat bath, and are not visible.}
	\label{f.Scenario2_Case1}
\end{figure}

In Figure~\ref{f.Scenario2_Case1} we show the linear response term $\E^\eps \Psi$ of the uniformly expanding map (\ref{e.Lanford}) for several finite $M$ heat bath sizes and for the thermodynamic limit $M=\infty$ for $\eps \in[27.5,30]$, clearly illustrating the breakdown of LRT. We recall that the same map exhibits LRT for small values of $\eps$, where the macroscopic mean field converges to a stable fixed point, for the same parameters (cf. Figure~\ref{f.Scenario2_Case1_trivial}).\\

%
% FIGURE
%
%

%
% FIGURE
%
%

%
% FIGURE
%
%

%%%%%%%%%%%%%%%%%%%%%%%%%%%%%%%%%%%%%%%%%%%%%%%%%%%%%%%%%%

%\subsubsection{The microscopic subsystems do not satisfy LRT but are appropriately heterogeneous}
%\label{s.coupled_NonTrivB}

\begin{figure}[htbp]
	\centering
	\includegraphics[width=0.97\columnwidth]{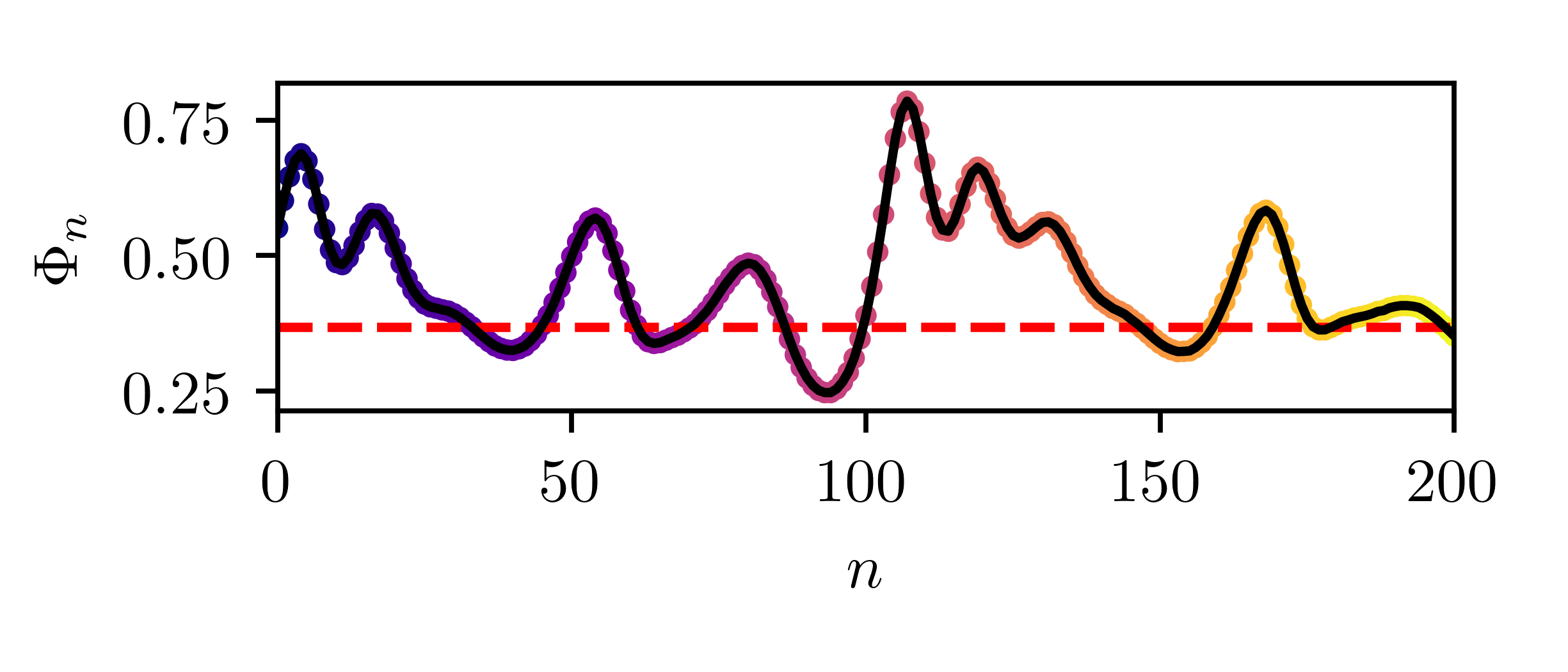}\\
	\includegraphics[width=0.97\columnwidth]{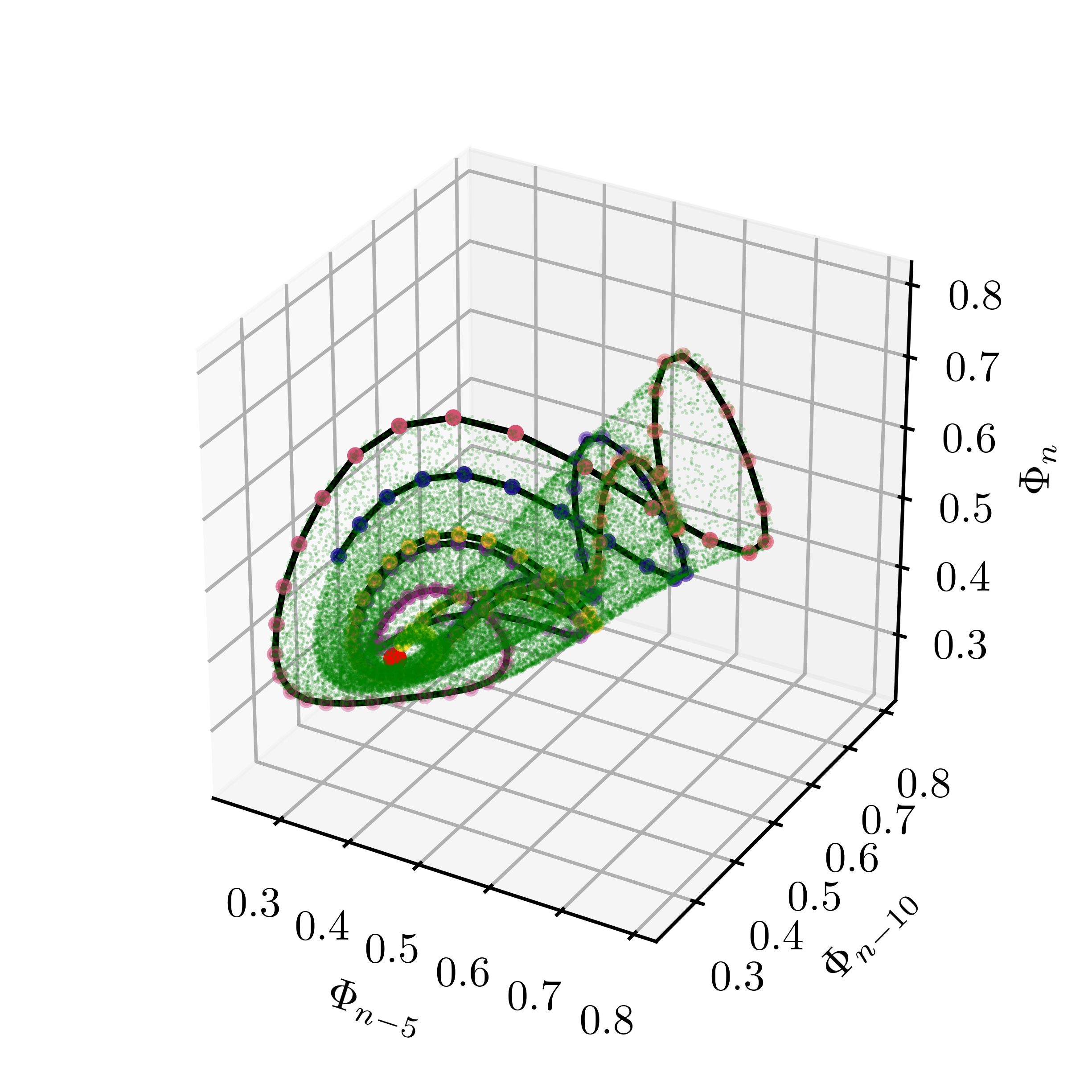}
	\caption{Top: Time series of the macroscopic map $\Phi_{n+1} = F(\Phi_n, \Phi_{n-1},\ldots)$ generated by the modified logistic map (\ref{e.logistic}) with mean field coupling (\ref{e.h})  for $\eps=0$, approximated by a finite ensemble of size $M = 10^7$. The red dotted line shows an unstable fixed point of the system. The parameters $\aj$ are drawn from the raised-cosine distribution (\ref{e.rosine}). Bottom: Projection onto delay coordinates of the attractor and dynamics of the same map. The red dot near the centre of the attractor denotes an unstable fixed point of the system.
	}
\label{f.Shilnikov}
\end{figure}

For curiosity and to further study the effect of the self-generated noise on the LRT behaviour of macroscopic observables in the mean field coupled case, we provide another example of nontrivial chaotic mean field dynamics which violates LRT. We revisit  the mean field coupled dynamics of  microscopic subsystems which do not satisfy LRT discussed in Section~\ref{s.coupled_TrivB}, and consider the modified logistic map (\ref{e.logistic}) with mean field coupling (\ref{e.h}) where the parameters  $\aj$ of the logistic map are drawn from the smooth raised-cosine distribution (\ref{e.rosine}). We recall that for $\varepsilon \approx [-0.2,-0.075]$ the macroscopic dynamics (\ref{e.PhiF}) was trivial and $\E^\varepsilon\psi$ satisfies LRT (cf. Figure~\ref{f.Scenario2_Case2}). The stable fixed point loses stability at $\eps\approx -0.075$ through a saddle-node bifurcation (not shown), from which emanates a stable limit cycle centred around an unstable fixed point with $\E \Psi \approx 0.615$: as $\eps$ increases, this bifurcates to chaos and for a wide range of values $\eps>-0.075$ nontrivial chaotic macroscopic dynamics of (\ref{e.PhiF}) is observed. Figure~\ref{f.Shilnikov} illustrates the macroscopic dynamics for $\eps=0$, which exhibit Shilnikov-type chaos. The associated response term $\E^\eps \Psi$ was shown in Figure~\ref{f.Scenario2_Case2} for several finite heat bath sizes $M$ and for the thermodynamic limit $M=\infty$ for $\eps \in[-0.2,0]$, clearly illustrating the transition to LRT violating macroscopic dynamics around $\eps=-0.075$. Note that the finite size response is smoothed due to the self-generated noise process $\zeta_n$.\\

%Note that for finite $M$ the response $\E^\eps \Psi$ as a function of $\eps$ is typically a smoothed version of that of the thermodynamic limit and the dense set of perturbation sizes inducing bifurcations are smeared out by the self-generated noise; this effect of smearing of bifurcations is also visible in Figure~\ref{f.Scenario2_Case2} for $\varepsilon >-0.075$ corresponding to nontrivial chaotic macroscopic dynamics.\\

The examples given above of high-dimensional system exhibiting non-uniformly hyperbolic chaotic collective behaviour are in disagreement with the often invoked assumption that macroscopic observables of high-dimensional systems obey linear response. 
%Gallavotti-Cohen hypothesis \cite{GallavottiCohen95a,GallavottiCohen95b, Gallavotti19} according to which the attracting dynamics of a high-dimensional system behaves for all practical purposes as an Anosov system. 
This is the more surprising as the non-uniformly hyperbolic chaotic behaviour is robust (modulo periodic windows) with respect to the external perturbation, different choices of the coupling function, different weightings in the coupling, etc.

%
% FIGURE
%
%

%%%%%%%%%%%%%%%%%%%%%%%%%%%%%%%%%%%%%%%%%%%%%%%%%%%%%%%%%%

\section{Discussion}
\label{s.discussion}

We established conditions under which macroscopic mean field observables enjoy linear response. We considered two scenarios, macroscopic observables of an uncoupled collection of microscopic subunits as well as macroscopic observables of  microscopic subunits which are coupled via their mean field. We found that linear response is possible even in the case when the microscopic systems individually violate LRT, provided the microscopic dynamics is heterogeneous with parameters drawn from a sufficiently smooth distribution. We also found that for back-coupled systems of finite size, LRT (for small enough perturbations) is expected for any kind of microscopic dynamics: this can be understood as the result of emergent, self-generated stochastic effects.  We further established that in the thermodynamic limit of infinite $M$, the mean field dynamics can exhibit attracting dynamics that appears non-uniformly hyperbolic and certainly fails to have LRT, even when the microscopic subunits are individually uniformly expanding; this presents a counter example to the widely believed hypothesis that macroscopic observables of high-dimensional systems typically obey linear response.  
%Gallavotti-Cohen hypothesis' prediction of Anosov dynamics in this situation \cite{Gallavotti19}. 
%{\bf At the conference Gallavotti said that because it the attractor is a $\sim$1D manifold then actually logistic, because most "Anosov", is the best you can do but I'm not sure about this: the Lorenz-63 induced map has a "1D" attractor and is actually uniformly hyperbolic. Presumably he realised this since.}

Our results rely on the existence of statistical limit laws such as the central limit theorem. These are proved for strongly chaotic systems, and in particular for uniformly expanding maps as well as for smooth unimodal maps. We follow here \citet{GottwaldMelbourne14} and assume that typical dynamical systems are strongly chaotic and hence enjoy good statistical properties, so that our results carry over to typical dynamical microscopic systems. \Rone{To ensure the existence of the CLT and the convergence of the deviations $\zeta_n$ to a Gaussian process with decay of correlations, we require the observables $\psi(q)$ to be at least H\"older continuous and the external forcing $\eps$ to be such that the perturbed system is mixing.}

We presented here results for mean field observables $\Psi$ of the form (\ref{e.Psi}). We remark that our results carry over for more general (e.g. weighted) mean field variables provided those weights are sufficiently smoothly distributed, and indeed we expect broadly similar results for more general ``macroscopic'' observables.\\

%We have worked under the assumption that the microscopic chaotic systems are all deterministically de-synchronised from each other (e.g. by using the hidden dynamics in (\ref{e.logistic}). Will probably get different, more complex, behaviours if they can synchronise. \\

In previous work on LRT in high-dimensional systems we considered the more specific case where $\eps$ was an additive perturbation of the logistic parameters $\aj$\cite{WormellGottwald18}. The case of homogeneous additive perturbations can be applied to the macroscopic observables treated here. This is readily seen, e.g. in the uncoupled case, by writing
\begin{align*}
	\langle \E^\eps \Psi_n \rangle &= \iint \psi(q) d\mu^{{a+\eps}}_n(q) d\nu(a)\\
	&= \iint \psi(q) d\mu^{{a}}_n(q) \nu(a-\eps)da.
\end{align*} 
The linear response term is then readily evaluated as
\begin{align*}
	\frac{d}{d\eps}\langle \E^\eps \Psi_n \rangle &= -\iint \psi(q) d\mu^{{a}}_n(q) \frac{d}{da}\nu(a)da,
\end{align*} 
which implies that LRT is valid provided that the system is appropriately heterogeneous with integrable distribution $d\nu(a)/da$.\\ %For global perturbations of the dynamics (\ref{e.scen2}) the linear response is more intricate as it involves the combined effect of the distribution of the parameters and the dynamics.

To reduce the complexity of expression we have enforced mixing dynamics, with no chaotic synchronisation, for example by including the hidden $r$-dynamics in (\ref{e.logistic}). It would be interesting to study the case when the microscopic dynamics is not restricted in this way, for example if periodic dynamics were allowed. We have only discussed the existence of LRT and have not considered fluctuation-dissipation formulae to provide a compact analytical formula for the response term. This may %, in particular, require analysis of the smoothness of the invariant measure with respect to the macroscopic state variable, as well as 
require treatment of the non-Markovian dynamics of the macroscopic variable as well as the interplay of the perturbed microscopic dynamics and the macroscopic dynamics, the latter having been studied in the context of slow-fast systems\cite{Abramov10}.\\

%We expect that %over exponentially long (in $M$) timescales, 
%periodic subsystems of a given length will synchronise with each other, however this synchronisation may occur over an exponentially long (in $M$) timescale, for example in trivial macroscopic dynamics. This would lead to two thermodynamic limits: a ``short-time'' limit where periodic subsystems are potentially desynchronised, and a ``long-time'' limit where they are synchronised.
%\\

We have corroborated our findings with detailed numerical simulations and have provided several heuristic arguments based on statistical limit laws: we hope these arguments can be made rigorous. \Rone{We remark that in the case of uncoupled systems, rather than averaging over the heat bath, one may express the invariant measures $\mu^{(a_{j}}$ as an infinite sum over the systems' unstable periodic orbits (UPOs) \cite{Ruelle04, Pollicott86,CvitanovicEckhardt91,EckhardtGrossmann94}. This approach may be more effective when trying to rigorously justify the convergence of the deviations $\zeta_n$ to a Gaussian process.} It seems likely also possible to apply ideas from a recent argument by Ruelle for linear response in non-hyperbolic systems, based on the statistical smearing-out of singularities in the physical measures\cite{Ruelle18}.

%%%%%%%%%%%%%%%%%%%%%%%%%%%%%%%%%%%%%%%%%%%%%%%%%%%%%%%%%%

\begin{acknowledgments}
GAG acknowledges support from the Australian Research Council, Grant No. DP180101385.
\end{acknowledgments}

\bibliography{bibliography}

\begin{thebibliography}{65}%
\makeatletter
\providecommand \@ifxundefined [1]{%
 \@ifx{#1\undefined}
}%
\providecommand \@ifnum [1]{%
 \ifnum #1\expandafter \@firstoftwo
 \else \expandafter \@secondoftwo
 \fi
}%
\providecommand \@ifx [1]{%
 \ifx #1\expandafter \@firstoftwo
 \else \expandafter \@secondoftwo
 \fi
}%
\providecommand \natexlab [1]{#1}%
\providecommand \enquote  [1]{``#1''}%
\providecommand \bibnamefont  [1]{#1}%
\providecommand \bibfnamefont [1]{#1}%
\providecommand \citenamefont [1]{#1}%
\providecommand \href@noop [0]{\@secondoftwo}%
\providecommand \href [0]{\begingroup \@sanitize@url \@href}%
\providecommand \@href[1]{\@@startlink{#1}\@@href}%
\providecommand \@@href[1]{\endgroup#1\@@endlink}%
\providecommand \@sanitize@url [0]{\catcode `\\12\catcode `\$12\catcode
  `\&12\catcode `\#12\catcode `\^12\catcode `\_12\catcode `\%12\relax}%
\providecommand \@@startlink[1]{}%
\providecommand \@@endlink[0]{}%
\providecommand \url  [0]{\begingroup\@sanitize@url \@url }%
\providecommand \@url [1]{\endgroup\@href {#1}{\urlprefix }}%
\providecommand \urlprefix  [0]{URL }%
\providecommand \Eprint [0]{\href }%
\providecommand \doibase [0]{http://dx.doi.org/}%
\providecommand \selectlanguage [0]{\@gobble}%
\providecommand \bibinfo  [0]{\@secondoftwo}%
\providecommand \bibfield  [0]{\@secondoftwo}%
\providecommand \translation [1]{[#1]}%
\providecommand \BibitemOpen [0]{}%
\providecommand \bibitemStop [0]{}%
\providecommand \bibitemNoStop [0]{.\EOS\space}%
\providecommand \EOS [0]{\spacefactor3000\relax}%
\providecommand \BibitemShut  [1]{\csname bibitem#1\endcsname}%
\let\auto@bib@innerbib\@empty
%</preamble>
\bibitem [{\citenamefont {Majda}, \citenamefont {Abramov},\ and\ \citenamefont
  {Gershgorin}(2010)}]{MajdaEtAl10}%
  \BibitemOpen
  \bibfield  {author} {\bibinfo {author} {\bibfnamefont {A.~J.}\ \bibnamefont
  {Majda}}, \bibinfo {author} {\bibfnamefont {R.}~\bibnamefont {Abramov}}, \
  and\ \bibinfo {author} {\bibfnamefont {B.}~\bibnamefont {Gershgorin}},\
  }\bibfield  {title} {\enquote {\bibinfo {title} {High skill in low-frequency
  climate response through fluctuation dissipation theorems despite structural
  instability},}\ }\href {\doibase 10.1073/pnas.0912997107} {\bibfield
  {journal} {\bibinfo  {journal} {Proceedings of the National Academy of
  Sciences}\ }\textbf {\bibinfo {volume} {107}},\ \bibinfo {pages} {581--586}
  (\bibinfo {year} {2010})}\BibitemShut {NoStop}%
\bibitem [{\citenamefont {Lucarini}\ and\ \citenamefont
  {Sarno}(2011)}]{LucariniSarno11}%
  \BibitemOpen
  \bibfield  {author} {\bibinfo {author} {\bibfnamefont {V.}~\bibnamefont
  {Lucarini}}\ and\ \bibinfo {author} {\bibfnamefont {S.}~\bibnamefont
  {Sarno}},\ }\bibfield  {title} {\enquote {\bibinfo {title} {A statistical
  mechanical approach for the computation of the climatic response to general
  forcings},}\ }\href {\doibase 10.5194/npg-18-7-2011} {\bibfield  {journal}
  {\bibinfo  {journal} {Nonlinear Processes in Geophysics}\ }\textbf {\bibinfo
  {volume} {18}},\ \bibinfo {pages} {7--28} (\bibinfo {year}
  {2011})}\BibitemShut {NoStop}%
\bibitem [{\citenamefont {Abramov}\ and\ \citenamefont
  {Majda}(2007)}]{AbramovMajda07}%
  \BibitemOpen
  \bibfield  {author} {\bibinfo {author} {\bibfnamefont {R.~V.}\ \bibnamefont
  {Abramov}}\ and\ \bibinfo {author} {\bibfnamefont {A.~J.}\ \bibnamefont
  {Majda}},\ }\bibfield  {title} {\enquote {\bibinfo {title} {Blended response
  algorithms for linear fluctuation-dissipation for complex nonlinear dynamical
  systems},}\ }\href {http://stacks.iop.org/0951-7715/20/i=12/a=004} {\bibfield
   {journal} {\bibinfo  {journal} {Nonlinearity}\ }\textbf {\bibinfo {volume}
  {20}},\ \bibinfo {pages} {2793} (\bibinfo {year} {2007})}\BibitemShut
  {NoStop}%
\bibitem [{\citenamefont {Abramov}\ and\ \citenamefont
  {Majda}(2008)}]{AbramovMajda08}%
  \BibitemOpen
  \bibfield  {author} {\bibinfo {author} {\bibfnamefont {R.~V.}\ \bibnamefont
  {Abramov}}\ and\ \bibinfo {author} {\bibfnamefont {A.~J.}\ \bibnamefont
  {Majda}},\ }\bibfield  {title} {\enquote {\bibinfo {title} {New
  approximations and tests of linear fluctuation-response for chaotic nonlinear
  forced-dissipative dynamical systems},}\ }\href {\doibase
  10.1007/s00332-007-9011-9} {\bibfield  {journal} {\bibinfo  {journal} {J.
  Nonlinear Sci.}\ }\textbf {\bibinfo {volume} {18}},\ \bibinfo {pages}
  {303--341} (\bibinfo {year} {2008})}\BibitemShut {NoStop}%
\bibitem [{\citenamefont {Cooper}\ and\ \citenamefont
  {Haynes}(2011)}]{CooperHaynes11}%
  \BibitemOpen
  \bibfield  {author} {\bibinfo {author} {\bibfnamefont {F.~C.}\ \bibnamefont
  {Cooper}}\ and\ \bibinfo {author} {\bibfnamefont {P.~H.}\ \bibnamefont
  {Haynes}},\ }\bibfield  {title} {\enquote {\bibinfo {title} {Climate
  sensitivity via a nonparametric fluctuation-dissipation theorem},}\ }\href
  {http://ezproxy.library.usyd.edu.au/login?url=http://search.proquest.com/docview/868716283?accountid=14757}
  {\bibfield  {journal} {\bibinfo  {journal} {Journal of the Atmospheric
  Sciences}\ }\textbf {\bibinfo {volume} {68}},\ \bibinfo {pages} {937--953}
  (\bibinfo {year} {2011})}\BibitemShut {NoStop}%
\bibitem [{\citenamefont {{Cooper}}, \citenamefont {{Esler}},\ and\
  \citenamefont {{Haynes}}(2013)}]{CooperEtAl13}%
  \BibitemOpen
  \bibfield  {author} {\bibinfo {author} {\bibfnamefont {F.~C.}\ \bibnamefont
  {{Cooper}}}, \bibinfo {author} {\bibfnamefont {J.~G.}\ \bibnamefont
  {{Esler}}}, \ and\ \bibinfo {author} {\bibfnamefont {P.~H.}\ \bibnamefont
  {{Haynes}}},\ }\bibfield  {title} {\enquote {\bibinfo {title} {Estimation of
  the local response to a forcing in a high dimensional system using the
  fluctuation-dissipation theorem},}\ }\href {\doibase 10.5194/npg-20-239-2013}
  {\bibfield  {journal} {\bibinfo  {journal} {Nonlin. Processes Geophys.}\
  }\textbf {\bibinfo {volume} {20}},\ \bibinfo {pages} {239--248} (\bibinfo
  {year} {2013})}\BibitemShut {NoStop}%
\bibitem [{\citenamefont {{Bell}}(1980)}]{Bell80}%
  \BibitemOpen
  \bibfield  {author} {\bibinfo {author} {\bibfnamefont {T.~L.}\ \bibnamefont
  {{Bell}}},\ }\bibfield  {title} {\enquote {\bibinfo {title} {Climate
  sensitivity from fluctuation dissipation: Some simple model tests},}\ }\href
  {\doibase 10.1175/1520-0469(1980)037<1700:CSFFDS>2.0.CO;2} {\bibfield
  {journal} {\bibinfo  {journal} {Journal of the Atmospheric Sciences}\
  }\textbf {\bibinfo {volume} {37}},\ \bibinfo {pages} {1700--1707} (\bibinfo
  {year} {1980})}\BibitemShut {NoStop}%
\bibitem [{\citenamefont {Gritsun}\ and\ \citenamefont
  {Dymnikov}(1999)}]{GritsunDymnikov99}%
  \BibitemOpen
  \bibfield  {author} {\bibinfo {author} {\bibfnamefont {A.}~\bibnamefont
  {Gritsun}}\ and\ \bibinfo {author} {\bibfnamefont {V.}~\bibnamefont
  {Dymnikov}},\ }\bibfield  {title} {\enquote {\bibinfo {title} {Barotropic
  atmosphere response to small external actions: {T}heory and numerical
  experiments},}\ }\href@noop {} {\bibfield  {journal} {\bibinfo  {journal}
  {Izv. Akad. Nauk. Fiz. Atmos. Okeana. Biol.}\ }\textbf {\bibinfo {volume}
  {35}},\ \bibinfo {pages} {565--581} (\bibinfo {year} {1999})}\BibitemShut
  {NoStop}%
\bibitem [{\citenamefont {Abramov}\ and\ \citenamefont
  {Majda}(2009)}]{AbramovMajda09}%
  \BibitemOpen
  \bibfield  {author} {\bibinfo {author} {\bibfnamefont {R.~V.}\ \bibnamefont
  {Abramov}}\ and\ \bibinfo {author} {\bibfnamefont {A.~J.}\ \bibnamefont
  {Majda}},\ }\bibfield  {title} {\enquote {\bibinfo {title} {A new algorithm
  for low-frequency climate response},}\ }\href@noop {} {\bibfield  {journal}
  {\bibinfo  {journal} {Journal of the Atmospheric Sciences}\ }\textbf
  {\bibinfo {volume} {66}},\ \bibinfo {pages} {286--309} (\bibinfo {year}
  {2009})}\BibitemShut {NoStop}%
\bibitem [{\citenamefont {Dymnikov}\ and\ \citenamefont
  {Gritsoun}(2001)}]{DymnikovGritsun01}%
  \BibitemOpen
  \bibfield  {author} {\bibinfo {author} {\bibfnamefont {V.~P.}\ \bibnamefont
  {Dymnikov}}\ and\ \bibinfo {author} {\bibfnamefont {A.~S.}\ \bibnamefont
  {Gritsoun}},\ }\bibfield  {title} {\enquote {\bibinfo {title} {Climate model
  attractors: chaos, quasi-regularity and sensitivity to small perturbations of
  external forcing},}\ }\href {\doibase 10.5194/npg-8-201-2001} {\bibfield
  {journal} {\bibinfo  {journal} {Nonlinear Processes in Geophysics}\ }\textbf
  {\bibinfo {volume} {8}},\ \bibinfo {pages} {201--209} (\bibinfo {year}
  {2001})}\BibitemShut {NoStop}%
\bibitem [{\citenamefont {{North}}, \citenamefont {{Bell}},\ and\ \citenamefont
  {{Hardin}}(1993)}]{NorthEtAl93}%
  \BibitemOpen
  \bibfield  {author} {\bibinfo {author} {\bibfnamefont {G.~R.}\ \bibnamefont
  {{North}}}, \bibinfo {author} {\bibfnamefont {R.~E.}\ \bibnamefont {{Bell}}},
  \ and\ \bibinfo {author} {\bibfnamefont {J.~W.}\ \bibnamefont {{Hardin}}},\
  }\bibfield  {title} {\enquote {\bibinfo {title} {Fluctuation dissipation in a
  general circulation model},}\ }\href {\doibase 10.1007/BF00209665} {\bibfield
   {journal} {\bibinfo  {journal} {Climate Dynamics}\ }\textbf {\bibinfo
  {volume} {8}},\ \bibinfo {pages} {259--264} (\bibinfo {year}
  {1993})}\BibitemShut {NoStop}%
\bibitem [{\citenamefont {Cionni}, \citenamefont {Visconti},\ and\
  \citenamefont {Sassi}(2004)}]{CionniEtAl04}%
  \BibitemOpen
  \bibfield  {author} {\bibinfo {author} {\bibfnamefont {I.}~\bibnamefont
  {Cionni}}, \bibinfo {author} {\bibfnamefont {G.}~\bibnamefont {Visconti}}, \
  and\ \bibinfo {author} {\bibfnamefont {F.}~\bibnamefont {Sassi}},\ }\bibfield
   {title} {\enquote {\bibinfo {title} {Fluctuation dissipation theorem in a
  general circulation model},}\ }\href {\doibase 10.1029/2004GL019739}
  {\bibfield  {journal} {\bibinfo  {journal} {Geophysical Research Letters}\
  }\textbf {\bibinfo {volume} {31}},\ \bibinfo {pages} {L09206} (\bibinfo
  {year} {2004})}\BibitemShut {NoStop}%
\bibitem [{\citenamefont {Gritsun}, \citenamefont {Branstator},\ and\
  \citenamefont {Dymnikov}(2002)}]{GritsunEtAl02}%
  \BibitemOpen
  \bibfield  {author} {\bibinfo {author} {\bibfnamefont {A.}~\bibnamefont
  {Gritsun}}, \bibinfo {author} {\bibfnamefont {G.}~\bibnamefont {Branstator}},
  \ and\ \bibinfo {author} {\bibfnamefont {V.}~\bibnamefont {Dymnikov}},\
  }\bibfield  {title} {\enquote {\bibinfo {title} {Construction of the linear
  response operator of an atmospheric general circulation model to small
  external forcing},}\ }\href@noop {} {\bibfield  {journal} {\bibinfo
  {journal} {Russ. J. Numer. Anal. Math. Modelling}\ }\textbf {\bibinfo
  {volume} {17}},\ \bibinfo {pages} {399--416} (\bibinfo {year}
  {2002})}\BibitemShut {NoStop}%
\bibitem [{\citenamefont {Gritsun}\ and\ \citenamefont
  {Branstator}(2007)}]{GritsunBranstator07}%
  \BibitemOpen
  \bibfield  {author} {\bibinfo {author} {\bibfnamefont {A.}~\bibnamefont
  {Gritsun}}\ and\ \bibinfo {author} {\bibfnamefont {G.}~\bibnamefont
  {Branstator}},\ }\bibfield  {title} {\enquote {\bibinfo {title} {Climate
  response using a three-dimensional operator based on the
  fluctuation-dissipation theorem},}\ }\href@noop {} {\bibfield  {journal}
  {\bibinfo  {journal} {Journal of the Atmospheric Sciences}\ }\textbf
  {\bibinfo {volume} {64}},\ \bibinfo {pages} {2558--2575} (\bibinfo {year}
  {2007})}\BibitemShut {NoStop}%
\bibitem [{\citenamefont {Gritsun}, \citenamefont {Branstator},\ and\
  \citenamefont {Majda}(2008)}]{GritsunEtAl08}%
  \BibitemOpen
  \bibfield  {author} {\bibinfo {author} {\bibfnamefont {A.}~\bibnamefont
  {Gritsun}}, \bibinfo {author} {\bibfnamefont {G.}~\bibnamefont {Branstator}},
  \ and\ \bibinfo {author} {\bibfnamefont {A.}~\bibnamefont {Majda}},\
  }\bibfield  {title} {\enquote {\bibinfo {title} {Climate response of linear
  and quadratic functionals using the fluctuation-dissipation theorem},}\
  }\href@noop {} {\bibfield  {journal} {\bibinfo  {journal} {Journal of the
  Atmospheric Sciences}\ }\textbf {\bibinfo {volume} {65}},\ \bibinfo {pages}
  {2824--2829} (\bibinfo {year} {2008})}\BibitemShut {NoStop}%
\bibitem [{\citenamefont {Ring}\ and\ \citenamefont
  {Plumb}(2008)}]{RingPlumb08}%
  \BibitemOpen
  \bibfield  {author} {\bibinfo {author} {\bibfnamefont {M.~J.}\ \bibnamefont
  {Ring}}\ and\ \bibinfo {author} {\bibfnamefont {R.~A.}\ \bibnamefont
  {Plumb}},\ }\bibfield  {title} {\enquote {\bibinfo {title} {The response of a
  simplified {GCM} to axisymmetric forcings: {A}pplicability of the
  fluctuation--dissipation theorem},}\ }\href@noop {} {\bibfield  {journal}
  {\bibinfo  {journal} {Journal of the Atmospheric Sciences}\ }\textbf
  {\bibinfo {volume} {65}},\ \bibinfo {pages} {3880--3898} (\bibinfo {year}
  {2008})}\BibitemShut {NoStop}%
\bibitem [{\citenamefont {{Gritsun}}(2010)}]{Gritsun10}%
  \BibitemOpen
  \bibfield  {author} {\bibinfo {author} {\bibfnamefont {A.~S.}\ \bibnamefont
  {{Gritsun}}},\ }\bibfield  {title} {\enquote {\bibinfo {title} {Construction
  of response operators to small external forcings for atmospheric general
  circulation models with time periodic right-hand sides},}\ }\href {\doibase
  10.1134/S000143381006006X} {\bibfield  {journal} {\bibinfo  {journal}
  {Izvestiya, Atmospheric and Oceanic Physics}\ }\textbf {\bibinfo {volume}
  {46}},\ \bibinfo {pages} {748--756} (\bibinfo {year} {2010})}\BibitemShut
  {NoStop}%
\bibitem [{\citenamefont {Langen}\ and\ \citenamefont
  {Alexeev}(2005)}]{LangenAlexeev05}%
  \BibitemOpen
  \bibfield  {author} {\bibinfo {author} {\bibfnamefont {P.~L.}\ \bibnamefont
  {Langen}}\ and\ \bibinfo {author} {\bibfnamefont {V.~A.}\ \bibnamefont
  {Alexeev}},\ }\bibfield  {title} {\enquote {\bibinfo {title} {Estimating 2
  $\times$ {$CO_2$} warming in an aquaplanet {GCM} using the
  fluctuation-dissipation theorem},}\ }\href {\doibase 10.1029/2005GL024136}
  {\bibfield  {journal} {\bibinfo  {journal} {Geophysical Research Letters}\
  }\textbf {\bibinfo {volume} {32}} (\bibinfo {year} {2005}),\
  10.1029/2005GL024136},\ \bibinfo {note} {l23708}\BibitemShut {NoStop}%
\bibitem [{\citenamefont {Kirk-Davidoff}(2009)}]{KirkDavidoff09}%
  \BibitemOpen
  \bibfield  {author} {\bibinfo {author} {\bibfnamefont {D.~B.}\ \bibnamefont
  {Kirk-Davidoff}},\ }\bibfield  {title} {\enquote {\bibinfo {title} {On the
  diagnosis of climate sensitivity using observations of fluctuations},}\
  }\href {\doibase 10.5194/acp-9-813-2009} {\bibfield  {journal} {\bibinfo
  {journal} {Atmos. Chem. Phys.}\ }\textbf {\bibinfo {volume} {9}},\ \bibinfo
  {pages} {813--822} (\bibinfo {year} {2009})}\BibitemShut {NoStop}%
\bibitem [{\citenamefont {Fuchs}, \citenamefont {Sherwood},\ and\ \citenamefont
  {Hernandez}(2014)}]{FuchsEtAl14}%
  \BibitemOpen
  \bibfield  {author} {\bibinfo {author} {\bibfnamefont {D.}~\bibnamefont
  {Fuchs}}, \bibinfo {author} {\bibfnamefont {S.}~\bibnamefont {Sherwood}}, \
  and\ \bibinfo {author} {\bibfnamefont {D.}~\bibnamefont {Hernandez}},\
  }\bibfield  {title} {\enquote {\bibinfo {title} {An exploration of
  multivariate fluctuation dissipation operators and their response to sea
  surface temperature perturbations},}\ }\href@noop {} {\bibfield  {journal}
  {\bibinfo  {journal} {Journal of the Atmospheric Sciences}\ }\textbf
  {\bibinfo {volume} {72}},\ \bibinfo {pages} {472--486} (\bibinfo {year}
  {2014})}\BibitemShut {NoStop}%
\bibitem [{\citenamefont {{Ragone}}, \citenamefont {{Lucarini}},\ and\
  \citenamefont {{Lunkeit}}(2016)}]{RagoneEtAl15}%
  \BibitemOpen
  \bibfield  {author} {\bibinfo {author} {\bibfnamefont {F.}~\bibnamefont
  {{Ragone}}}, \bibinfo {author} {\bibfnamefont {V.}~\bibnamefont
  {{Lucarini}}}, \ and\ \bibinfo {author} {\bibfnamefont {F.}~\bibnamefont
  {{Lunkeit}}},\ }\bibfield  {title} {\enquote {\bibinfo {title} {A new
  framework for climate sensitivity and prediction: a modelling perspective},}\
  }\href {\doibase 10.1007/s00382-015-2657-3} {\bibfield  {journal} {\bibinfo
  {journal} {Climate Dynamics}\ }\textbf {\bibinfo {volume} {46}},\ \bibinfo
  {pages} {1459--1471} (\bibinfo {year} {2016})}\BibitemShut {NoStop}%
\bibitem [{\citenamefont {Ruelle}(1997)}]{Ruelle97}%
  \BibitemOpen
  \bibfield  {author} {\bibinfo {author} {\bibfnamefont {D.}~\bibnamefont
  {Ruelle}},\ }\bibfield  {title} {\enquote {\bibinfo {title} {Differentiation
  of {SRB} states},}\ }\href@noop {} {\bibfield  {journal} {\bibinfo  {journal}
  {Communications in Mathematical Physics}\ }\textbf {\bibinfo {volume}
  {187}},\ \bibinfo {pages} {227--241} (\bibinfo {year} {1997})}\BibitemShut
  {NoStop}%
\bibitem [{\citenamefont {Ruelle}(1998)}]{Ruelle98}%
  \BibitemOpen
  \bibfield  {author} {\bibinfo {author} {\bibfnamefont {D.}~\bibnamefont
  {Ruelle}},\ }\bibfield  {title} {\enquote {\bibinfo {title} {General linear
  response formula in statistical mechanics, and the fluctuation-dissipation
  theorem far from equilibrium},}\ }\href {\doibase
  10.1016/S0375-9601(98)00419-8} {\bibfield  {journal} {\bibinfo  {journal}
  {Phys. Lett. A}\ }\textbf {\bibinfo {volume} {245}},\ \bibinfo {pages}
  {220--224} (\bibinfo {year} {1998})}\BibitemShut {NoStop}%
\bibitem [{\citenamefont {Ruelle}(2009{\natexlab{a}})}]{Ruelle09a}%
  \BibitemOpen
  \bibfield  {author} {\bibinfo {author} {\bibfnamefont {D.}~\bibnamefont
  {Ruelle}},\ }\bibfield  {title} {\enquote {\bibinfo {title} {A review of
  linear response theory for general differentiable dynamical systems},}\
  }\href@noop {} {\bibfield  {journal} {\bibinfo  {journal} {Nonlinearity}\
  }\textbf {\bibinfo {volume} {22}},\ \bibinfo {pages} {855--870} (\bibinfo
  {year} {2009}{\natexlab{a}})}\BibitemShut {NoStop}%
\bibitem [{\citenamefont {Ruelle}(2009{\natexlab{b}})}]{Ruelle09b}%
  \BibitemOpen
  \bibfield  {author} {\bibinfo {author} {\bibfnamefont {D.}~\bibnamefont
  {Ruelle}},\ }\bibfield  {title} {\enquote {\bibinfo {title} {Structure and
  f-dependence of the a.c.i.m. for a unimodal map f of {M}isiurewicz type},}\
  }\href {\doibase 10.1007/s00220-008-0637-8} {\bibfield  {journal} {\bibinfo
  {journal} {Communications in Mathematical Physics}\ }\textbf {\bibinfo
  {volume} {287}},\ \bibinfo {pages} {1039--1070} (\bibinfo {year}
  {2009}{\natexlab{b}})}\BibitemShut {NoStop}%
\bibitem [{\citenamefont {Baladi}\ and\ \citenamefont
  {Smania}(2008)}]{BaladiSmania08}%
  \BibitemOpen
  \bibfield  {author} {\bibinfo {author} {\bibfnamefont {V.}~\bibnamefont
  {Baladi}}\ and\ \bibinfo {author} {\bibfnamefont {D.}~\bibnamefont
  {Smania}},\ }\bibfield  {title} {\enquote {\bibinfo {title} {Linear response
  formula for piecewise expanding unimodal maps},}\ }\href@noop {} {\bibfield
  {journal} {\bibinfo  {journal} {Nonlinearity}\ }\textbf {\bibinfo {volume}
  {21}},\ \bibinfo {pages} {677--711} (\bibinfo {year} {2008})}\BibitemShut
  {NoStop}%
\bibitem [{\citenamefont {Baladi}\ and\ \citenamefont
  {Smania}(2010)}]{BaladiSmania10}%
  \BibitemOpen
  \bibfield  {author} {\bibinfo {author} {\bibfnamefont {V.}~\bibnamefont
  {Baladi}}\ and\ \bibinfo {author} {\bibfnamefont {D.}~\bibnamefont
  {Smania}},\ }\bibfield  {title} {\enquote {\bibinfo {title} {Alternative
  proofs of linear response for piecewise expanding unimodal maps},}\
  }\href@noop {} {\bibfield  {journal} {\bibinfo  {journal} {Ergodic Theory and
  Dynamical Systems}\ }\textbf {\bibinfo {volume} {30}},\ \bibinfo {pages}
  {1--20} (\bibinfo {year} {2010})}\BibitemShut {NoStop}%
\bibitem [{\citenamefont {Baladi}(2014)}]{Baladi14}%
  \BibitemOpen
  \bibfield  {author} {\bibinfo {author} {\bibfnamefont {V.}~\bibnamefont
  {Baladi}},\ }\bibfield  {title} {\enquote {\bibinfo {title} {Linear response,
  or else},}\ }in\ \href@noop {} {\emph {\bibinfo {booktitle} {ICM Seoul 2014,
  Proceedings, Volume III}}}\ (\bibinfo {year} {2014})\ pp.\ \bibinfo {pages}
  {525--545},\ \Eprint {http://arxiv.org/abs/1408.2937} {arXiv:1408.2937
  [math.DS]} \BibitemShut {NoStop}%
\bibitem [{\citenamefont {Baladi}, \citenamefont {Benedicks},\ and\
  \citenamefont {Schnellmann}(2015)}]{BaladiEtAl15}%
  \BibitemOpen
  \bibfield  {author} {\bibinfo {author} {\bibfnamefont {V.}~\bibnamefont
  {Baladi}}, \bibinfo {author} {\bibfnamefont {M.}~\bibnamefont {Benedicks}}, \
  and\ \bibinfo {author} {\bibfnamefont {D.}~\bibnamefont {Schnellmann}},\
  }\bibfield  {title} {\enquote {\bibinfo {title} {Whitney-{H}\"older
  continuity of the {SRB} measure for transversal families of smooth unimodal
  maps},}\ }\href {\doibase 10.1007/s00222-014-0554-8} {\bibfield  {journal}
  {\bibinfo  {journal} {Invent. Math.}\ }\textbf {\bibinfo {volume} {201}},\
  \bibinfo {pages} {773--844} (\bibinfo {year} {2015})}\BibitemShut {NoStop}%
\bibitem [{\citenamefont {De~Lima}\ and\ \citenamefont
  {Smania}(2018)}]{DeLimaSmania18}%
  \BibitemOpen
  \bibfield  {author} {\bibinfo {author} {\bibfnamefont {A.}~\bibnamefont
  {De~Lima}}\ and\ \bibinfo {author} {\bibfnamefont {D.}~\bibnamefont
  {Smania}},\ }\bibfield  {title} {\enquote {\bibinfo {title} {Central limit
  theorem for the modulus of continuity of averages of observables on
  transversal families of piecewise expanding unimodal maps},}\ }\href@noop {}
  {\bibfield  {journal} {\bibinfo  {journal} {Journal of the Institute of
  Mathematics of Jussieu}\ }\textbf {\bibinfo {volume} {17}},\ \bibinfo {pages}
  {673--733} (\bibinfo {year} {2018})}\BibitemShut {NoStop}%
\bibitem [{\citenamefont {Gallavotti}\ and\ \citenamefont
  {Cohen}(1995{\natexlab{a}})}]{GallavottiCohen95a}%
  \BibitemOpen
  \bibfield  {author} {\bibinfo {author} {\bibfnamefont {G.}~\bibnamefont
  {Gallavotti}}\ and\ \bibinfo {author} {\bibfnamefont {E.~G.~D.}\ \bibnamefont
  {Cohen}},\ }\bibfield  {title} {\enquote {\bibinfo {title} {Dynamical
  ensembles in nonequilibrium statistical mechanics},}\ }\href {\doibase
  10.1103/PhysRevLett.74.2694} {\bibfield  {journal} {\bibinfo  {journal}
  {Phys. Rev. Lett.}\ }\textbf {\bibinfo {volume} {74}},\ \bibinfo {pages}
  {2694--2697} (\bibinfo {year} {1995}{\natexlab{a}})}\BibitemShut {NoStop}%
\bibitem [{\citenamefont {Gallavotti}\ and\ \citenamefont
  {Cohen}(1995{\natexlab{b}})}]{GallavottiCohen95b}%
  \BibitemOpen
  \bibfield  {author} {\bibinfo {author} {\bibfnamefont {G.}~\bibnamefont
  {Gallavotti}}\ and\ \bibinfo {author} {\bibfnamefont {E.}~\bibnamefont
  {Cohen}},\ }\bibfield  {title} {\enquote {\bibinfo {title} {Dynamical
  ensembles in stationary states},}\ }\href@noop {} {\bibfield  {journal}
  {\bibinfo  {journal} {Journal of Statistical Physics}\ }\textbf {\bibinfo
  {volume} {80}},\ \bibinfo {pages} {931--970} (\bibinfo {year}
  {1995}{\natexlab{b}})}\BibitemShut {NoStop}%
\bibitem [{\citenamefont {Gallavotti}(2019)}]{Gallavotti19}%
  \BibitemOpen
  \bibfield  {author} {\bibinfo {author} {\bibfnamefont {G.}~\bibnamefont
  {Gallavotti}},\ }\bibfield  {title} {\enquote {\bibinfo {title}
  {Nonequilibrium and fluctuation relation},}\ }\href@noop {} {\bibfield
  {journal} {\bibinfo  {journal} {arXiv preprint arXiv:1906.10069}\ } (\bibinfo
  {year} {2019})}\BibitemShut {NoStop}%
\bibitem [{\citenamefont {H\"{a}nggi}(1978)}]{Haenggi78}%
  \BibitemOpen
  \bibfield  {author} {\bibinfo {author} {\bibfnamefont {P.}~\bibnamefont
  {H\"{a}nggi}},\ }\bibfield  {title} {\enquote {\bibinfo {title} {Stochastic
  processes 2: response theory and fluctuation theorems},}\ }\href@noop {}
  {\bibfield  {journal} {\bibinfo  {journal} {Helvetica Physica Acta}\ }\textbf
  {\bibinfo {volume} {51}},\ \bibinfo {pages} {202--219} (\bibinfo {year}
  {1978})}\BibitemShut {NoStop}%
\bibitem [{\citenamefont {Hairer}\ and\ \citenamefont
  {Majda}(2010)}]{HairerMajda10}%
  \BibitemOpen
  \bibfield  {author} {\bibinfo {author} {\bibfnamefont {M.}~\bibnamefont
  {Hairer}}\ and\ \bibinfo {author} {\bibfnamefont {A.~J.}\ \bibnamefont
  {Majda}},\ }\bibfield  {title} {\enquote {\bibinfo {title} {A simple
  framework to justify linear response theory},}\ }\href
  {http://stacks.iop.org/0951-7715/23/i=4/a=008} {\bibfield  {journal}
  {\bibinfo  {journal} {Nonlinearity}\ }\textbf {\bibinfo {volume} {23}},\
  \bibinfo {pages} {909} (\bibinfo {year} {2010})}\BibitemShut {NoStop}%
\bibitem [{\citenamefont {Wormell}\ and\ \citenamefont
  {Gottwald}(2018)}]{WormellGottwald18}%
  \BibitemOpen
  \bibfield  {author} {\bibinfo {author} {\bibfnamefont {C.~L.}\ \bibnamefont
  {Wormell}}\ and\ \bibinfo {author} {\bibfnamefont {G.~A.}\ \bibnamefont
  {Gottwald}},\ }\bibfield  {title} {\enquote {\bibinfo {title} {On the
  validity of linear response theory in high-dimensional deterministic
  dynamical systems},}\ }\href {\doibase 10.1007/s10955-018-2106-x} {\bibfield
  {journal} {\bibinfo  {journal} {Journal of Statistical Physics}\ }\textbf
  {\bibinfo {volume} {172}},\ \bibinfo {pages} {1479--1498} (\bibinfo {year}
  {2018})}\BibitemShut {NoStop}%
\bibitem [{\citenamefont {Kaneko}(1990)}]{Kaneko90}%
  \BibitemOpen
  \bibfield  {author} {\bibinfo {author} {\bibfnamefont {K.}~\bibnamefont
  {Kaneko}},\ }\bibfield  {title} {\enquote {\bibinfo {title} {Globally coupled
  chaos violates the law of large numbers but not the central-limit theorem},}\
  }\href@noop {} {\bibfield  {journal} {\bibinfo  {journal} {Physical review
  letters}\ }\textbf {\bibinfo {volume} {65}},\ \bibinfo {pages} {1391}
  (\bibinfo {year} {1990})}\BibitemShut {NoStop}%
\bibitem [{\citenamefont {Shibata}, \citenamefont {Chawanya},\ and\
  \citenamefont {Kaneko}(1999)}]{Shibata99}%
  \BibitemOpen
  \bibfield  {author} {\bibinfo {author} {\bibfnamefont {T.}~\bibnamefont
  {Shibata}}, \bibinfo {author} {\bibfnamefont {T.}~\bibnamefont {Chawanya}}, \
  and\ \bibinfo {author} {\bibfnamefont {K.}~\bibnamefont {Kaneko}},\
  }\bibfield  {title} {\enquote {\bibinfo {title} {Noiseless collective motion
  out of noisy chaos},}\ }\href@noop {} {\bibfield  {journal} {\bibinfo
  {journal} {Physical review letters}\ }\textbf {\bibinfo {volume} {82}},\
  \bibinfo {pages} {4424} (\bibinfo {year} {1999})}\BibitemShut {NoStop}%
\bibitem [{\citenamefont {Pikovsky}\ and\ \citenamefont
  {Kurths}(1994)}]{Pikovsky94}%
  \BibitemOpen
  \bibfield  {author} {\bibinfo {author} {\bibfnamefont {A.~S.}\ \bibnamefont
  {Pikovsky}}\ and\ \bibinfo {author} {\bibfnamefont {J.}~\bibnamefont
  {Kurths}},\ }\bibfield  {title} {\enquote {\bibinfo {title} {Do globally
  coupled maps really violate the law of large numbers?}}\ }\href {\doibase
  10.1103/PhysRevLett.72.1644} {\bibfield  {journal} {\bibinfo  {journal}
  {Phys. Rev. Lett.}\ }\textbf {\bibinfo {volume} {72}},\ \bibinfo {pages}
  {1644--1646} (\bibinfo {year} {1994})}\BibitemShut {NoStop}%
\bibitem [{\citenamefont {Ershov}\ and\ \citenamefont
  {Potapov}(1995)}]{Ershov95}%
  \BibitemOpen
  \bibfield  {author} {\bibinfo {author} {\bibfnamefont {S.~V.}\ \bibnamefont
  {Ershov}}\ and\ \bibinfo {author} {\bibfnamefont {A.~B.}\ \bibnamefont
  {Potapov}},\ }\bibfield  {title} {\enquote {\bibinfo {title} {On mean field
  fluctuations in globally coupled maps},}\ }\href@noop {} {\bibfield
  {journal} {\bibinfo  {journal} {Physica D: Nonlinear Phenomena}\ }\textbf
  {\bibinfo {volume} {86}},\ \bibinfo {pages} {523--558} (\bibinfo {year}
  {1995})}\BibitemShut {NoStop}%
\bibitem [{\citenamefont {Ershov}\ and\ \citenamefont
  {Potapov}(1997)}]{Ershov97}%
  \BibitemOpen
  \bibfield  {author} {\bibinfo {author} {\bibfnamefont {S.~V.}\ \bibnamefont
  {Ershov}}\ and\ \bibinfo {author} {\bibfnamefont {A.~B.}\ \bibnamefont
  {Potapov}},\ }\bibfield  {title} {\enquote {\bibinfo {title} {On mean field
  fluctuations in globally coupled logistic-type maps},}\ }\href@noop {}
  {\bibfield  {journal} {\bibinfo  {journal} {Physica D: Nonlinear Phenomena}\
  }\textbf {\bibinfo {volume} {106}},\ \bibinfo {pages} {9--38} (\bibinfo
  {year} {1997})}\BibitemShut {NoStop}%
\bibitem [{\citenamefont {S\'elley}\ and\ \citenamefont
  {B\'alint}(2016)}]{Selley16}%
  \BibitemOpen
  \bibfield  {author} {\bibinfo {author} {\bibfnamefont {F.}~\bibnamefont
  {S\'elley}}\ and\ \bibinfo {author} {\bibfnamefont {P.}~\bibnamefont
  {B\'alint}},\ }\bibfield  {title} {\enquote {\bibinfo {title} {Mean-field
  coupling of identical expanding circle maps},}\ }\href {\doibase
  10.1007/s10955-016-1568-y} {\bibfield  {journal} {\bibinfo  {journal}
  {Journal of Statistical Physics}\ }\textbf {\bibinfo {volume} {164}},\
  \bibinfo {pages} {858--889} (\bibinfo {year} {2016})}\BibitemShut {NoStop}%
\bibitem [{\citenamefont {Gottwald}, \citenamefont {Wormell},\ and\
  \citenamefont {Wouters}(2016)}]{GottwaldEtAl16}%
  \BibitemOpen
  \bibfield  {author} {\bibinfo {author} {\bibfnamefont {G.~A.}\ \bibnamefont
  {Gottwald}}, \bibinfo {author} {\bibfnamefont {J.~P.}\ \bibnamefont
  {Wormell}}, \ and\ \bibinfo {author} {\bibfnamefont {J.}~\bibnamefont
  {Wouters}},\ }\bibfield  {title} {\enquote {\bibinfo {title} {On spurious
  detection of linear response and misuse of the fluctuation-dissipation
  theorem in finite time series},}\ }\href
  {https://doi.org/10.1016/j.physd.2016.05.010} {\bibfield  {journal} {\bibinfo
   {journal} {Phys. D}\ }\textbf {\bibinfo {volume} {331}},\ \bibinfo {pages}
  {89--101} (\bibinfo {year} {2016})}\BibitemShut {NoStop}%
\bibitem [{\citenamefont {Lyubich}(2002)}]{Lyubich02}%
  \BibitemOpen
  \bibfield  {author} {\bibinfo {author} {\bibfnamefont {M.}~\bibnamefont
  {Lyubich}},\ }\bibfield  {title} {\enquote {\bibinfo {title} {Almost every
  real quadratic map is either regular or stochastic},}\ }\href {\doibase
  10.2307/3597183} {\bibfield  {journal} {\bibinfo  {journal} {Ann. of Math.
  (2)}\ }\textbf {\bibinfo {volume} {156}},\ \bibinfo {pages} {1--78} (\bibinfo
  {year} {2002})}\BibitemShut {NoStop}%
\bibitem [{\citenamefont {Collet}\ and\ \citenamefont
  {Eckmann}(1983)}]{colletEckmann83}%
  \BibitemOpen
  \bibfield  {author} {\bibinfo {author} {\bibfnamefont {P.}~\bibnamefont
  {Collet}}\ and\ \bibinfo {author} {\bibfnamefont {J.-P.}\ \bibnamefont
  {Eckmann}},\ }\bibfield  {title} {\enquote {\bibinfo {title} {Positive
  {L}iapunov exponents and absolute continuity for maps of the interval},}\
  }\href {\doibase 10.1017/S0143385700001802} {\bibfield  {journal} {\bibinfo
  {journal} {Ergodic Theory Dynam. Systems}\ }\textbf {\bibinfo {volume} {3}},\
  \bibinfo {pages} {13--46} (\bibinfo {year} {1983})}\BibitemShut {NoStop}%
\bibitem [{\citenamefont {Alves}, \citenamefont {Luzzatto},\ and\ \citenamefont
  {Pinheiro}(2004)}]{AlvesEtAl04}%
  \BibitemOpen
  \bibfield  {author} {\bibinfo {author} {\bibfnamefont {J.~F.}\ \bibnamefont
  {Alves}}, \bibinfo {author} {\bibfnamefont {S.}~\bibnamefont {Luzzatto}}, \
  and\ \bibinfo {author} {\bibfnamefont {V.}~\bibnamefont {Pinheiro}},\
  }\bibfield  {title} {\enquote {\bibinfo {title} {{Lyapunov exponents and
  rates of mixing for one-dimensional maps}},}\ }\href@noop {} {\bibfield
  {journal} {\bibinfo  {journal} {Ergodic Theory Dynam. Systems}\ }\textbf
  {\bibinfo {volume} {24}},\ \bibinfo {pages} {637--657} (\bibinfo {year}
  {2004})}\BibitemShut {NoStop}%
\bibitem [{\citenamefont {Melbourne}\ and\ \citenamefont
  {Nicol}(2008)}]{MelbourneNicol08}%
  \BibitemOpen
  \bibfield  {author} {\bibinfo {author} {\bibfnamefont {I.}~\bibnamefont
  {Melbourne}}\ and\ \bibinfo {author} {\bibfnamefont {M.}~\bibnamefont
  {Nicol}},\ }\bibfield  {title} {\enquote {\bibinfo {title} {Large deviations
  for nonuniformly hyperbolic systems},}\ }\href {\doibase
  10.1090/S0002-9947-08-04520-0} {\bibfield  {journal} {\bibinfo  {journal}
  {Trans. Amer. Math. Soc.}\ }\textbf {\bibinfo {volume} {360}},\ \bibinfo
  {pages} {6661--6676} (\bibinfo {year} {2008})}\BibitemShut {NoStop}%
\bibitem [{\citenamefont {Trefethen}(2013)}]{Trefethen13}%
  \BibitemOpen
  \bibfield  {author} {\bibinfo {author} {\bibfnamefont {L.~N.}\ \bibnamefont
  {Trefethen}},\ }\href@noop {} {\emph {\bibinfo {title} {Approximation theory
  and approximation practice}}}\ (\bibinfo  {publisher} {Siam},\ \bibinfo
  {address} {Philadelphia, PA},\ \bibinfo {year} {2013})\BibitemShut {NoStop}%
\bibitem [{\citenamefont {Rice}(2006)}]{Rice}%
  \BibitemOpen
  \bibfield  {author} {\bibinfo {author} {\bibfnamefont {J.}~\bibnamefont
  {Rice}},\ }\href@noop {} {\emph {\bibinfo {title} {Mathematical statistics
  and data analysis}}}\ (\bibinfo  {publisher} {Thomson Learning},\ \bibinfo
  {address} {Belmont, CA},\ \bibinfo {year} {2006})\BibitemShut {NoStop}%
\bibitem [{\citenamefont {Ruelle}(2018)}]{Ruelle18}%
  \BibitemOpen
  \bibfield  {author} {\bibinfo {author} {\bibfnamefont {D.}~\bibnamefont
  {Ruelle}},\ }\bibfield  {title} {\enquote {\bibinfo {title} {Linear response
  theory for diffeomorphisms with tangencies of stable and unstable
  manifolds---a contribution to the {G}allavotti-{C}ohen chaotic hypothesis},}\
  }\href@noop {} {\bibfield  {journal} {\bibinfo  {journal} {Nonlinearity}\
  }\textbf {\bibinfo {volume} {31}},\ \bibinfo {pages} {5683} (\bibinfo {year}
  {2018})}\BibitemShut {NoStop}%
\bibitem [{\citenamefont {Avila}, \citenamefont {Lyubich},\ and\ \citenamefont
  {de~Melo}(2003)}]{AvilaEtAl03}%
  \BibitemOpen
  \bibfield  {author} {\bibinfo {author} {\bibfnamefont {A.}~\bibnamefont
  {Avila}}, \bibinfo {author} {\bibfnamefont {M.}~\bibnamefont {Lyubich}}, \
  and\ \bibinfo {author} {\bibfnamefont {W.}~\bibnamefont {de~Melo}},\
  }\bibfield  {title} {\enquote {\bibinfo {title} {Regular or stochastic
  dynamics in real analytic families of unimodal maps},}\ }\href {\doibase
  10.1007/s00222-003-0307-6} {\bibfield  {journal} {\bibinfo  {journal}
  {Inventiones mathematicae}\ }\textbf {\bibinfo {volume} {154}},\ \bibinfo
  {pages} {451--550} (\bibinfo {year} {2003})}\BibitemShut {NoStop}%
\bibitem [{\citenamefont {Baladi}\ and\ \citenamefont
  {Smania}(2012)}]{BaladiSmania12}%
  \BibitemOpen
  \bibfield  {author} {\bibinfo {author} {\bibfnamefont {V.}~\bibnamefont
  {Baladi}}\ and\ \bibinfo {author} {\bibfnamefont {D.}~\bibnamefont
  {Smania}},\ }\bibfield  {title} {\enquote {\bibinfo {title} {Linear response
  for smooth deformations of generic nonuniformly hyperbolic unimodal maps},}\
  }in\ \href@noop {} {\emph {\bibinfo {booktitle} {Annales scientifiques de
  l'{\'E}cole Normale Sup{\'e}rieure}}},\ Vol.~\bibinfo {volume} {45}\
  (\bibinfo {year} {2012})\ pp.\ \bibinfo {pages} {861--926}\BibitemShut
  {NoStop}%
\bibitem [{\citenamefont {Ruelle}(2004)}]{Ruelle04}%
  \BibitemOpen
  \bibfield  {author} {\bibinfo {author} {\bibfnamefont {D.}~\bibnamefont
  {Ruelle}},\ }\href {\doibase 10.1017/CBO9780511617546} {\emph {\bibinfo
  {title} {Thermodynamic Formalism: {T}he Mathematical Structure of Equilibrium
  Statistical Mechanics}}},\ \bibinfo {edition} {2nd}\ ed.,\ Cambridge
  Mathematical Library\ (\bibinfo  {publisher} {Cambridge University Press},\
  \bibinfo {year} {2004})\BibitemShut {NoStop}%
\bibitem [{\citenamefont {Froyland}, \citenamefont {Lloyd},\ and\ \citenamefont
  {Quas}(2013)}]{Froyland10}%
  \BibitemOpen
  \bibfield  {author} {\bibinfo {author} {\bibfnamefont {G.}~\bibnamefont
  {Froyland}}, \bibinfo {author} {\bibfnamefont {S.}~\bibnamefont {Lloyd}}, \
  and\ \bibinfo {author} {\bibfnamefont {A.}~\bibnamefont {Quas}},\ }\bibfield
  {title} {\enquote {\bibinfo {title} {A semi-invertible {O}seledets theorem
  with applications to transfer operator cocycles},}\ }\href@noop {} {\bibfield
   {journal} {\bibinfo  {journal} {Discrete Contin. Dyn. Syst.}\ }\textbf
  {\bibinfo {volume} {33}},\ \bibinfo {pages} {3835--3860} (\bibinfo {year}
  {2013})}\BibitemShut {NoStop}%
\bibitem [{\citenamefont {Buzzi}(1999)}]{Buzzi99}%
  \BibitemOpen
  \bibfield  {author} {\bibinfo {author} {\bibfnamefont {J.}~\bibnamefont
  {Buzzi}},\ }\bibfield  {title} {\enquote {\bibinfo {title} {Exponential decay
  of correlations for random {L}asota--{Y}orke maps},}\ }\href@noop {}
  {\bibfield  {journal} {\bibinfo  {journal} {Communications in mathematical
  physics}\ }\textbf {\bibinfo {volume} {208}},\ \bibinfo {pages} {25--54}
  (\bibinfo {year} {1999})}\BibitemShut {NoStop}%
\bibitem [{\citenamefont {Wormell}(2019)}]{Wormell19}%
  \BibitemOpen
  \bibfield  {author} {\bibinfo {author} {\bibfnamefont {C.~L.}\ \bibnamefont
  {Wormell}},\ }\bibfield  {title} {\enquote {\bibinfo {title} {Spectral
  {G}alerkin methods for transfer operators in uniformly expanding dynamics},}\
  }\href@noop {} {\bibfield  {journal} {\bibinfo  {journal} {Numerische
  Mathematik}\ }\textbf {\bibinfo {volume} {142}},\ \bibinfo {pages} {421--463}
  (\bibinfo {year} {2019})}\BibitemShut {NoStop}%
\bibitem [{\citenamefont {Wormell}()}]{Poltergeist}%
  \BibitemOpen
  \bibfield  {author} {\bibinfo {author} {\bibfnamefont {C.~L.}\ \bibnamefont
  {Wormell}},\ }\href@noop {} {\enquote {\bibinfo {title} {{{P}oltergeist}},}\
  }\bibinfo {note} {Available at \url{http://github.com/wormell/Poltergeist.jl}
  and in the Julia package repository}\BibitemShut {NoStop}%
\bibitem [{Note1()}]{Note1}%
  \BibitemOpen
  \bibinfo {note} {When randomly searching for equilibria it is important to
  make sure that, as well as randomly initialising the $q^{(j)}_0$, the
  {\protect \it distribution} from which the $q^{(j)}_0$ are sampled is also
  randomly initialised, as up to an error term of $\protect \mathcal
  {O}(M^{-1/2})$ the macroscopic dynamics are deterministic functions of the
  initial measures of the microscopic variables $\mu ^a_0$.}\BibitemShut
  {Stop}%
\bibitem [{\citenamefont {Collet}\ and\ \citenamefont
  {Eckmann}(2007)}]{ColletEckmann07}%
  \BibitemOpen
  \bibfield  {author} {\bibinfo {author} {\bibfnamefont {P.}~\bibnamefont
  {Collet}}\ and\ \bibinfo {author} {\bibfnamefont {J.-P.}\ \bibnamefont
  {Eckmann}},\ }\href@noop {} {\emph {\bibinfo {title} {Concepts and results in
  chaotic dynamics: a short course}}}\ (\bibinfo  {publisher} {Springer Science
  \& Business Media},\ \bibinfo {address} {Berlin},\ \bibinfo {year}
  {2007})\BibitemShut {NoStop}%
\bibitem [{\citenamefont {Gottwald}\ and\ \citenamefont
  {Melbourne}(2014)}]{GottwaldMelbourne14}%
  \BibitemOpen
  \bibfield  {author} {\bibinfo {author} {\bibfnamefont {G.~A.}\ \bibnamefont
  {Gottwald}}\ and\ \bibinfo {author} {\bibfnamefont {I.}~\bibnamefont
  {Melbourne}},\ }\bibfield  {title} {\enquote {\bibinfo {title} {A test for a
  conjecture on the nature of attractors for smooth dynamical systems},}\
  }\href {\doibase 10.1063/1.4868984} {\bibfield  {journal} {\bibinfo
  {journal} {Chaos: An Interdisciplinary Journal of Nonlinear Science}\
  }\textbf {\bibinfo {volume} {24}},\ \bibinfo {pages} {024403} (\bibinfo
  {year} {2014})}\BibitemShut {NoStop}%
\bibitem [{\citenamefont {Abramov}(2010)}]{Abramov10}%
  \BibitemOpen
  \bibfield  {author} {\bibinfo {author} {\bibfnamefont {R.~V.}\ \bibnamefont
  {Abramov}},\ }\bibfield  {title} {\enquote {\bibinfo {title} {Approximate
  linear response for slow variables of dynamics with explicit time scale
  separation},}\ }\href@noop {} {\bibfield  {journal} {\bibinfo  {journal}
  {Journal of Computational Physics}\ }\textbf {\bibinfo {volume} {229}},\
  \bibinfo {pages} {7739--7746} (\bibinfo {year} {2010})}\BibitemShut {NoStop}%
\bibitem [{\citenamefont {Pollicott}(1986)}]{Pollicott86}%
  \BibitemOpen
  \bibfield  {author} {\bibinfo {author} {\bibfnamefont {M.}~\bibnamefont
  {Pollicott}},\ }\bibfield  {title} {\enquote {\bibinfo {title} {Meromorphic
  extensions of generalised zeta functions},}\ }\href {\doibase
  10.1007/BF01388795} {\bibfield  {journal} {\bibinfo  {journal} {Invent.
  Math.}\ }\textbf {\bibinfo {volume} {85}},\ \bibinfo {pages} {147--164}
  (\bibinfo {year} {1986})}\BibitemShut {NoStop}%
\bibitem [{\citenamefont {Cvitanovic}\ and\ \citenamefont
  {Eckhardt}(1991)}]{CvitanovicEckhardt91}%
  \BibitemOpen
  \bibfield  {author} {\bibinfo {author} {\bibfnamefont {P.}~\bibnamefont
  {Cvitanovic}}\ and\ \bibinfo {author} {\bibfnamefont {B.}~\bibnamefont
  {Eckhardt}},\ }\bibfield  {title} {\enquote {\bibinfo {title} {Periodic orbit
  expansions for classical smooth flows},}\ }\href {\doibase
  10.1088/0305-4470/24/5/005} {\bibfield  {journal} {\bibinfo  {journal}
  {Journal of Physics A: Mathematical and General}\ }\textbf {\bibinfo {volume}
  {24}},\ \bibinfo {pages} {L237--L241} (\bibinfo {year} {1991})}\BibitemShut
  {NoStop}%
\bibitem [{\citenamefont {Eckhardt}\ and\ \citenamefont
  {Grossmann}(1994)}]{EckhardtGrossmann94}%
  \BibitemOpen
  \bibfield  {author} {\bibinfo {author} {\bibfnamefont {B.}~\bibnamefont
  {Eckhardt}}\ and\ \bibinfo {author} {\bibfnamefont {S.}~\bibnamefont
  {Grossmann}},\ }\bibfield  {title} {\enquote {\bibinfo {title} {Correlation
  functions in chaotic systems from periodic orbits},}\ }\href {\doibase
  10.1103/PhysRevE.50.4571} {\bibfield  {journal} {\bibinfo  {journal} {Phys.
  Rev. E}\ }\textbf {\bibinfo {volume} {50}},\ \bibinfo {pages} {4571--4576}
  (\bibinfo {year} {1994})}\BibitemShut {NoStop}%
\bibitem [{\citenamefont {Gottwald}\ and\ \citenamefont
  {Melbourne}(2013)}]{GottwaldMelbourne13c}%
  \BibitemOpen
  \bibfield  {author} {\bibinfo {author} {\bibfnamefont {G.~A.}\ \bibnamefont
  {Gottwald}}\ and\ \bibinfo {author} {\bibfnamefont {I.}~\bibnamefont
  {Melbourne}},\ }\bibfield  {title} {\enquote {\bibinfo {title}
  {Homogenization for deterministic maps and multiplicative noise},}\
  }\href@noop {} {\bibfield  {journal} {\bibinfo  {journal} {Proceedings of the
  Royal Society A: Mathematical, Physical and Engineering Science}\ }\textbf
  {\bibinfo {volume} {469}} (\bibinfo {year} {2013})}\BibitemShut {NoStop}%
\end{thebibliography}%

\appendix

\section{Statistical test for linear response given time series}
\label{a.lrttest}

A statistical test to probe the linear or higher-order response of a chaotic system from time series data at various parameter values was proposed in \cite{GottwaldEtAl16}. In this appendix we summarise the principles of this test, which we use in the body of this paper.

Suppose that a system has response $\E^\eps \Psi$ which has certain regularity properties for $\eps \in [\eps_1,\eps_2]$ (or around some $\eps_0$), and suppose this regularity property means that there exist functions $\varphi_i(\eps), i = 1,\ldots, I$ and (unknown) coefficients $\beta_i \in \R$ such that 
\begin{equation} \E^\eps \Psi \approx \sum_{i=1}^I \beta_i \varphi_i. \label{e.LRapprox}\end{equation}
is a good approximation in $L^2([\eps_1,\eps_2])$.

For example, if $\E^\eps \Psi$ is $C^1$ and $\eps_2 -\eps_1$ is sufficiently small, then $\E^\eps \Psi$ can be well-approximated with a Taylor expansion about $\eps_1$: thus, $\varphi_0 \equiv 1$ and $\varphi_1(\eps) = \eps - \eps_1$ form a good basis for approximation, and we would expect an $L^2$ error of size $o(\eps_1-\eps_2)$. If instead $\E^\eps \Psi$ is smooth (e.g. $C^r$) on a larger interval, then we could choose Chebyshev polynomials as a basis for approximation $\varphi_i(\eps) = T_{i-1}((2\eps - \eps_1 - \eps_2)/(\eps_2 - \eps_1))$ for $i = 1,\ldots, I$, with an $L^2$ error of $O(I^r)$.

Suppose that for perturbation values $\eps_j,\, j = 1,\ldots,J$ we have time series of the observable's dynamics $(\Psi_{j,n})_{n=1,\ldots,N}$ where the time series length $N$ is sufficiently large. Supposing that (as is typical for many systems \cite{GottwaldMelbourne13c}) $\Psi_n$ obeys a CLT for each selected parameter, then for large enough $N$ the Birkhoff averages for each $\eps_j$ have Gaussian approximations
\begin{equation} \bar{\Psi}^N_j := \frac{1}{N} \sum_{j=1}^N \Psi_{j,n} = \E^{\eps_j} \Psi + \sigma(\eps_j) \xi_j / \sqrt{N}, \label{e.MeanCLT}\end{equation}
where $\xi_j$ are {\it i.i.d.} standard normal variables and the Birkhoff variance $\sigma^2(\eps_j)$ can be estimated by various means, including taking multiple time series for each $\eps_j$, or subsampling.

If we define the vector with coefficients $y_j = \bar{\Psi}^N_j$ and the matrix with coefficients $X_{ji} = \sqrt{N} \sigma(\eps_j)^{-1} \varphi_i(\eps_j)$, then we can write (\ref{e.LRapprox}-\ref{e.MeanCLT}) as the linear equation
\[ y = X \beta + \xi, \]
where $\xi \sim \mathcal{N}(0,I_{J\times J}).$ This is of course just a standard linear statistical model, and we can use the theory of these models\cite{Rice} to test the null hypothesis that the approximation of the response by the $\phi_i$ (\ref{e.LRapprox}) is an equality, i.e. that $\E^\eps \Psi$ has linear (or smooth) response. 

Defining the least-squares projection matrix
\[ H = X (X^T X)^{-1} X^T \]
and the Pearson chi-square test statistic
\[ \chi^2 = y^T (I - H) y, \]
we have that if the approximate equality in (\ref{e.LRapprox}) is exact, then $\chi^2$ has chi-squared distribution $\chi^2_{J-I}$ where $I$ is the number of basis functions $\varphi_i$.

If $\chi^2_{\textrm{obs}}$ is the observed value of the test statistic, the p-value for the test for linear (or higher-order) response is then given by 
\[ p = P(\chi^2_{J-I} \geq \chi^2_{\textrm{obs}}), \]
provided the error associated with the non-exact nature of the approximation (\ref{e.LRapprox}) is appropriately small. This error is small if
\[ \E \chi^2 - \E \chi^2_{J-I} = N \|(I-H) (\E^{\eps_j} \Psi / \sigma(\eps_j))_{j = 1,\ldots,J}\|_{\ell^2}, \]
which, supposing $\sigma$ is a reasonably smooth function of $\eps$ and the $\eps_j$ are uniformly spaced, estimates the minimum possible $L^2(\sigma^2)$ error in approximations of the response of the form in (\ref{e.LRapprox}), multiplied by the sample sizes $N$.

\section{Numerical method to compute the thermodynamic limit $M\to\infty$ for uniformly expanding maps}
\label{a.poltergeist}

In the thermodynamic limit of infinite $M$ the strong law of large numbers holds and 
\[ \Phi_n = \langle \E\Phi_n \rangle = \int \phi(q) d\mu_n(q), \]
where $\mu_n$ is the (time-varying) physical measure of the system, which evolves as
\[ \mu_{n+1} = \mathcal{L}_{K_n} \mu_n, \]
where $\mathcal{L}_{K_n}$ is the transfer operator of the system (\ref{e.Lanford}) and recalling that $K_n = \tanh(\eps \Phi_n - 2)$.

Because for all fixed $K_n=K$ the map (\ref{e.Lanford}) is uniformly expanding, the physical measures $\mu_n$ are absolutely continuous with respect to Lebesgue, and we can write them as $\mu_n(q) dq$. Furthermore, because the map (\ref{e.Lanford}) is analytic and hence infinitely many times differentiable, it is possible to approximate the measure density and transfer operator dynamics very accurately using Chebyshev spectral Galerkin methods \cite{Wormell19}. We have implemented an adaptive-order spectral approximation of the measure density in the Julia package Poltergeist.jl \cite{Poltergeist}, which allows us to simulate the dynamics of $\mu_n$. The core routine, which outputs $\mu_{n+1}$ and $\Phi_{n+1}$ given inputs $\mu_n$, $\eps$ and the driving $d_n$ (by default $\Phi_n$), is defined as follows (note that Julia recognises Unicode characters):
\begin{lstlisting}
function F($\mu$_n, $\eps$, d_n = sum($\phi$*$\mu$_n))
    K_n = tanh(d_n*$\eps$-2)
    f_n = f_map(K_n)
        # create MarkovMap object
    $\mu$_n1 = transfer(f_n,$\mu$_n) 
        # compute $\mu_{n+1}$
    return $\mu$_n1, sum($\phi$*$\mu$_n1)
end
\end{lstlisting}
More details of the algorithm and some examples of its use may be found at \url{https://github.com/wormell/PoltergeistExamples/blob/master/WeakSelfCoupling-LimitingSystem.ipynb}. 

After the first time the routine is called (during which Julia compiles the code), the algorithm takes around $8 \times 10^{-4}$ seconds on a standard laptop to compute each $\mu_{n+1}$ from $\mu_n$, and has an approximation error of only around $10^{-13}$: by comparison, if one aims to estimate $\mu_n$ as a Monte-Carlo approximation with a large ensemble of $M = \O(10^{8})$, a relatively large approximation error of $10^{-3}$ is incurred.

%In \cite{Wormell19} an adaptive spectral method was proposed to 

\end{document}